\documentclass{Article2022}

\newcommand{\VarX}{\mathsf{x}}
\newcommand{\VarZ}{\mathsf{z}}

\newcommand{\TreeT}{\mathfrak{t}}
\newcommand{\TreeS}{\mathfrak{s}}
\newcommand{\TreeR}{\mathfrak{r}}

\newcommand{\ForestF}{\mathfrak{f}}
\newcommand{\ForestG}{\mathfrak{g}}

\newcommand{\Length}{\ell}
\newcommand{\SetVariables}{\mathbb{X}}
\newcommand{\GeneratingSet}{\mathfrak{G}}
\newcommand{\Product}{\operatorname{\bullet}}
\newcommand{\App}{\operatorname{\star}}
\newcommand{\Combinator}[1]{\operatorname{\mathbf{#1}}}
\newcommand{\ATRS}{\mathcal{T}}
\newcommand{\CLS}{\mathcal{C}}
\newcommand{\Monoid}{\mathcal{M}}

\newcommand{\Rew}{\to}
\newcommand{\RewContext}{\Rightarrow}
\newcommand{\Leq}{\preccurlyeq}
\newcommand{\Equiv}{\equiv}
\newcommand{\SetTerms}{\mathfrak{T}}
\newcommand{\Deg}{\mathrm{deg}}
\newcommand{\Height}{\mathrm{ht}}
\newcommand{\RewGraph}{\mathrm{G}}
\newcommand{\SetDuplicativeTrees}{\mathcal{D}}
\newcommand{\SetDuplicativeForests}{\SetDuplicativeTrees^*}
\newcommand{\LeqDuplicative}{\ll}
\newcommand{\Meet}{\wedge}
\newcommand{\JJoin}{\vee}
\newcommand{\MockingbirdLattice}{\mathrm{M}}
\newcommand{\Poset}{\mathcal{P}}
\newcommand{\Model}{\mathcal{M}}
\newcommand{\ModelMin}{\Model^{\mathrm{min}}}
\newcommand{\ModelMax}{\Model^{\mathrm{max}}}
\newcommand{\Ladder}{\mathfrak{l}}
\newcommand{\RightComb}{\mathfrak{r}}
\newcommand{\ToForest}{\texttt{fr}}
\newcommand{\Pruned}{\texttt{pr}}
\newcommand{\MaxLength}{\texttt{ml}}
\newcommand{\T}[1]{{\operatorname{\mathrm T}^{#1}}}
\newcommand{\TT}[2]{\operatorname{{#1}^{\otimes_{#2}}}}
\newcommand{\M}{{\Combinator{M}}}
\newcommand{\Roots}{\mathrm{rts}}
\newcommand{\Support}{\mathrm{Supp}}
\newcommand{\Isolated}{\mathrm{iso}}
\newcommand{\GreaterLadders}{\mathcal{L}}

\newcommand{\SequenceMinDeg}{\mathbf{d}_{\min}}
\newcommand{\SequenceMaxDeg}{\mathbf{d}_{\max}}
\newcommand{\SequenceMinMaxHt}{\mathbf{h}_{\mathrm{mm}}}
\newcommand{\SequenceIsolatedDeg}{\mathbf{d}_{\Isolated}}
\newcommand{\SequenceIsolatedHt}{\mathbf{h}_{\Isolated}}
\newcommand{\SequenceGreater}{\mathbf{h}_{\mathrm{gr}}}
\newcommand{\SequenceNbInputs}{\mathbf{h}_{\mathrm{ni}}}
\newcommand{\SequenceNbSmaller}{\mathbf{h}_{\mathrm{ns}}}

\newcommand{\SeriesF}{\mathbf{F}}
\newcommand{\SeriesMin}{\mathbf{F}_{\min}}
\newcommand{\SeriesMax}{\mathbf{F}_{\max}}
\newcommand{\SeriesIsolated}{\mathbf{F}_{\Isolated}}
\newcommand{\CharacteristicSeries}{\mathbf{c}}
\newcommand{\SeriesLadders}{\mathbf{ld}}

\newcommand{\GeneratingSeriesMinDeg}{\mathsf{D}_{\min}}
\newcommand{\GeneratingSeriesMaxDeg}{\mathsf{D}_{\max}}
\newcommand{\GeneratingSeriesIsolatedDeg}{\mathsf{D}_{\Isolated}}
\newcommand{\GeneratingSeriesMinMaxHt}{\mathsf{H}_{\mathrm{mm}}}
\newcommand{\GeneratingSeriesIsolatedHt}{\mathsf{H}_{\Isolated}}
\newcommand{\GeneratingSeriesGreaterHt}{\mathsf{H}_{\mathrm{gr}}}
\newcommand{\GeneratingSeriesNbInputsHt}{\mathsf{H}_{\mathrm{ni}}}
\newcommand{\GeneratingSeriesNbSmallerHt}{\mathsf{H}_{\mathrm{ns}}}

\newcommand{\EnumerationMap}{\mathrm{en}}
\newcommand{\SeriesGreater}{\mathbf{gr}}
\newcommand{\SeriesNbSmaller}{\mathbf{ns}}
\newcommand{\SeriesMeetDecomposition}{\mathbf{md}}
\newcommand{\SeriesCoverings}{\mathbf{cv}}
\newcommand{\SeriesNbInputs}{\mathbf{ni}}

\newcommand{\HadamardProduct}{\operatorname{\boxtimes}}
\newcommand{\MaxProduct}{\operatorname{\uparrow}}
\newcommand{\ConcatenateForests}{\operatorname{\centerdot}}
\newcommand{\Merge}{\mathrm{mg}}

\newcommand{\RewDuplicative}{%
    \!%
    \begin{tikzpicture}[Centering,xscale=.08]
        \node at(0,0){$\RewContext$};
        \node at(1,0){$\RewContext$};
    \end{tikzpicture}%
    \!}

\newcommand{\WhiteNode}{
    \begin{tikzpicture}[Centering]
        \node[Node,WhiteNode,minimum size=1.5mm]at(0,0){};
    \end{tikzpicture}}

\newcommand{\BlackNode}{
    \begin{tikzpicture}[Centering]
        \node[Node,BlackNode,minimum size=1.5mm]at(0,0){};
    \end{tikzpicture}}

\newcommand{\AnyNode}{
    \begin{tikzpicture}[Centering]
        \node[Node,AnyNode,minimum size=1.5mm]at(0,0){};
    \end{tikzpicture}}

\tikzstyle{WhiteNode}=[Node,draw=ColB!80,fill=ColB!4]
\tikzstyle{BlackNode}=[Node,draw=ColA!90,fill=ColA!60]
\tikzstyle{AnyNode}=[Node,rectangle,draw=ColA!90,fill=ColB!20]

\title[The Mockingbird lattice]{The combinator $\M$ and the Mockingbird lattice}
\date{\today}
\author[S. Giraudo]{Samuele Giraudo}
\address{%
    \scriptsize
    LIGM, Université Gustave Eiffel, CNRS, ESIEE Paris, F-$77454$
    Marne-la-Vallée, France.
    {\tt \href{mailto:samuele.giraudo@univ-eiffel.fr}{samuele.giraudo@univ-eiffel.fr}}}
\thanks{%
    This research has been partially supported by the projects CARPLO (ANR-20-CE40-0007)
    and LambdaComb (ANR-21-CE48-0017) of the Agence nationale de la recherche.}

\begin{document}

\begin{abstract}
    We study combinatorial and order theoretic structures arising from the fragment of
    combinatory logic spanned by the basic combinator $\M$. This basic combinator, named as
    the Mockingbird by Smullyan, is defined by the rewrite rule $\M \VarX_1 \Rew \VarX_1
    \VarX_1$. We prove that the reflexive and transitive closure of this rewrite relation is
    a partial order on terms on $\M$ and that all connected components of its rewrite graph
    are Hasse diagram of lattices. This last result is based on the introduction of new
    lattices on duplicative forests, which are sorts of treelike structures. These lattices
    are not graded, not self-dual, and not semidistributive. We present some enumerative
    properties of these lattices like the enumeration of their elements, of the edges of
    their Hasse diagrams, and of their intervals. These results are derived from formal
    power series on terms and on duplicative forests endowed with particular operations.
\vspace{-3ex}
\end{abstract}

\MakeFirstPage

\section*{Introduction}
Combinatory logic is a model of computation introduced by Schönfinkel~\cite{Sch24} and
developed by Curry~\cite{Cur30} with the objective to abstain from the need of bound
variables specific to the $\lambda$-calculus. Its combinatorial heart is formed by terms,
which are binary trees with labeled leaves, and rules to compute a result from a term, which
are rewrite relations on trees~\cite{BN98,BKVT03} (clear and complete modern references
about combinatory logic are~\cite{Bar81,HS08,Bim11,Wol21}). An important instance of
combinatory logic is the system containing the basic combinators $\Combinator{K}$ and
$\Combinator{S}$ together with the two rewrite rules
\begin{equation}
    \begin{tikzpicture}[Centering,xscale=0.3,yscale=0.3]
        \node[LeafST](0)at(0.00,-3.33){$\Combinator{K}$};
        \node[LeafST](2)at(2.00,-3.33){$\VarX_1$};
        \node[LeafST](4)at(4.00,-1.67){$\VarX_2$};
        \node[NodeST](1)at(1.00,-1.67){$\App$};
        \node[NodeST](3)at(3.00,0.00){$\App$};
        \draw[Edge](0)--(1);
        \draw[Edge](1)--(3);
        \draw[Edge](2)--(1);
        \draw[Edge](4)--(3);
        \node(r)at(3.00,1.25){};
        \draw[Edge](r)--(3);
    \end{tikzpicture}
    \Rew
    \begin{tikzpicture}[Centering,xscale=0.3,yscale=0.3]
        \node[LeafST](0)at(0.00,0.00){$\VarX_1$};
        \node(r)at(0.00,1.75){};
        \draw[Edge](r)--(0);
    \end{tikzpicture}
    \qquad \mbox{ and } \qquad
    \begin{tikzpicture}[Centering,xscale=0.26,yscale=0.28]
        \node[LeafST](0)at(0.00,-5.25){$\Combinator{S}$};
        \node[LeafST](2)at(2.00,-5.25){$x_1$};
        \node[LeafST](4)at(4.00,-3.50){$x_2$};
        \node[LeafST](6)at(6.00,-1.75){$x_3$};
        \node[NodeST](1)at(1.00,-3.50){$\App$};
        \node[NodeST](3)at(3.00,-1.75){$\App$};
        \node[NodeST](5)at(5.00,0.00){$\App$};
        \draw[Edge](0)--(1);
        \draw[Edge](1)--(3);
        \draw[Edge](2)--(1);
        \draw[Edge](3)--(5);
        \draw[Edge](4)--(3);
        \draw[Edge](6)--(5);
        \node(r)at(5.00,1.31){};
        \draw[Edge](r)--(5);
    \end{tikzpicture}
    \Rew
    \begin{tikzpicture}[Centering,xscale=0.26,yscale=0.22,font=\scriptsize]
        \node[LeafST](0)at(0.00,-4.67){$\VarX_1$};
        \node[LeafST](2)at(2.00,-4.67){$\VarX_3$};
        \node[LeafST](4)at(4.00,-4.67){$\VarX_2$};
        \node[LeafST](6)at(6.00,-4.67){$\VarX_3$};
        \node[NodeST](1)at(1.00,-2.33){$\App$};
        \node[NodeST](3)at(3.00,0.00){$\App$};
        \node[NodeST](5)at(5.00,-2.33){$\App$};
        \draw[Edge](0)--(1);
        \draw[Edge](1)--(3);
        \draw[Edge](2)--(1);
        \draw[Edge](4)--(5);
        \draw[Edge](5)--(3);
        \draw[Edge](6)--(5);
        \node(r)at(3.00,1.75){};
        \draw[Edge](r)--(3);
    \end{tikzpicture}.
\end{equation}
In this system, we have for instance the sequence of computation
\begin{equation}
    \begin{tikzpicture}[Centering,xscale=0.21,yscale=0.22]
        \node[LeafST](0)at(0.00,-6.50){$\Combinator{S}$};
        \node[LeafST](10)at(10.00,-4.33){$\Combinator{S}$};
        \node[LeafST](12)at(12.00,-4.33){$\Combinator{S}$};
        \node[LeafST](2)at(2.00,-10.83){$\Combinator{K}$};
        \node[LeafST](4)at(4.00,-10.83){$\Combinator{K}$};
        \node[LeafST](6)at(6.00,-8.67){$\Combinator{S}$};
        \node[LeafST](8)at(8.00,-4.33){$\Combinator{K}$};
        \node[NodeST](1)at(1.00,-4.33){$\App$};
        \node[NodeST](11)at(11.00,-2.17){$\App$};
        \node[NodeST](3)at(3.00,-8.67){$\App$};
        \node[NodeST](5)at(5.00,-6.50){$\App$};
        \node[NodeST](7)at(7.00,-2.17){$\App$};
        \node[NodeST](9)at(9.00,0.00){$\App$};
        \draw[Edge](0)--(1);
        \draw[Edge](1)--(7);
        \draw[Edge](10)--(11);
        \draw[Edge](11)--(9);
        \draw[Edge](12)--(11);
        \draw[Edge](2)--(3);
        \draw[Edge](3)--(5);
        \draw[Edge](4)--(3);
        \draw[Edge](5)--(1);
        \draw[Edge](6)--(5);
        \draw[Edge](7)--(9);
        \draw[Edge](8)--(7);
        \node(r)at(9.00,1.62){};
        \draw[Edge](r)--(9);
    \end{tikzpicture}
    \enspace \RewContext \enspace
    \begin{tikzpicture}[Centering,xscale=0.23,yscale=0.23]
        \node[LeafST](0)at(0.00,-6.75){$\Combinator{S}$};
        \node[LeafST](2)at(2.00,-6.75){$\Combinator{K}$};
        \node[LeafST](4)at(4.00,-4.50){$\Combinator{K}$};
        \node[LeafST](6)at(6.00,-4.50){$\Combinator{S}$};
        \node[LeafST](8)at(8.00,-4.50){$\Combinator{S}$};
        \node[NodeST](1)at(1.00,-4.50){$\App$};
        \node[NodeST](3)at(3.00,-2.25){$\App$};
        \node[NodeST](5)at(5.00,0.00){$\App$};
        \node[NodeST](7)at(7.00,-2.25){$\App$};
        \draw[Edge](0)--(1);
        \draw[Edge](1)--(3);
        \draw[Edge](2)--(1);
        \draw[Edge](3)--(5);
        \draw[Edge](4)--(3);
        \draw[Edge](6)--(7);
        \draw[Edge](7)--(5);
        \draw[Edge](8)--(7);
        \node(r)at(5.00,1.69){};
        \draw[Edge](r)--(5);
    \end{tikzpicture}
    \enspace \RewContext \enspace
    \begin{tikzpicture}[Centering,xscale=0.2,yscale=0.19]
        \node[LeafST](0)at(0.00,-5.50){$\Combinator{K}$};
        \node[LeafST](10)at(10.00,-8.25){$\Combinator{S}$};
        \node[LeafST](2)at(2.00,-8.25){$\Combinator{S}$};
        \node[LeafST](4)at(4.00,-8.25){$\Combinator{S}$};
        \node[LeafST](6)at(6.00,-5.50){$\Combinator{K}$};
        \node[LeafST](8)at(8.00,-8.25){$\Combinator{S}$};
        \node[NodeST](1)at(1.00,-2.75){$\App$};
        \node[NodeST](3)at(3.00,-5.50){$\App$};
        \node[NodeST](5)at(5.00,0.00){$\App$};
        \node[NodeST](7)at(7.00,-2.75){$\App$};
        \node[NodeST](9)at(9.00,-5.50){$\App$};
        \draw[Edge](0)--(1);
        \draw[Edge](1)--(5);
        \draw[Edge](10)--(9);
        \draw[Edge](2)--(3);
        \draw[Edge](3)--(1);
        \draw[Edge](4)--(3);
        \draw[Edge](6)--(7);
        \draw[Edge](7)--(5);
        \draw[Edge](8)--(9);
        \draw[Edge](9)--(7);
        \node(r)at(5.00,2.06){};
        \draw[Edge](r)--(5);
    \end{tikzpicture}
    \enspace \RewContext \enspace
    \begin{tikzpicture}[Centering,xscale=0.35,yscale=0.35]
        \node[LeafST](0)at(0.00,-1.50){$\Combinator{S}$};
        \node[LeafST](2)at(2.00,-1.50){$\Combinator{S}$};
        \node[NodeST](1)at(1.00,0.00){$\App$};
        \draw[Edge](0)--(1);
        \draw[Edge](2)--(1);
        \node(r)at(1.00,1.12){};
        \draw[Edge](r)--(1);
    \end{tikzpicture}
\end{equation}
obtained by performing at each step one rewrite according to the previous rules. This system
is important because it is combinatorially complete: each $\lambda$-term can be translated,
by so-called bracket abstraction algorithms~\cite{Sch24,CF58}, into a term over
$\Combinator{K}$ and $\Combinator{S}$ emulating it.

A lot of other basic combinators with their own rewrite rules have been introduced by
Smullyan in~\cite{Smu85} after ---now widely used--- bird names, forming the enchanted
forest of combinator birds. For instance, $\Combinator{K}$ is the Kestrel and
$\Combinator{S}$ is the Starling. Usual computer science-oriented questions consist in
considering a fragment of combinatory logic, that is a finite set of combinators with their
rewrite rules which is not necessarily combinatorially complete, and ask for the following
questions:
\begin{enumerate}[label={\bf (\alph*)}]
    \item Given two terms $\TreeT$ and $\TreeT'$, can we decide if $\TreeT$ and $\TreeT'$
    can be rewritten eventually in a same term? This is known as the word
    problem~\cite{BN98,Sta00}. This question admits a positive answer for some basic
    combinators like among others the Lark~\cite{Sta89,SWB93} and the Warbler~\cite{SWB93}
    but it still open for the Starling~\cite{BEJW17};

    \item Given a term $\TreeT$, can we decide if all rewrite sequences starting from
    $\TreeT$ are finite? This is known as the strong normalization problem. This question
    admits a positive answer for, among others, the Starling~\cite{Wal00} and the
    Jay~\cite{PS01}. Related to this, see also~\cite{DGKRTZ13,BGZ17} for a probabilistic and
    asymptotic study of strong normalizing terms.
\end{enumerate}

Here, we decide to pursue this study in a different direction by asking questions from a
combinatorial point of view, including the study of order theoretic structures and adopting
an enumerative approach. In particular, by denoting by $\Leq$ (resp.\ by $\Equiv$) the
reflexive and transitive (resp.\ reflexive, symmetric, and transitive) closure of the
rewrite relation, we ask the following questions:
\begin{enumerate}[label={\bf (\alph*')}]
    \item Determine if $\Leq$ is a partial order relation;

    \item In this case, determine if each interval of this poset is a lattice;

    \item Enumerate the $\Equiv$-equivalence classes of terms w.r.t.\ the minimal degrees of
    their terms.
\end{enumerate}
This work fits in this general project consisting in mixing combinatory logic with
combinatorics. These enumerative questions, including the enumeration w.r.t.\ some size
notions (like the height or the number of basic combinators of the terms) of some particular
terms of a system (like its normal forms or the minimal or maximal elements of its partial
order) lead to density properties of the counted terms among the terms of the same size.
This gives information about the probability for a term of a given size generated uniformly
at random to satisfy some property. Besides, some of the mentioned questions can be
addressed by working with rewrite graphs, that are graphs labeled by terms and where there
is a directed edge between two terms if the first can be obtained from the second by a
rewriting step. Such graphs appear naturally in the study of term rewrite
systems~\cite{BL79}. Some different graphs have been considered in the literature, where
vertices are no longer terms but rather sequences of rewriting steps (see for
instance~\cite{BL79,VZ84} and~\cite[Exercise 12.4.2.]{Bar81}). These graphs, as well as some
of their semilattice properties~\cite[Theorem 2.3.1.]{BL79}, have no direct links with the
present work.

We choose here to start this project by studying the system made of a single combinator, the
combinator $\M$, known as the Mockingbird~\cite{Smu85} or as the little omega. This
combinator is a very simple and an important one (see for instance~\cite{Sta11,Sta17}). By
drawing the portions of the rewrite graph starting from terms on $\M$, the first properties
that stand out are that the graph does not contain any nontrivial loops and that its
connected components are finite and have exactly one minimal and one maximal element. At
this stage, driven by computer exploration, we conjecture that the relation $\Leq$ on the
terms on $\M$ is a partial order relation and that each $\Equiv$-equivalence class is a
lattice w.r.t.\ this partial order relation. This lattice property is for us a good clue for
the fact that this system contains rather rich combinatorial properties.

To prove this last property, we introduce new lattices on duplicative forests, that are
kinds of treelike structures, and show that each maximal interval from any term on $\M$ is
isomorphic as a poset to a maximal interval of a lattice of duplicative forests
(Theorem~\ref{thm:mockingbird_lattices}). We define the Mockingbird lattice of order $d \geq
0$ as the lattice $\MockingbirdLattice(d)$ consisting in the closed terms on $\M$ equal as
or greater than the right comb closed term on $\M$ of degree $d$. We prove that each
interval of the poset of duplicative forests is contained as an interval in
$\MockingbirdLattice(d)$ for a certain $d \geq 0$
(Theorem~\ref{thm:universality_mockingbird_lattices}). Since any closed term on $\M$ can be
seen as a binary tree, this provides a new lattice structure on these objects. A lot of
similar lattices have been studied on binary trees such as, among others the very famous
Tamari lattice~\cite{Tam62}, the Kreweras lattice~\cite{Kre72}, the Stanley
lattice~\cite{Sta75}, the phagocyte lattice~\cite{BP06}, and the pruning-grafting
lattice~\cite{BP08}. However, unlike these lattices having for each order $d \geq 0$ a
cardinality equal to the $d$-th Catalan number, the lattices $\MockingbirdLattice(d)$ are
enumerated by a different integer sequence. To obtain enumerative results about the
Mockingbird lattices and all the posets of terms on $\M$ in general, we use formal power
series on terms and on duplicative forests, and several products on these. In this way, we
enumerate the minimal and maximal elements of the infinite poset of the closed terms on $\M$
(Propositions~\ref{prop:maximal_elements_enumeration}
and~\ref{prop:minimal_elements_enumeration}), the lengths of the shortest and longest
saturated chains of $\MockingbirdLattice(d)$ (Proposition~\ref{prop:paths_lengths}), the
cardinality of $\MockingbirdLattice(d)$ (Proposition~\ref{prop:number_elements}), the number
of edges of the Hasse diagram of $\MockingbirdLattice(d)$
(Proposition~\ref{prop:number_coverings}), and the number of intervals of
$\MockingbirdLattice(d)$ (Proposition~\ref{prop:number_intervals}). We also provide general
results for systems of combinatory logic: we give in particular a necessary condition on the
rewrite rules in order to have only finite connected components in the rewrite graph
(Proposition~\ref{prop:finite_equivalence_classes}). This relies on a combinatorial property
on basic combinators, called the hierarchical property. We also discuss some consequences of
this fact in order to construct models for systems having this property. These are in fact
the algebras over a certain abstract clone~\cite{Tay93} defined from the system. When the
relation $\Leq$ is a partial order relation, such models have a potential nice algorithmic
computational complexity.

This paper is organized as follows. Section~\ref{sec:terms_rewrite_cls} contains the
preliminary notions and definitions about terms, rewrite relations, and combinatory logic
systems. We show here some general properties of these systems. In
Section~\ref{sec:mockingbird_lattice}, we study the combinatory logic system on $\M$ and the
Mockingbird lattices. Finally, Section~\ref{sec:enumerative_properties} contains our
enumerative results about Mockingbird lattices. This text ends with the presentation of some
open questions raised by this work.

This paper is an extended version of~\cite{Gir22}, announcing the main results without any
proofs. The present version contains all proofs of the stated results, more examples, and
presents new results as discussions about some models of the combinatory logic system on
$\M$ and the enumeration of isolated elements of the poset of the closed terms on~$\M$. This
version uses also the formalism of abstract clones in order to work with terms and rewrite
systems.

\subsubsection*{General notations and conventions}
If $S$ is a finite set $\# S$ is the cardinality of $S$. For any integers $i$ and $j$, $[i,
j]$ denotes the set $\{i, i + 1, \dots, j\}$. For any integer $i$, $[i]$ denotes the set
$[1, i]$ and $\HanL{i}$ denotes the set $[0, i]$. For any set $A$, $A^*$ is the set of all
words on $A$. For any $w \in A^*$, $\Length(w)$ is the length of $w$, for any $i \in
[\Length(w)]$, $w(i)$ is the $i$-th letter of $w$, and $|w|_a$ is the number of occurrences
of the letter $a \in A$ in $w$. The only word of length $0$ is the empty word $\epsilon$.
For any statement $P$, the Iverson bracket $\Iverson{P}$ takes $1$ as value if $P$ is true
and $0$ otherwise.

\section{Terms, rewrite relations, and combinatory logic systems}
\label{sec:terms_rewrite_cls}
In this preliminary part, we set the main definitions and notations used in the sequel. We
introduce also the central notion of combinatory logic systems, which are defined from
special kinds of rewrite relations. We show also some general properties of these systems.

\subsection{Terms, compositions, and rewrite relations}
Let us start by presenting here the notions of terms, composition, and rewrite relations.

\subsubsection{Terms}
An \Def{alphabet} is a finite set $\GeneratingSet$. Its elements are called \Def{constants}
or \Def{basic combinators}. Any element of the set $\SetVariables := \bigcup_{n \geq 1}
\SetVariables_n$, where $\SetVariables_n := \Bra{\VarX_1, \dots, \VarX_n}$, is a
\Def{variable}. The set $\SetTerms(\GeneratingSet)$ of \Def{$\GeneratingSet$-terms} (or
simply \Def{terms} when the context is clear) is so that any variable of $\SetVariables$ is
a $\GeneratingSet$-term, any constant of $\GeneratingSet$ is a $\GeneratingSet$-term, and if
$\TreeT_1$ and $\TreeT_2$ are two $\GeneratingSet$-terms, then $\Par{\TreeT_1 \App
\TreeT_2}$ is a $\GeneratingSet$-term. For any $n \geq 0$, we denote by
$\SetTerms_n(\GeneratingSet)$ the set of $\GeneratingSet$-terms $\TreeT$ having all
variables belonging to $\SetVariables_n$.

From this definition, any term is a rooted planar binary tree where leaves are decorated by
variables or by constants. We shall express terms concisely by removing superfluous
parentheses by considering that $\App$ associates to the left and also by removing the
symbols $\App$ (see the example given in~\eqref{equ:example_term_1},
\eqref{equ:example_term_2}, and~\eqref{equ:example_term}).

Let $\TreeT$ be a $\GeneratingSet$-term. The \Def{degree} $\Deg(\TreeT)$ of $\TreeT$ is the
number of internal nodes (the nodes decorated by $\App$) of $\TreeT$ seen as a binary tree.
For any $a \in \SetVariables \cup \GeneratingSet$, $\Deg_a(\TreeT)$ is the number of nodes
labeled by $a$ in $\TreeT$. The \Def{depth} of a node $u$ of $\TreeT$ is the number of
internal nodes in the path connecting the root of $\TreeT$ and $u$. The \Def{height}
$\Height(\TreeT)$ of $\TreeT$ is the maximal depth among all the nodes of $\TreeT$. The
\Def{frontier} of $\TreeT$ is the sequence of all the variables appearing in $\TreeT$ from
the left to the right. The term $\TreeT$ is \Def{planar} if its frontier is of the form
$\VarX_1 \VarX_2 \dots \VarX_n$ where $n \geq 0$. The term $\TreeT$ is \Def{linear} if there
are no multiple occurrences of the same variable in its frontier. A \Def{closed term} or
\Def{combinator} is a term of $\SetTerms_0(\GeneratingSet)$ having thus no occurrence of any
variable.

For instance, by setting $\GeneratingSet := \Bra{\Combinator{A}, \Combinator{B}}$, the
$\GeneratingSet$-term
\begin{equation} \label{equ:example_term_1}
    \TreeT := \Par{\Par{\Par{\Combinator{A} \App \VarX_3} \App \Par{\Combinator{B}
    \App \Par{\Par{\VarX_1 \App \VarX_3} \App \VarX_1}}} \App
    \Par{\Combinator{A} \App \Combinator{A}}}
\end{equation}
draws as the binary tree
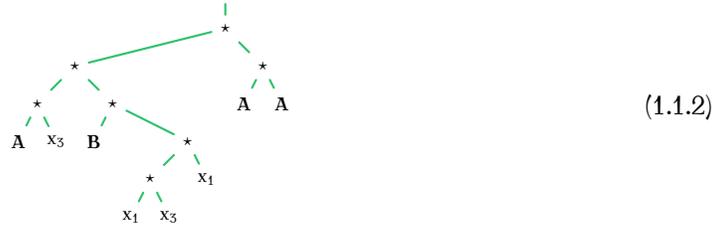
\begin{equation} \label{equ:example_term_2}
    \begin{tikzpicture}[Centering,xscale=0.25,yscale=0.2]
        \node[LeafST](0)at(0.00,-7.50){$\Combinator{A}$};
        \node[LeafST](10)at(10.00,-10.00){$\VarX_1$};
        \node[LeafST](12)at(12.00,-5.00){$\Combinator{A}$};
        \node[LeafST](14)at(14.00,-5.00){$\Combinator{A}$};
        \node[LeafST](2)at(2.00,-7.50){$\VarX_3$};
        \node[LeafST](4)at(4.00,-7.50){$\Combinator{B}$};
        \node[LeafST](6)at(6.00,-12.50){$\VarX_1$};
        \node[LeafST](8)at(8.00,-12.50){$\VarX_3$};
        \node[NodeST](1)at(1.00,-5.00){$\App$};
        \node[NodeST](11)at(11.00,0.00){$\App$};
        \node[NodeST](13)at(13.00,-2.50){$\App$};
        \node[NodeST](3)at(3.00,-2.50){$\App$};
        \node[NodeST](5)at(5.00,-5.00){$\App$};
        \node[NodeST](7)at(7.00,-10.00){$\App$};
        \node[NodeST](9)at(9.00,-7.50){$\App$};
        \draw[Edge](0)--(1);
        \draw[Edge](1)--(3);
        \draw[Edge](10)--(9);
        \draw[Edge](12)--(13);
        \draw[Edge](13)--(11);
        \draw[Edge](14)--(13);
        \draw[Edge](2)--(1);
        \draw[Edge](3)--(11);
        \draw[Edge](4)--(5);
        \draw[Edge](5)--(3);
        \draw[Edge](6)--(7);
        \draw[Edge](7)--(9);
        \draw[Edge](8)--(7);
        \draw[Edge](9)--(5);
        \node(r)at(11.00,2.25){};
        \draw[Edge](r)--(11);
    \end{tikzpicture}
\end{equation}
and expresses concisely as
\begin{equation} \label{equ:example_term}
    \TreeT
    = \Par{\Combinator{A} \VarX_3} \Par{\Combinator{B} \Par{\VarX_1 \VarX_3 \VarX_1}}
    \Par{\Combinator{A} \Combinator{A}}.
\end{equation}
Its degree is $7$, its height is $5$, its frontier is $\VarX_3 \VarX_1 \VarX_3 \VarX_1$, and
is an element $\SetTerms_n(\GeneratingSet)$ for any $n \geq 3$.

\subsubsection{Compositions of terms and clones}
Let $\TreeT$ and $\TreeT'_1$, \dots, $\TreeT'_n$, $n \geq 0$, be $\GeneratingSet$-terms. The
\Def{composition} of $\TreeT$ with $\TreeT'_1$, \dots, $\TreeT'_n$ is the
$\GeneratingSet$-term $\TreeT \Han{\TreeT'_1, \dots, \TreeT'_n}$ obtained by simultaneously
replacing for all $i \in [n]$ all occurrences of the variables $\VarX_i$ in $\TreeT$
by~$\TreeT'_i$. For instance
\begin{equation}
    \Combinator{A} \Par{\VarX_3 \VarX_1} \VarX_3
    \Han{\VarX_1 \VarX_4, \VarX_1 \VarX_1 \VarX_2,
    \Combinator{A} \Par{\Combinator{A} \VarX_1}}
    =
    \Combinator{A} \Par{\Combinator{A} \Par{\Combinator{A} \VarX_1} \Par{\VarX_1 \VarX_4}}
    \Par{\Combinator{A} \Par{\Combinator{A} \VarX_1}}.
\end{equation}
Given two $\GeneratingSet$-terms $\TreeT$ and $\TreeS$, $\TreeS$ is a \Def{factor} of
$\TreeT$ if
\begin{equation} \label{equ:factor}
    \TreeT
    = \TreeT' \Han{\TreeS_1, \dots, \TreeS_{i - 1}, \TreeS, \TreeS_{i + 1}, \dots, \TreeS_n}
    \Han{\TreeR_1, \dots, \TreeR_m}
\end{equation}
for some integers $n, m \geq 0$ and $\GeneratingSet$-terms $\TreeT'$, $\TreeS_1$, \dots,
$\TreeS_{i - 1}$, $\TreeS_{i + 1}$, \dots, $\TreeS_n$, $\TreeR_1$, \dots, $\TreeR_m$,
where $\VarX_i$ appears in $\TreeT'$. When this property does not hold, $\TreeT$
\Def{avoids} $\TreeS$. When in~\eqref{equ:factor}, for all $i \in [m]$ the $\TreeR_i$
are variables, $\TreeS$ is a \Def{suffix} of $\TreeT$. When this property does not hold,
$\TreeT$ \Def{suffix-avoids} $\TreeS$. For instance, the term defined
in~\eqref{equ:example_term} admits $\Combinator{B} \Par{\VarX_1 \VarX_2}$ as factor,
$\VarX_1 \VarX_2 \VarX_3$ as suffix, but suffix-avoids $\VarX_1 \VarX_1 \VarX_1$.

Observe that we have, for any $n, m \geq 0$, $i \in [n]$, and $\TreeT, \TreeS_1, \dots,
\TreeS_n, \TreeR_1, \dots, \TreeR_m \in \SetTerms(\GeneratingSet)$,
\begin{subequations}
\begin{equation}
    \VarX_i \Han{\TreeT_1, \dots, \TreeT_n} = \TreeT_i,
\end{equation}
\begin{equation}
    \TreeT \Han{\VarX_1, \dots \VarX_n} = \TreeT,
\end{equation}
\begin{equation}
    \TreeT \Han{\TreeS_1, \dots, \TreeS_n} \Han{\TreeR_1, \dots, \TreeR_m}
    =
    \TreeT \Han{\TreeS_1 \Han{\TreeR_1, \dots, \TreeR_m}, \dots,
        \TreeS_n \Han{\TreeR_1, \dots, \TreeR_m}}.
\end{equation}
\end{subequations}
These three relations imply that the set $\SetTerms(\GeneratingSet)$ together with the
composition operation is an abstract clone~\cite{Tay93}. This is in fact the free abstract
clone generated by $\GeneratingSet \sqcup \{\App\}$ where $\App$ is a binary generator.

\subsubsection{Rewrite relations}
A \Def{rewrite relation} on $\SetTerms(\GeneratingSet)$ is a binary relation $\Rew$ on
$\SetTerms(\GeneratingSet)$. The \Def{context closure} of $\Rew$ is the binary relation
$\RewContext$ on $\SetTerms(\GeneratingSet)$ satisfying
\begin{equation} \label{equ:context_closure}
    \TreeT \Han{\TreeS_1, \dots, \TreeS_i, \dots, \TreeS_n}
    \Han{\TreeR_1, \dots, \TreeR_m}
    \RewContext
    \TreeT \Han{\TreeS_1, \dots, \TreeS_i', \dots, \TreeS_n}
    \Han{\TreeR_1, \dots, \TreeR_m}
\end{equation}
where $\TreeT$ is a planar $\GeneratingSet$-term, $\VarX_i$ appears in $\TreeT$, and
$\TreeS_1$, \dots, $\TreeS_i$, $\TreeS_i'$, \dots, $\TreeS_n$, $n \geq 1$, $i \in [n]$, and
$\TreeR_1$, \dots, $\TreeR_m$, $m \geq 0$, are $\GeneratingSet$-terms, and $\TreeS_i \Rew
\TreeS_i'$. Intuitively, $\TreeT \RewContext \TreeT'$ if $\TreeT'$ can be obtained from
$\TreeT$ by replacing a factor $\TreeS$ of $\TreeT$ by $\TreeS'$ and by identifying
variables with terms in a coherent way, whenever $\TreeS \Rew \TreeS'$. For instance, for
$\GeneratingSet := \Bra{\Combinator{A}, \Combinator{B}, \Combinator{C}}$ and the rewrite
relation $\Rew$ satisfying
\begin{math}
    \VarX_1 \Par{\Combinator{A} \VarX_1} \Rew \VarX_1 \VarX_2,
\end{math}
we have
\begin{equation}
    \VarX_3 \Combinator{B} \Par{\Combinator{A} \Par{\VarX_3 \Combinator{B}}}
    \Par{\VarX_4 \Combinator{A}}
    \RewContext
    \VarX_3 \Combinator{B} \Combinator{C} \Par{\VarX_4 \Combinator{A}}
\end{equation}
because
\begin{equation}
    \VarX_3 \Combinator{B}
    \Par{\Combinator{A} \Par{\VarX_3 \Combinator{B}}}
    \Par{\VarX_4 \Combinator{A}}
    =
    \VarX_1 \Par{\VarX_2 \Combinator{A}}
    \
    \Han{\VarX_1 \Par{\Combinator{A} \VarX_1}, \VarX_4}
    \
    \Han{\VarX_3 \Combinator{B}, \Combinator{C}, \VarX_3, \VarX_4}
\end{equation}
and
\begin{equation}
    \VarX_3 \Combinator{B} \Combinator{C} \Par{\VarX_4 \Combinator{A}}
    =
    \VarX_1 \Par{\VarX_2 \Combinator{A}}
    \
    \Han{\VarX_1 \VarX_2, \VarX_4}
    \
    \Han{\VarX_3 \Combinator{B}, \Combinator{C}, \VarX_3, \VarX_4}.
\end{equation}
Observe in this example that the variable $\VarX_2$ occurs in the right member of $\Rew$ and
not in the left member. For this reason, this variable can be replaced by any term in the
context closure of $\Rew$ (it is replaced by $\Combinator{C}$ in the example).

\subsubsection{Applicative term rewrite systems}
An \Def{applicative term rewrite system} (or \Def{ATRS} for short) is a pair
$(\GeneratingSet, \Rew)$ such that $\GeneratingSet$ is an alphabet and $\Rew$ is a rewrite
relation on $\SetTerms(\GeneratingSet)$. Let $\ATRS := (\GeneratingSet, \Rew)$ be an ATRS.
We denote by $\Leq$ the reflexive and transitive closure of $\RewContext$. This relation
$\Leq$ is by construction a preorder. We denote by $\Equiv$ the reflexive, symmetric, and
transitive closure of $\RewContext$ and by $\Han{\TreeT}_\Equiv$ the $\Equiv$-equivalence
class of $\TreeT \in \SetTerms(\GeneratingSet)$. This relation $\Equiv$ is by construction
an equivalence relation. The \Def{rewrite graph} $\RewGraph_\ATRS$ of $\ATRS$ is the digraph
on the set $\SetTerms(\GeneratingSet)$ of vertices where there is an arc from $\TreeT \in
\SetTerms(\GeneratingSet)$ to $\TreeT' \in \SetTerms(\GeneratingSet)$ if $\TreeT \RewContext
\TreeT'$. For any $\TreeT \in \SetTerms(\GeneratingSet)$, we also denote by
$\RewGraph_\ATRS(\TreeT)$ the subgraph of $\RewGraph_\ATRS$ restrained on the set
\begin{math}
    \Bra{\TreeT' \in \SetTerms(\GeneratingSet) : \TreeT \Leq \TreeT'}
\end{math}
of vertices. The ATRS $\ATRS$ is \Def{locally finite} if for any term $\TreeT$,
$\Han{\TreeT}_{\Equiv}$ is finite. This is equivalent to the fact that all the connected
components of $\RewGraph_\ATRS$ are finite. A $\GeneratingSet$-term $\TreeT$ is a
\Def{normal form} of $\ATRS$ if there is no arc of source $\TreeT$ in $\RewGraph_\ATRS$. We
say that $\TreeT$ is \Def{weakly normalizing} if there is at least one normal form in
$\RewGraph_\ATRS(\TreeT)$ and that $\TreeT$ is \Def{strongly normalizing} if
$\RewGraph_\ATRS(\TreeT)$ is finite and acyclic. When all the $\GeneratingSet$-terms are
strongly normalizing, $\ATRS$ is \Def{terminating}. Besides, if for any $\TreeT, \TreeS_1,
\TreeS_2 \in \SetTerms(\GeneratingSet)$, $\TreeT \Leq \TreeS_1$ and $\TreeT \Leq \TreeS_2$
implies the existence of $\TreeT' \in \SetTerms(\GeneratingSet)$ such that $\TreeS_1 \Leq
\TreeT'$ and $\TreeS_2 \Leq \TreeT'$, then $\ATRS$ is \Def{confluent}.

\subsubsection{Partial orders}
When $\ATRS := (\GeneratingSet, \Rew)$ is an ATRS such that $\Leq$ is antisymmetric, $\Leq$
is a partial order relation. In this case, $\ATRS$ has the \Def{poset property} and we shall
denote by $\Poset_\ATRS$ the poset $\Par{\SetTerms(\GeneratingSet), \Leq}$. For any $\TreeT
\in \SetTerms(\GeneratingSet)$, let $\Poset_\ATRS(\TreeT)$ be the subposet of $\Poset_\ATRS$
having $\TreeT$ as least element. Again in this case, a $\GeneratingSet$-term $\TreeT$ is
\Def{minimal} (resp.\ \Def{maximal}) if $\TreeT$ is a minimal (resp.\ maximal) element of
$\Poset_\ATRS$. Observe that any normal form is maximal but the converse is false because a
maximal element $\TreeT$ could satisfy $\TreeT \RewContext \TreeT$. A $\GeneratingSet$-term
$\TreeT$ is \Def{isolated} if $\TreeT$ is both minimal and maximal. When $\ATRS$ has the
poset property, $\ATRS$ is \Def{rooted} if for any term $\TreeT$, the subposet
$\Han{\TreeT}_{\Equiv}$ of $\Poset_\ATRS$ has a unique minimal element. Finally, $\ATRS$ has
the \Def{lattice property} if $\ATRS$ has the poset property and for all $\TreeT \in
\SetTerms(\GeneratingSet)$, all posets $\Poset_\ATRS(\TreeT)$ are lattices.

\subsection{Combinatory logic systems}
We begin by defining combinatory logic systems as particular ATRS and present some
consequence of nonerasing and hierarchical combinatory logic systems (notions defined
thereafter). We explain some of the consequences of the local finiteness, the poset
property, and the lattice property on the models of combinatory logic systems.

\subsubsection{Main definitions}
A \Def{combinatory logic system} (or \Def{CLS} for short) is an ATRS $\CLS :=
(\GeneratingSet, \Rew)$ such that for each constant $\Combinator{X}$ of $\GeneratingSet$,
there is exactly one rewrite rule where $\Combinator{X}$ appears, and this rule is of the
form
\begin{equation} \label{equ:rule_CLS}
    \Combinator{X} \VarX_1 \dots \VarX_n \Rew \TreeT_{\Combinator{X}}
\end{equation}
where $n \geq 1$ and $\TreeT_{\Combinator{X}}$ is a term having no constants and having all
variables in $\SetVariables_n$. Moreover, all rewrite rules of $\CLS$ must be of the
form~\eqref{equ:rule_CLS}. The integer $n$ is the \Def{order} of $\Combinator{X}$ in $\CLS$.
Some well-known terms $\TreeT_{\Combinator{X}}$ appearing among others in~\cite{Smu85} are,
\begin{itemize}
    \item with order $1$, $\TreeT_{\Combinator{I}} := \VarX_1$ (\Def{Identity bird}),
    $\TreeT_\M := \VarX_1 \VarX_1$ (\Def{Mockingbird});
    \item with order $2$, $\TreeT_{\Combinator{K}} := \VarX_1$ (\Def{Kestrel}),
    $\TreeT_{\Combinator{T}} := \VarX_2 \VarX_1$ (\Def{Thrush}), $\TreeT_{\Combinator{M_1}}
    := \VarX_1 \VarX_1 \VarX_2$ (\Def{Mockingbird 1}), $\TreeT_{\Combinator{W}} := \VarX_1
    \VarX_1 \VarX_2$ (\Def{Warbler}), $\TreeT_{\Combinator{L}} := \VarX_1 \Par{\VarX_2
    \VarX_2}$ (\Def{Lark}), $\TreeT_{\Combinator{O}} := \VarX_2 \Par{\VarX_1 \VarX_2}$
    (\Def{Owl}), $\TreeT_{\Combinator{U}} := \VarX_2 \Par{\VarX_1 \VarX_1 \VarX_2}$
    (\Def{Turing Bird});
    \item with order $3$, $\TreeT_{\Combinator{C}} := \VarX_1 \VarX_3 \VarX_2$
    (\Def{Cardinal}), $\TreeT_{\Combinator{V}} := \VarX_3 \VarX_1 \VarX_2$ (\Def{Vireo}),
    $\TreeT_{\Combinator{B}} := \VarX_1 \Par{\VarX_2 \VarX_3}$ (\Def{Bluebird}),
    $\TreeT_{\Combinator{S}} := \VarX_1 \VarX_3 \Par{\VarX_2 \VarX_3}$ (\Def{Starling});
    \item with order $4$, $\TreeT_{\Combinator{J}} := \VarX_1 \VarX_2 \Par{\VarX_1 \VarX_4
    \VarX_3}$ (\Def{Jay}).
\end{itemize}
The constant $\Combinator{X}$ is \Def{nonerasing} if $\TreeT_{\Combinator{X}}$ contains at
least one occurrence of each variable $\VarX_i$ for any $i \in [n]$. The constant
$\Combinator{X}$ is \Def{hierarchical} if for any $i \in [n]$, $\VarX_i$ appears in
$\TreeT_{\Combinator{X}}$ at depth $n + 1 - i$. For instance, the terms
$\TreeT_{\Combinator{X}}$ such that $\Combinator{X}$ are hierarchical and of order $3$ or
less are
\begin{equation}
    \VarX_1 \VarX_1, \enspace
    \VarX_1 \VarX_1 \VarX_2, \enspace
    \VarX_2 \Par{\VarX_1 \VarX_1}, \enspace
    \VarX_1 \VarX_1 \VarX_2 \VarX_3, \enspace
    \VarX_2 \Par{\VarX_1 \VarX_1} \VarX_3, \enspace
    \VarX_3 \Par{\VarX_1 \VarX_1 \VarX_2}, \enspace
    \VarX_3 \Par{\VarX_2 \Par{\VarX_1 \VarX_1}}.
\end{equation}
In particular, if $\Combinator{X}$ is hierarchical, then $\Combinator{X}$ is nonerasing. We
say that $\CLS$ is \Def{nonerasing} (resp.\ \Def{hierarchical}) if all constants of
$\GeneratingSet$ are nonerasing (resp.\ hierarchical). Other interesting properties of basic
combinators are introduced in~\cite{Bim11} but they do not intervene directly in this work.

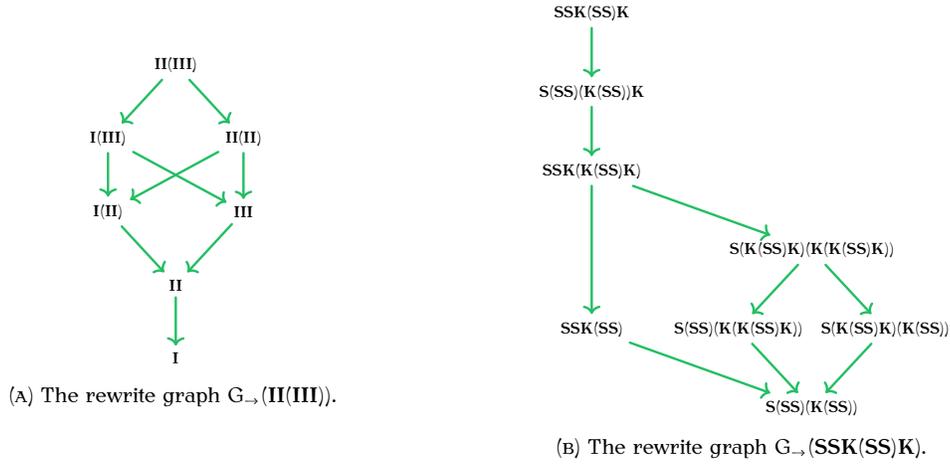
\begin{figure}[ht]
    \centering
    \newcommand{\I}{{\Combinator{I}}}
    \newcommand{\Kk}{{\Combinator{K}}}
    \newcommand{\Ss}{{\Combinator{S}}}
    \subfloat[][The rewrite graph $\RewGraph_\Rew(\I \I (\I \I \I))$.]{
    \begin{Page}{.46}
        \centering
        \scalebox{.75}{
        \begin{tikzpicture}[Centering,xscale=1.2,yscale=1.3]
            \node[GraphLabeledVertex](0)at(0,0){$\I \I (\I \I \I)$};
            \node[GraphLabeledVertex](1)at(-1,-1){$\I (\I \I \I)$};
            \node[GraphLabeledVertex](2)at(1,-1){$\I \I (\I \I)$};
            \node[GraphLabeledVertex](3)at(-1,-2){$\I (\I \I)$};
            \node[GraphLabeledVertex](4)at(1,-2){$\I \I \I$};
            \node[GraphLabeledVertex](5)at(0,-3){$\I \I$};
            \node[GraphLabeledVertex](6)at(0,-4){$\I$};
            \draw[GraphArc](0)--(1);
            \draw[GraphArc](0)--(2);
            \draw[GraphArc](1)--(3);
            \draw[GraphArc](1)--(4);
            \draw[GraphArc](2)--(3);
            \draw[GraphArc](2)--(4);
            \draw[GraphArc](3)--(5);
            \draw[GraphArc](4)--(5);
            \draw[GraphArc](5)--(6);
        \end{tikzpicture}}
    \end{Page}
    \label{subfig:rewrite_graph_I}}
    \hfill
    \subfloat[][The rewrite graph $\RewGraph_\Rew(\Ss\Ss\Kk(\Ss\Ss)\Kk)$.]{
    \begin{Page}{.5}
        \centering
        \scalebox{.75}{
        \begin{tikzpicture}[Centering,xscale=1.3,yscale=1.4]
            \node[GraphLabeledVertex](0)at(-4,0){$\Ss\Ss\Kk(\Ss\Ss)\Kk$};
            \node[GraphLabeledVertex](1)at(-4,-1){$\Ss(\Ss\Ss)(\Kk(\Ss\Ss))\Kk$};
            \node[GraphLabeledVertex](2)at(-4,-2){$\Ss\Ss\Kk(\Kk(\Ss\Ss)\Kk)$};
            \node[GraphLabeledVertex](3)at(-1,-3)
                {$\Ss(\Kk(\Ss\Ss)\Kk)(\Kk(\Kk(\Ss\Ss)\Kk))$};
            \node[GraphLabeledVertex](4)at(-4,-4){$\Ss\Ss\Kk(\Ss\Ss)$};
            \node[GraphLabeledVertex](5)at(-2,-4){$\Ss(\Ss\Ss)(\Kk(\Kk(\Ss\Ss)\Kk))$};
            \node[GraphLabeledVertex](6)at(0,-4){$\Ss(\Kk(\Ss\Ss)\Kk)(\Kk(\Ss\Ss))$};
            \node[GraphLabeledVertex](7)at(-1,-5){$\Ss(\Ss\Ss)(\Kk(\Ss\Ss))$};
            \draw[GraphArc](0)--(1);
            \draw[GraphArc](1)--(2);
            \draw[GraphArc](2)--(3);
            \draw[GraphArc](2)--(4);
            \draw[GraphArc](3)--(5);
            \draw[GraphArc](3)--(6);
            \draw[GraphArc](5)--(7);
            \draw[GraphArc](6)--(7);
            \draw[GraphArc](4)--(7);
        \end{tikzpicture}}
    \end{Page}
    \label{subfig:rewrite_graph_KS}}
    \caption{Some subgraphs of rewrite graphs of some CLS.}
    \label{fig:rewrite_graphs}
\end{figure}

Consider for instance the CLS $\CLS := (\GeneratingSet, \Rew)$ where $\GeneratingSet$
contains only the constant $\Combinator{I}$ where $\TreeT_{\Combinator{I}}$ is defined
above. It is straightforward to show that $\CLS$ has the poset property. Nevertheless,
$\CLS$ has not the lattice property, as suggested by the Hasse diagram shown in
Figure~\ref{subfig:rewrite_graph_I}.
On the other side, the CLS $\CLS := (\GeneratingSet, \Rew)$ where $\GeneratingSet$ contains
only the constants $\Combinator{K}$ and $\Combinator{S}$ where $\TreeT_{\Combinator{K}}$ and
$\TreeT_{\Combinator{S}}$ are defined above does not have the poset property. Indeed, $\Leq$
is not antisymmetric because by setting
\begin{math}
    \newcommand{\Kk}{{\Combinator{K}}}
    \newcommand{\Ss}{{\Combinator{S}}}
    \TreeT :=
    \Ss (\Ss (\Kk \Ss) (\Kk (\Ss \Kk \Kk))) (\Kk (\Ss \Kk \Kk)),
\end{math}
we have $\TreeS \Leq \TreeS'$ and $\TreeS' \Leq \TreeS$ where $\TreeS := \TreeT \TreeT
(\TreeT \TreeT \TreeT)$ and $\TreeS' := \TreeT \TreeT \TreeT (\TreeT \TreeT \TreeT)$. This
can be seen by noticing that $\TreeT \, \VarX_1 \VarX_2 \Leq \VarX_2 \VarX_2$.
Figure~\ref{subfig:rewrite_graph_KS} shows a part of the rewrite graph of~$\CLS$.

\subsubsection{Properties of CLS} \label{subsubsec:properties_CLS}
Let us state some properties of CLS related to their properties of confluence, termination,
and local finiteness.

\begin{Proposition} \label{prop:confluence}
    Any CLS is confluent.
\end{Proposition}
\begin{proof}
    By definition, the underlying ATRS of any CLS is orthogonal~\cite[Chapter 6]{BN98}. The
    statement follows from the fact that all orthogonal rewrite systems are
    confluent~\cite{Ros73}.
\end{proof}

\begin{Lemma} \label{lem:nonerasing_equivalence_preservation}
    Let $\CLS := (\GeneratingSet, \Rew)$ be a CLS. If $\CLS$ is nonerasing, then for any $n
    \geq 0$ and $\TreeT \in \SetTerms_n(\GeneratingSet)$, $\Han{\TreeT}_\Equiv \subseteq
    \SetTerms_n(\GeneratingSet)$.
\end{Lemma}
\begin{proof}
    Assume that $\TreeT$ and $\TreeT'$ are two $\GeneratingSet$-terms such that $\TreeT
    \RewContext \TreeT'$. From the definition of $\RewContext$ from $\Rew$ provided
    by~\eqref{equ:context_closure}, the fact that all constants of $\GeneratingSet$ are
    nonerasing implies that for any $\VarX_i \in \SetVariables$, $\Deg_{\VarX_i}(\TreeT)
    \geq 1$ if and only if $\Deg_{\VarX_i}\Par{\TreeT'} \geq 1$. Therefore, both $\TreeT$
    and $\TreeT'$ belong to $\SetTerms_n(\GeneratingSet)$ where $n$ is an integer nonsmaller
    than the greatest index of the variables appearing in $\TreeT$ and $\TreeT'$. Since
    $\Equiv$ is the reflexive, symmetric, and transitive closure of $\RewContext$, the
    result follows.
\end{proof}

\begin{Proposition} \label{prop:finite_equivalence_classes}
    Let $\CLS := (\GeneratingSet, \Rew)$ be a CLS. If $\CLS$ is hierarchical, then $\CLS$ is
    locally finite and all the $\GeneratingSet$-terms of a same connected component of
    $\RewGraph_\CLS$ have the same height.
\end{Proposition}
\begin{proof}
    Assume that $\TreeT$ and $\TreeT'$ are two $\GeneratingSet$-terms such that $\TreeT
    \RewContext \TreeT'$. From the definition of $\RewContext$ from $\Rew$ provided
    by~\eqref{equ:context_closure}, the fact that all constants of $\GeneratingSet$ are
    hierarchical implies that $\Height(\TreeT) = \Height\Par{\TreeT'}$. Moreover, since
    $\CLS$ is also nonerasing, by Lemma~\ref{lem:nonerasing_equivalence_preservation}, for
    any $\VarX_i \in \SetVariables$, $\Deg_{\VarX_i}(\TreeT) \geq 1$ if and only if
    $\Deg_{\VarX_i}\Par{\TreeT'} \geq 1$. Observe moreover that for any $n \geq 0$, the
    number of terms of the same height in $\SetTerms_n(\GeneratingSet)$ is finite. Since
    $\Equiv$ is the reflexive, symmetric, and transitive closure of $\RewContext$, and since
    two terms $\TreeT$ and $\TreeT'$ are $\Equiv$-equivalent if and only if $\TreeT$ and
    $\TreeT'$ belong to the same connected component of $\RewGraph(\CLS)$, the result
    follows.
\end{proof}

Let us use these results to state some properties of a CLS $\CLS$. By
Proposition~\ref{prop:confluence} and~\cite[Theorem 1.2.2]{BKVT03}, each
$\Equiv$-equivalence class of terms of $\CLS$ admits at most one normal form. If $\CLS$ is
hierarchical and has the poset property, it follows by
Proposition~\ref{prop:finite_equivalence_classes} that for each term $\TreeT$ of $\CLS$, the
subposet on $\Han{\TreeT}_\Equiv$ of $\Poset_\CLS$ admits exactly one maximal element. If
additionally $\CLS$ is rooted, then for each term of $\CLS$, the subposet on
$\Han{\TreeT}_\Equiv$ of $\Poset_\CLS$ admits exactly one minimal element. In this case,
there is a one-to-one correspondence between the set of the minimal terms and the set of the
maximal terms of $\CLS$ and this correspondence preserves the $\Equiv$-equivalence classes.

\subsubsection{Models} \label{subsubsec:models}
Let $\CLS := (\GeneratingSet, \Rew)$ be a CLS. A \Def{model} of $\CLS$ is an algebra over
the abstract clone defined as the quotient of $\SetTerms(\GeneratingSet)$ by the clone
congruence $\Equiv$ (see for instance~\cite{Tay93} for a general description of algebras
over abstract clones). In more concrete terms, for this particular case of abstract clone,
this is the data of a magma $(\Model, \App)$ such that $\GeneratingSet \subseteq \Model$
and, for each rule $\Combinator{X} \VarX_1 \dots \VarX_n \Rew \TreeT_{\Combinator{X}}$, the
axiom
\begin{math}
    \Par{\dots \Par{\Par{\Combinator{X} \App \VarX_1} \App \VarX_2} \dots} \App \VarX_n
    = e_{\Combinator{X}}
\end{math}
holds, where $e_{\Combinator{X}}$ is the expression having $\TreeT_{\Combinator{X}}$ as
syntax tree. For instance, a model of the CLS on the constants $\Combinator{K}$ and
$\Combinator{S}$ is a set $\Model$ containing $\Combinator{K}$ and $\Combinator{S}$, and
endowed with an operation $\App : \Model^2 \to \Model$ satisfying $\Par{\Combinator{K}
\App \VarX_1} \App \VarX_2 = \VarX_1$ and
\begin{math}
    \Par{\Par{\Combinator{S} \App \VarX_1} \App \VarX_2} \App \VarX_3
    = \Par{\VarX_1 \App \VarX_3} \App \Par{\VarX_2 \App \VarX_3}.
\end{math}
Such structures are known as \Def{combinatory algebras}~\cite{HS08}.

The \Def{trivial model}) of $\CLS$ is the magma $(\Model, \App)$ where $\Model$ is the set
$\SetTerms(\GeneratingSet)/_\Equiv$ and, for any $\Han{\TreeT_1}_\Equiv,
\Han{\TreeT_2}_\Equiv \in \Model$, $\Han{\TreeT_1}_\Equiv \App \Han{\TreeT_2}_\Equiv$ is the
$\Equiv$-equivalence class of $\TreeT_1 \TreeT_2$ where $\TreeT_1$ is any element of
$\Han{\TreeT_1}_\Equiv$ and $\TreeT_2$ is any element of $\Han{\TreeT_2}_\Equiv$. Moreover,
when $\CLS$ is nonerasing, for any $n \geq 0$, the \Def{$n$-trivial model} of $\CLS$ is the
trivial model restrained to $\Equiv$-equivalence classes of terms of
$\SetTerms_n(\GeneratingSet)$. Thanks to
Lemma~\ref{lem:nonerasing_equivalence_preservation}, the $n$-trivial model is well-defined.

When $\CLS$ is hierarchical and has the poset property, as already noticed, each
$\Equiv$-equivalence class contains exactly one maximal element. Therefore, the set
$\ModelMax$ of the maximal terms of $\SetTerms(\GeneratingSet)$ is a set of representatives
of the set of the $\Equiv$-equivalence classes and forms a model of $\CLS$ isomorphic to the
trivial model. When additionally $\CLS$ is rooted, as already noticed, each
$\Equiv$-equivalence class contains exactly one minimal element. For this reason, the set
$\ModelMin$ of the minimal terms of $\SetTerms(\GeneratingSet)$ is a set of representatives
of the set of the $\Equiv$-equivalence classes and forms a model of $\CLS$ isomorphic to the
previous ones. Here also, for any $n \geq 0$, we define the model $\ModelMax_n$ (resp.\
$\ModelMin_n$) as the restriction of $\ModelMax$ (resp.\ $\ModelMin$) on
$\SetTerms_n(\GeneratingSet)$. The interest of these two models $\ModelMax$ and $\ModelMin$
relies on considerations about algorithmic complexity for the computation of the product
$\TreeT_1 \App \TreeT_2$ where $\TreeT_1$ and $\TreeT_2$ are two terms of the models.

\section{The Mockingbird lattice} \label{sec:mockingbird_lattice}
This central section of this work concerns the study of the poset associated with the CLS
containing the basic combinator $\M$. We shall prove that this CLS has the poset property,
is rooted, and has also the lattice property by introducing new lattices on some kind of
treelike structures, called duplicative forests.

\subsection{The combinator $\M$ and its poset}
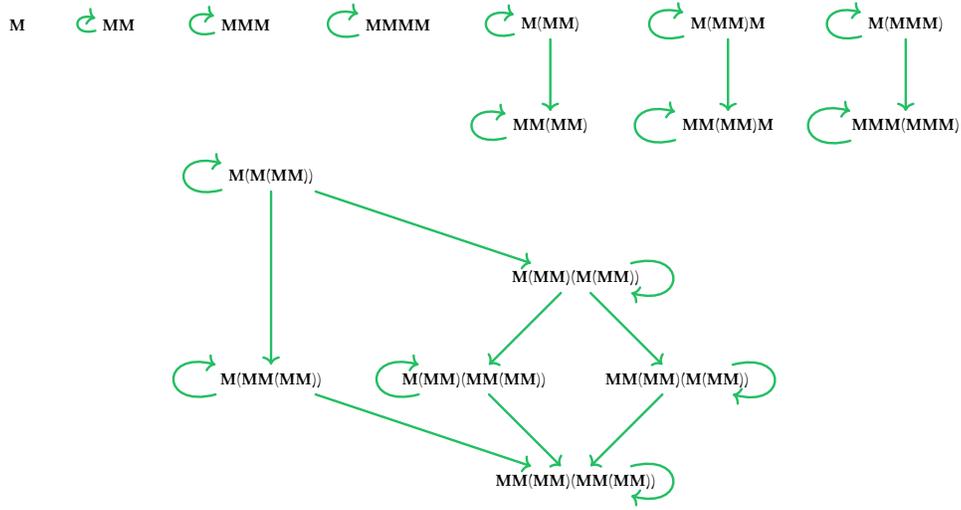
\begin{figure}[ht]
    \centering
    \scalebox{0.75}{
    \begin{tikzpicture}[Centering,xscale=0.9,yscale=0.9]
        \begin{scope}[xshift=0cm]
            \node[GraphLabeledVertex](M)at(0,0){$\M$};
        \end{scope}
        \begin{scope}[xshift=2cm]
            \node[GraphLabeledVertex](MM)at(0,0){$\M\M$};
            \draw[GraphArc](MM)edge[loop left,looseness=5](MM);
        \end{scope}
        \begin{scope}[xshift=4.5cm]
            \node[GraphLabeledVertex](MMM)at(0,0){$\M\M\M$};
            \draw[GraphArc](MMM)edge[loop left,looseness=5](MMM);
        \end{scope}
        \begin{scope}[xshift=7.5cm]
            \node[GraphLabeledVertex](MMMM)at(0,0){$\M\M\M\M$};
            \draw[GraphArc](MMMM)edge[loop left,looseness=5](MMMM);
        \end{scope}
        \begin{scope}[xshift=10.5cm]
            \node[GraphLabeledVertex](M[MM])at(0,0){$\M(\M\M)$};
            \node[GraphLabeledVertex](MM[MM])at(0,-2){$\M\M(\M\M)$};
            \draw[GraphArc](M[MM])edge[loop left,looseness=5](M[MM]);
            \draw[GraphArc](MM[MM])edge[loop left,looseness=5](MM[MM]);
            \draw[GraphArc](M[MM])--(MM[MM]);
        \end{scope}
        \begin{scope}[xshift=14cm]
            \node[GraphLabeledVertex](M[MM]M)at(0,0){$\M(\M\M)\M$};
            \node[GraphLabeledVertex](MM[MM]M)at(0,-2){$\M\M(\M\M)\M$};
            \draw[GraphArc](M[MM]M)edge[loop left,looseness=5](M[MM]M);
            \draw[GraphArc](MM[MM]M)edge[loop left,looseness=5](MM[MM]M);
            \draw[GraphArc](M[MM]M)--(MM[MM]M);
        \end{scope}
        \begin{scope}[xshift=17.5cm]
            \node[GraphLabeledVertex](M[MMM])at(0,0){$\M(\M\M\M)$};
            \node[GraphLabeledVertex](MMM[MMM])at(0,-2){$\M\M\M(\M\M\M)$};
            \draw[GraphArc](M[MMM])edge[loop left,looseness=5](M[MMM]);
            \draw[GraphArc](MMM[MMM])edge[loop left,looseness=5](MMM[MMM]);
            \draw[GraphArc](M[MMM])--(MMM[MMM]);
        \end{scope}
        \begin{scope}[xshift=7cm,yshift=-3cm]
            \node[GraphLabeledVertex](M[M[MM]])at(-2,0){$\M(\M(\M\M))$};
            \node[GraphLabeledVertex](M[MM][M[MM]])at(4,-2){$\M(\M\M)(\M(\M\M))$};
            \node[GraphLabeledVertex](M[MM][MM[MM]])at(2,-4){$\M(\M\M)(\M\M(\M\M))$};
            \node[GraphLabeledVertex](MM[MM][M[MM]])at(6,-4){$\M\M(\M\M)(\M(\M\M))$};
            \node[GraphLabeledVertex](MM[MM][MM[MM]])at(4,-6){$\M\M(\M\M)(\M\M(\M\M))$};
            \node[GraphLabeledVertex](M[MM[MM]])at(-2,-4){$\M(\M\M(\M\M))$};
            \draw[GraphArc](M[M[MM]])--(M[MM][M[MM]]);
            \draw[GraphArc](M[MM][M[MM]])--(M[MM][MM[MM]]);
            \draw[GraphArc](M[MM][M[MM]])--(MM[MM][M[MM]]);
            \draw[GraphArc](M[MM][MM[MM]])--(MM[MM][MM[MM]]);
            \draw[GraphArc](MM[MM][M[MM]])--(MM[MM][MM[MM]]);
            \draw[GraphArc](M[M[MM]])--(M[MM[MM]]);
            \draw[GraphArc](M[MM[MM]])--(MM[MM][MM[MM]]);
            \draw[GraphArc](M[M[MM]])edge[loop left,looseness=5](M[M[MM]]);
            \draw[GraphArc](M[MM][M[MM]])edge[loop right,looseness=5](M[MM][M[MM]]);
            \draw[GraphArc](M[MM][MM[MM]])edge[loop left,looseness=5](M[MM][MM[MM]]);
            \draw[GraphArc](MM[MM][M[MM]])edge[loop right,looseness=5](MM[MM][M[MM]]);
            \draw[GraphArc](MM[MM][MM[MM]])edge[loop right,looseness=5](MM[MM][MM[MM]]);
            \draw[GraphArc](M[MM[MM]])edge[loop left,looseness=5](M[MM[MM]]);
        \end{scope}
    \end{tikzpicture}}
    \caption{The fragment of the rewrite graph of $\CLS$ restrained to the terms reachable
    from closed terms of degrees $3$ or less.}
    \label{fig:rewrite_graph_l3}
\end{figure}

Let the CLS $\CLS := (\GeneratingSet, \Rew)$ such that $\GeneratingSet := \{\M\}$ and
$\TreeT_{\M} = \VarX_1 \VarX_1$. Since $\M$ is the Mockingbird basic combinator, we call
$\CLS$ the \Def{Mockingbird CLS}. Observe that $\M$ is hierarchical so that $\CLS$ satisfies
the properties stated by Proposition~\ref{prop:finite_equivalence_classes} and is in
particular locally finite. From now, we shall simply write $\RewGraph$ instead of
$\RewGraph_\CLS$. Figure~\ref{fig:rewrite_graph_l3} provides an example of a fragment
of~$\RewGraph$.

\subsubsection{First properties}
Let us present some first properties of $\CLS$. We begin by providing a recursive
description of the rewrite relation~$\RewContext$.

\begin{Lemma} \label{lem:rewrite}
    For any $\TreeT, \TreeT' \in \SetTerms(\GeneratingSet)$, we have $\TreeT \RewContext
    \TreeT'$ if and only if $\TreeT = \TreeT_1 \TreeT_2$, $\TreeT' = \TreeT_1' \TreeT_2'$
    with $\TreeT_1, \TreeT_2, \TreeT_1', \TreeT_2' \in \SetTerms(\GeneratingSet)$, and at
    least one of the following assertions holds:
    \begin{enumerate}[label=(\roman*)]
        \item $\TreeT_1 \RewContext \TreeT_1'$ and $\TreeT_2 = \TreeT_2'$;
        \item $\TreeT_2 \RewContext \TreeT_2'$ and $\TreeT_1 = \TreeT_1'$;
        \item $\TreeT_1 = \M$ and $\TreeT_2 = \TreeT_1' = \TreeT_2'$.
    \end{enumerate}
\end{Lemma}
\begin{proof}
    This is a direct consequence of the definitions of the context closure $\RewContext$
    from the rewrite relation $\Rew$ and of the term~$\TreeT_{\M}$.
\end{proof}

\begin{Proposition} \label{prop:first_graph_properties}
    The CLS $\CLS$
    \begin{enumerate}[label=(\roman*)]
        \item \label{item:first_graph_properties_1}
        is locally finite;
        \item \label{item:first_graph_properties_2}
        has the poset property;
        \item \label{item:first_graph_properties_3}
        is rooted.
    \end{enumerate}
\end{Proposition}
\begin{proof}
    First, since $\M$ is hierarchical, by Proposition~\ref{prop:finite_equivalence_classes},
    $\CLS$ is locally finite. Therefore, \ref{item:first_graph_properties_1} holds.

    Let us prove that $\CLS$ has the poset property. For this, let $\theta :
    \SetTerms(\GeneratingSet) \to \Z^2$ be the map defined for any $\TreeT \in
    \SetTerms(\GeneratingSet)$ by
    \begin{math}
        \theta(\TreeT) := \Par{\Deg(\TreeT), - \Deg_{\M}(\TreeT)}.
    \end{math}
    We denote by $\leq$ the lexicographic order on $\Z^2$ and by $\dot{+}$ the pointwise
    addition on $\Z^2$. Let us prove by structural induction on $\SetTerms(\GeneratingSet)$
    that for any $\GeneratingSet$-terms $\TreeT$ and $\TreeT'$, $\TreeT \ne \TreeT'$ and
    $\TreeT \RewContext \TreeT'$ implies $\theta(\TreeT) < \theta\Par{\TreeT'}$. By
    Lemma~\ref{lem:rewrite}, we have $\TreeT = \TreeT_1 \TreeT_2$ and $\TreeT' = \TreeT_1'
    \TreeT_2'$ for some $\GeneratingSet$-terms $\TreeT_1$, $\TreeT_2$, $\TreeT_1'$, and
    $\TreeT_2'$, and we have the three following cases. First, if $\TreeT_1 \RewContext
    \TreeT_1'$ and $\TreeT_2 = \TreeT_2'$, then
    \begin{math}
        \theta(\TreeT) = \theta\Par{\TreeT_1} \dot{+} \theta\Par{\TreeT_2} \dot{+} (1, 0)
    \end{math}
    and
    \begin{math}
        \theta\Par{\TreeT'}
        = \theta\Par{\TreeT_1'} \dot{+} \theta\Par{\TreeT_2} \dot{+} (1, 0).
    \end{math}
    Since, by induction hypothesis, $\theta\Par{\TreeT_1} < \theta\Par{\TreeT_1'}$, we have
    $\theta(\TreeT) < \theta\Par{\TreeT'}$. Second, if $\TreeT_2 \RewContext \TreeT_2'$ and
    $\TreeT_1 = \TreeT_1'$, the same arguments as the previous ones (by interchanging the
    indices $1$ and $2$ in the concerned terms) apply. Finally, if $\TreeT_1 = \M$ and
    $\TreeT_2 = \TreeT_1' = \TreeT_2'$, then
    \begin{math}
        \theta(\TreeT) = \theta\Par{\TreeT_2} \dot{+} (1, -1)
    \end{math}
    and
    \begin{math}
        \theta\Par{\TreeT'}
        = \theta\Par{\TreeT_2} \dot{+} \theta\Par{\TreeT_2} \dot{+} (1, 0).
    \end{math}
    Since $\TreeT \ne \TreeT'$, we have $\TreeT_2 \ne \M$ so that $\theta(\TreeT) <
    \theta\Par{\TreeT'}$. This implies that $\Leq$ is antisymmetric, so
    that~\ref{item:first_graph_properties_2} holds.

    To prove that $\CLS$ is rooted, we consider the rewrite relation $\Rew'$ obtained by
    inverting~$\Rew$. Thus, $\Rew'$ satisfies $\VarX_1 \VarX_1 \Rew' \M \VarX_1$. This ATRS
    $\Par{\SetTerms(\GeneratingSet), \Rew'}$ does not admit any overlapping term
    (see~\cite[Chapter 2]{BKVT03}). Therefore, it is confluent. A consequence
    of~\ref{item:first_graph_properties_1}, is that this ATRS has only finite
    $\Equiv'$-equivalence classes. By using additionally~\cite[Theorem 1.2.2]{BKVT03}, each
    $\Equiv'$-equivalence class admits exactly one term $\TreeT$ such that for any $\TreeT'
    \in \Han{\TreeT}_{\Equiv'}$, $\TreeT' \Leq' \TreeT$. This implies that in $\CLS$, each
    $\Equiv$-equivalence class admits exactly one minimal term. Therefore,
    \ref{item:first_graph_properties_3} holds.
\end{proof}

By Proposition~\ref{prop:first_graph_properties}, $\Poset_\CLS$ is a well-defined poset,
called \Def{Mockingbird poset}. From now, we shall simply write $\Poset$ instead
of~$\Poset_\CLS$. We have the following recursive description of~$\Leq$.

\begin{Lemma} \label{lem:comparison_terms}
    For any $\TreeT, \TreeT' \in \SetTerms(\GeneratingSet)$, we have $\TreeT \Leq \TreeT'$
    if and only if at least one of the following assertions holds:
    \begin{enumerate}[label=(\roman*)]
        \item $\TreeT = \M = \TreeT'$;
        \item $\TreeT = \VarX_i = \TreeT'$ for an $\VarX_i \in \SetVariables$;
        \item $\TreeT = \TreeT_1 \TreeT_2$, $\TreeT' = \TreeT_1' \TreeT_2'$, $\TreeT_1 \Leq
        \TreeT_1'$, and $\TreeT_2 \Leq \TreeT_2'$ where $\TreeT_1, \TreeT_1', \TreeT_2,
        \TreeT_2' \in \SetTerms(\GeneratingSet)$;
        \item $\TreeT = \M \TreeT_2$, $\TreeT' = \TreeT_1' \TreeT_2'$, $\TreeT_2 \Leq
        \TreeT_1'$, and $\TreeT_2 \Leq \TreeT_2'$ where $\TreeT_1', \TreeT_2, \TreeT_2' \in
        \SetTerms(\GeneratingSet)$.
    \end{enumerate}
\end{Lemma}
\begin{proof}
    This is a direct consequence of Lemma~\ref{lem:rewrite} and of the fact that $\Leq$ is
    the reflexive and transitive closure of~$\RewContext$.
\end{proof}

\begin{Proposition} \label{prop:minimal_maximal}
    Let $\TreeT \in \SetTerms(\GeneratingSet)$.
    \begin{enumerate}[label=(\roman*)]
        \item \label{item:minimal_maximal_1}
        The term $\TreeT$ is a maximal element of $\Poset$ if and only if $\TreeT$ avoids
        $\M \Par{\VarX_1 \VarX_2}$ and suffix-avoids $\M \VarX_1$.
        \item \label{item:minimal_maximal_2}
        The term $\TreeT$ is a minimal element of $\Poset$ if and only if $\TreeT$ avoids
        $\Par{\VarX_1 \VarX_2} \Par{\VarX_1 \VarX_2}$ and suffix-avoids $\VarX_1 \VarX_1$.
    \end{enumerate}
\end{Proposition}
\begin{proof}
    The term $\TreeT$ is maximal in $\Poset$ if and only if $\TreeT \RewContext \TreeT'$
    implies that $\TreeT' = \TreeT$. By Lemma~\ref{lem:rewrite}, this is equivalent to the
    fact $\TreeT$ is maximal in $\Poset$ if and only if when $\TreeT$ has an internal node
    having $\M$ as left subtree, this internal node admits necessarily $\M$ as right
    subtree. Therefore, \ref{item:minimal_maximal_1} holds.

    The term $\TreeT$ is minimal in $\Poset$ if and only if $\TreeT' \RewContext \TreeT$
    implies that $\TreeT' = \TreeT$. By Lemma~\ref{lem:rewrite}, this is equivalent to the
    fact $\TreeT$ is minimal in $\Poset$ if and only if when $\TreeT$ has an internal node
    such that its two children are equal, these two children are necessarily equal to $\M$.
    Therefore, \ref{item:minimal_maximal_2} holds.
\end{proof}

\subsubsection{Some models}
We begin this discussion about models of $\CLS$ by a very simple observation: any idempotent
monoid $(\Monoid, \cdot, e)$ is a model of $\CLS$. Indeed, by identifying $\M$ with $e$ (the
unit of~$\Monoid$) and by identifying $\App$ with~$\cdot$, we have for any $x \in \Monoid$,
\begin{math}
    \M \App x = e \cdot x = x = x \cdot x = x \App x.
\end{math}

Besides, since $\CLS$ is hierarchical and has, by
Proposition~\ref{prop:first_graph_properties}, the poset property, this CLS admits the model
$\ModelMax$ described in Section~\ref{subsubsec:models}. More specifically, this model is
such that $\ModelMax$ is the set of the maximal elements of $\Poset$, and, by
Proposition~\ref{prop:minimal_maximal}, for any $\TreeT_1, \TreeT_2 \in \ModelMax$,
\begin{equation}
    \TreeT_1 \App \TreeT_2 =
    \begin{cases}
        \TreeT_2 \TreeT_2 & \mbox{ if } \TreeT_1 = \M, \\
        \TreeT_1 \TreeT_2 & \mbox{otherwise}.
    \end{cases}
\end{equation}
Since by Proposition~\ref{prop:first_graph_properties}, $\CLS$ is also rooted, $\CLS$ admits
the model $\ModelMin$ described in Section~\ref{subsubsec:models}. This model is such that
$\ModelMin$ is the set of the minimal elements of $\Poset$, and, by
Proposition~\ref{prop:minimal_maximal}, for any $\TreeT_1, \TreeT_2 \in \ModelMin$,
\begin{equation}
    \TreeT_1 \App \TreeT_2 =
    \begin{cases}
        \M \TreeT_2 & \mbox{ if } \TreeT_1 = \TreeT_2, \\
        \TreeT_1 \TreeT_2 & \mbox{otherwise}.
    \end{cases}
\end{equation}
By representing $\GeneratingSet$-terms directly has binary trees, the computations in
$\ModelMax$ have a better time complexity than in $\ModelMin$ because to compute $\TreeT_1
\App \TreeT_2$ in $\ModelMin$, we need to decide if $\TreeT_1$ and $\TreeT_2$ are equal.
Nevertheless, in return, the term $\TreeT_1 \App \TreeT_2$ produced by a computation in
$\ModelMin$ has a number of internal nodes equal as or smaller than the analogous
computation in~$\ModelMax$.

\subsection{Duplicative lattices}
We introduce here duplicative forests and a partial order relation on these objects. We show
that this partial relation endow all intervals of this poset with lattice structures.

\subsubsection{Duplicative forests}
A \Def{duplicative tree} is a planar rooted tree such that each internal node is either a
\Def{black node} $\BlackNode$ or a \Def{white node} $\WhiteNode$. A \Def{duplicative forest}
is a word $\ForestF = \ForestF(1) \dots \ForestF(\ell)$, $\ell \geq 0$, of duplicative
trees. In particular, the empty word $\epsilon$ is the \Def{empty forest}. We denote by
$\SetDuplicativeTrees$ (resp.\ $\SetDuplicativeForests$) the set of such trees (resp.\
forests). The \Def{height} $\Height(\ForestF)$ of $\ForestF$ is the number of internal nodes
in a longest path following edges connecting a node to one of its child. In the sequel, we
shall write expressions using some occurrences of $\AnyNode$: each such expression denotes
the two expressions obtained by replacing simultaneously all $\AnyNode$ either by
$\WhiteNode$ or by~$\BlackNode$. The \Def{grafting product} is the unary operation
$\AnyNode$ on $\SetDuplicativeForests$ such that for any $\ForestF \in
\SetDuplicativeForests$, $\AnyNode(\ForestF)$ is the duplicative tree obtained by grafting
the roots of the duplicative trees of $\ForestF$ on a common root node~$\AnyNode$. The
\Def{concatenation product} is the binary operation $\ConcatenateForests$ on
$\SetDuplicativeForests$ such that for any $\ForestF_1, \ForestF_2 \in
\SetDuplicativeForests$, $\ForestF_1 \ConcatenateForests \ForestF_2$ is the duplicative
forest made of the trees of $\ForestF_1$ and then of the trees of~$\ForestF_2$.

\subsubsection{A poset on duplicative forests}
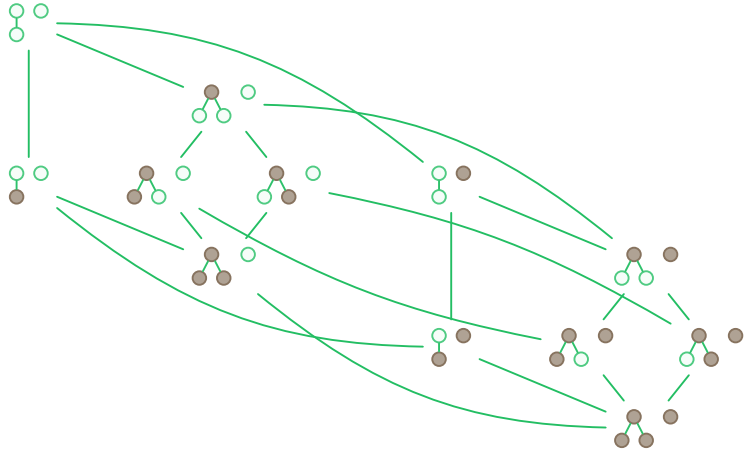
\begin{figure}[ht]
    \centering
    \scalebox{.9}{
    \begin{tikzpicture}[Centering,xscale=0.24,yscale=0.3]
        \tikzstyle{GraphEdge2}=[GraphEdge,thick]
        \node(WWW)at(-12,0){
            \begin{tikzpicture}[Centering,xscale=0.18,yscale=0.26]
                \node[WhiteNode](1)at(0.00,-2.67){};
                \node[WhiteNode](3)at(2.00,-1.33){};
                \node[WhiteNode](0)at(0.00,-1.33){};
                \draw[Edge](1)--(0);
            \end{tikzpicture}};
        \node(WBW)at(-12,-8){
            \begin{tikzpicture}[Centering,xscale=0.18,yscale=0.26]
                \node[BlackNode](1)at(0.00,-2.67){};
                \node[WhiteNode](3)at(2.00,-1.33){};
                \node[WhiteNode](0)at(0.00,-1.33){};
                \draw[Edge](1)--(0);
            \end{tikzpicture}};
        \node(BWWW)at(0,-4){
            \begin{tikzpicture}[Centering,xscale=0.18,yscale=0.26]
                \node[WhiteNode](0)at(0.00,-2.67){};
                \node[WhiteNode](2)at(2.00,-2.67){};
                \node[WhiteNode](4)at(4.00,-1.33){};
                \node[BlackNode](1)at(1.00,-1.33){};
                \draw[Edge](0)--(1);
                \draw[Edge](2)--(1);
            \end{tikzpicture}};
        \node(BBWW)at(-4,-8){
            \begin{tikzpicture}[Centering,xscale=0.18,yscale=0.26]
                \node[BlackNode](0)at(0.00,-2.67){};
                \node[WhiteNode](2)at(2.00,-2.67){};
                \node[WhiteNode](4)at(4.00,-1.33){};
                \node[BlackNode](1)at(1.00,-1.33){};
                \draw[Edge](0)--(1);
                \draw[Edge](2)--(1);
            \end{tikzpicture}};
        \node(BWBW)at(4,-8){
            \begin{tikzpicture}[Centering,xscale=0.18,yscale=0.26]
                \node[WhiteNode](0)at(0.00,-2.67){};
                \node[BlackNode](2)at(2.00,-2.67){};
                \node[WhiteNode](4)at(4.00,-1.33){};
                \node[BlackNode](1)at(1.00,-1.33){};
                \draw[Edge](0)--(1);
                \draw[Edge](2)--(1);
            \end{tikzpicture}};
        \node(BBBW)at(0,-12){
            \begin{tikzpicture}[Centering,xscale=0.18,yscale=0.26]
                \node[BlackNode](0)at(0.00,-2.67){};
                \node[BlackNode](2)at(2.00,-2.67){};
                \node[WhiteNode](4)at(4.00,-1.33){};
                \node[BlackNode](1)at(1.00,-1.33){};
                \draw[Edge](0)--(1);
                \draw[Edge](2)--(1);
            \end{tikzpicture}};
        \draw[GraphEdge2](WWW)--(WBW);
        \draw[GraphEdge2](WWW)--(BWWW);
        \draw[GraphEdge2](WBW)--(BBBW);
        \draw[GraphEdge2](BWWW)--(BBWW);
        \draw[GraphEdge2](BWWW)--(BWBW);
        \draw[GraphEdge2](BBWW)--(BBBW);
        \draw[GraphEdge2](BWBW)--(BBBW);
        \node(WWB)at(14,-8){
            \begin{tikzpicture}[Centering,xscale=0.18,yscale=0.26]
                \node[WhiteNode](1)at(0.00,-2.67){};
                \node[BlackNode](3)at(2.00,-1.33){};
                \node[WhiteNode](0)at(0.00,-1.33){};
                \draw[Edge](1)--(0);
            \end{tikzpicture}};
        \node(WBB)at(14,-16){
            \begin{tikzpicture}[Centering,xscale=0.18,yscale=0.26]
                \node[BlackNode](1)at(0.00,-2.67){};
                \node[BlackNode](3)at(2.00,-1.33){};
                \node[WhiteNode](0)at(0.00,-1.33){};
                \draw[Edge](1)--(0);
            \end{tikzpicture}};
        \node(BWWB)at(26,-12){
            \begin{tikzpicture}[Centering,xscale=0.18,yscale=0.26]
                \node[WhiteNode](0)at(0.00,-2.67){};
                \node[WhiteNode](2)at(2.00,-2.67){};
                \node[BlackNode](4)at(4.00,-1.33){};
                \node[BlackNode](1)at(1.00,-1.33){};
                \draw[Edge](0)--(1);
                \draw[Edge](2)--(1);
            \end{tikzpicture}};
        \node(BBWB)at(22,-16){
            \begin{tikzpicture}[Centering,xscale=0.18,yscale=0.26]
                \node[BlackNode](0)at(0.00,-2.67){};
                \node[WhiteNode](2)at(2.00,-2.67){};
                \node[BlackNode](4)at(4.00,-1.33){};
                \node[BlackNode](1)at(1.00,-1.33){};
                \draw[Edge](0)--(1);
                \draw[Edge](2)--(1);
            \end{tikzpicture}};
        \node(BWBB)at(30,-16){
            \begin{tikzpicture}[Centering,xscale=0.18,yscale=0.26]
                \node[WhiteNode](0)at(0.00,-2.67){};
                \node[BlackNode](2)at(2.00,-2.67){};
                \node[BlackNode](4)at(4.00,-1.33){};
                \node[BlackNode](1)at(1.00,-1.33){};
                \draw[Edge](0)--(1);
                \draw[Edge](2)--(1);
            \end{tikzpicture}};
        \node(BBBB)at(26,-20){
            \begin{tikzpicture}[Centering,xscale=0.18,yscale=0.26]
                \node[BlackNode](0)at(0.00,-2.67){};
                \node[BlackNode](2)at(2.00,-2.67){};
                \node[BlackNode](4)at(4.00,-1.33){};
                \node[BlackNode](1)at(1.00,-1.33){};
                \draw[Edge](0)--(1);
                \draw[Edge](2)--(1);
            \end{tikzpicture}};
        \draw[GraphEdge2](WWB)--(WBB);
        \draw[GraphEdge2](WWB)--(BWWB);
        \draw[GraphEdge2](WBB)--(BBBB);
        \draw[GraphEdge2](BWWB)--(BBWB);
        \draw[GraphEdge2](BWWB)--(BWBB);
        \draw[GraphEdge2](BBWB)--(BBBB);
        \draw[GraphEdge2](BWBB)--(BBBB);
        \draw[GraphEdge2](WWW)edge[bend left=16](WWB);
        \draw[GraphEdge2](WBW)edge[bend right=16](WBB);
        \draw[GraphEdge2](BWWW)edge[bend left=16](BWWB);
        \draw[GraphEdge2](BBWW)edge[bend right=8](BBWB);
        \draw[GraphEdge2](BWBW)edge[bend left=8](BWBB);
        \draw[GraphEdge2](BBBW)edge[bend right=16](BBBB);
    \end{tikzpicture}}
    \caption{The Hasse diagram of a maximal interval of the duplicative forest poset.}
    \label{fig:example_poset_duplicative_forests}
\end{figure}

Let $\RewDuplicative$ be the binary relation on $\SetDuplicativeForests$ defined recursively
as follows. For any $\ForestF, \ForestF' \in \SetDuplicativeForests$, we have $\ForestF
\RewDuplicative \ForestF'$ if $\ForestF$ and $\ForestF'$ have the same length $\ell \geq 1$
and
\begin{itemize}
    \item either $\ell = 1$, $\ForestF = \WhiteNode(\ForestG)$ and $\ForestF' =
    \BlackNode(\ForestG \ConcatenateForests \ForestG)$ where $\ForestG \in
    \SetDuplicativeForests$;

    \item or $\ell = 1$, $\ForestF = \AnyNode(\ForestG)$, $\ForestF' = \AnyNode
    \Par{\ForestG'}$ where $\ForestG, \ForestG' \in \SetDuplicativeForests$ and $\ForestG
    \RewDuplicative \ForestG'$;

    \item or $\ell \geq 2$ and there is a $j \in [\ell]$ such that $\ForestF(j)
    \RewDuplicative \ForestF'(j)$, and for all $i \in [\ell] \setminus \{j\}$, $\ForestF(i)
    = \ForestF'(i)$.
\end{itemize}
We observe that $\ForestF \RewDuplicative \ForestF'$ if and only if $\ForestF'$ can be
obtained from $\ForestF$ by selecting a white node of $\ForestF$, by turning it into black,
and by duplicating its sequence of descendants. For instance, we have
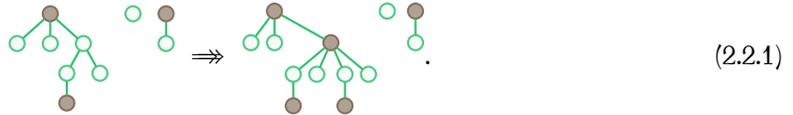
\begin{equation}
    \begin{tikzpicture}[Centering,xscale=0.22,yscale=0.18]
        \node[WhiteNode](0)at(-1.00,-4.40){};
        \node[WhiteNode](10)at(8.00,-4.40){};
        \node[WhiteNode](2)at(1.00,-4.40){};
        \node[BlackNode](4)at(2.00,-8.80){};
        \node[WhiteNode](6)at(4.00,-6.60){};
        \node[WhiteNode](8)at(6.00,-2.20){};
        \node[BlackNode](1)at(1.00,-2.20){};
        \node[WhiteNode](3)at(2.00,-6.60){};
        \node[WhiteNode](5)at(3.00,-4.40){};
        \node[BlackNode](9)at(8.00,-2.20){};
        \draw[Edge](0)--(1);
        \draw[Edge](10)--(9);
        \draw[Edge](2)--(1);
        \draw[Edge](3)--(5);
        \draw[Edge](4)--(3);
        \draw[Edge](5)--(1);
        \draw[Edge](6)--(5);
    \end{tikzpicture}
    \enspace \RewDuplicative \enspace
    \begin{tikzpicture}[Centering,xscale=0.25,yscale=0.15]
        \node[WhiteNode](0)at(-0.25,-5.60){};
        \node[WhiteNode](11)at(7.00,-2.80){};
        \node[WhiteNode](13)at(8.50,-5.60){};
        \node[WhiteNode](2)at(1.00,-5.60){};
        \node[BlackNode](4)at(2.00,-11.20){};
        \node[WhiteNode](5)at(3.25,-8.40){};
        \node[BlackNode](8)at(4.75,-11.20){};
        \node[WhiteNode](9)at(6.00,-8.40){};
        \node[BlackNode](1)at(1.00,-2.80){};
        \node[BlackNode](12)at(8.50,-2.80){};
        \node[WhiteNode](3)at(2.00,-8.40){};
        \node[BlackNode](6)at(4.00,-5.60){};
        \node[WhiteNode](7)at(4.75,-8.40){};
        \draw[Edge](0)--(1);
        \draw[Edge](13)--(12);
        \draw[Edge](2)--(1);
        \draw[Edge](3)--(6);
        \draw[Edge](4)--(3);
        \draw[Edge](5)--(6);
        \draw[Edge](6)--(1);
        \draw[Edge](7)--(6);
        \draw[Edge](8)--(7);
        \draw[Edge](9)--(6);
    \end{tikzpicture}.
\end{equation}
Observe also that in this case, there are more black nodes in $\ForestF'$ than in
$\ForestF$. Hence, the reflexive and transitive closure $\LeqDuplicative$ of
$\RewDuplicative$ is antisymmetric so that $(\SetDuplicativeForests, \LeqDuplicative)$ is a
poset. We call this poset the \Def{duplicative forest poset}. For any $\ForestF \in
\SetDuplicativeForests$, we denote by $\SetDuplicativeForests(\ForestF)$ the subposet of the
duplicative forest poset on the set $\Bra{\ForestF' \in \SetDuplicativeForests : \ForestF
\LeqDuplicative \ForestF'}$. Figure~\ref{fig:example_poset_duplicative_forests} shows the
Hasse diagram of the poset $\SetDuplicativeForests(\ForestF)$ where $\ForestF$ is a certain
duplicative forest.
According to this Hasse diagram, the duplicative forest poset is not graded. We have the
following recursive description of~$\LeqDuplicative$.

\begin{Lemma} \label{lem:comparison_duplicative_forests}
    For any two duplicative forests $\ForestF$ and $\ForestF'$, we have $\ForestF
    \LeqDuplicative \ForestF'$ if and only if $\ForestF$ and $\ForestF'$ have the same
    length $\ell \geq 0$ and one of the following assertions holds:
    \begin{enumerate}[label=(\roman*)]
        \item $\ell = 0$;
        \item $\ell = 1$, $\ForestF = \WhiteNode(\ForestG)$ and $\ForestF' =
        \BlackNode\Par{\ForestG' \ConcatenateForests \ForestG''}$ where $\ForestG,
        \ForestG', \ForestG'' \in \SetDuplicativeForests$, $\ForestG \LeqDuplicative
        \ForestG'$, and $\ForestG \LeqDuplicative \ForestG''$;
        \item $\ell = 1$, $\ForestF = \AnyNode(\ForestG)$, $\ForestF' = \AnyNode
        \Par{\ForestG'}$ where $\ForestG, \ForestG' \in \SetDuplicativeForests$ and
        $\ForestG \LeqDuplicative \ForestG'$;
        \item $\ell \geq 2$ and $\ForestF(i) \LeqDuplicative \ForestF'(i)$ for all $i \in
        [\ell]$.
    \end{enumerate}
\end{Lemma}
\begin{proof}
    This is a direct consequence of the definitions of $\RewDuplicative$ and of
    $\LeqDuplicative$.
\end{proof}

A first consequence of Lemma~\ref{lem:comparison_duplicative_forests} is that for any
$\ForestF, \ForestF' \in \SetDuplicativeForests$, $\ForestF \LeqDuplicative \ForestF'$
implies $\Height(\ForestF) = \Height\Par{\ForestF'}$.

\subsubsection{Lattices on duplicative forests}
\label{subsubsec:lattices_duplicative_forests}
We use in the sequel the usual notions and notations about lattices (see for
instance~\cite{Sta11b}). Let $\Meet$ and $\JJoin$ be the two binary, commutative, and
associative partial operations on $\SetDuplicativeForests$ defined recursively, for any
$\ell \geq 0$, $\ForestF_1, \dots, \ForestF_\ell \in \SetDuplicativeTrees$, $\ForestF_1',
\dots, \ForestF_\ell' \in \SetDuplicativeTrees$, and $\ForestF, \ForestF', \ForestF'' \in
\SetDuplicativeForests$, by
\begin{subequations}
\begin{equation}
    \ForestF_1 \ConcatenateForests \dots \ConcatenateForests \ForestF_\ell
    \Meet
    \ForestF_1' \ConcatenateForests \dots \ConcatenateForests \ForestF_\ell'
    := \Par{\ForestF_1 \Meet \ForestF_1'}
    \ConcatenateForests \dots \ConcatenateForests
    \Par{\ForestF_\ell \Meet \ForestF_\ell'},
\end{equation}
\begin{equation}
    \AnyNode(\ForestF) \Meet \AnyNode\Par{\ForestF'}
    := \AnyNode\Par{\ForestF \Meet \ForestF'},
\end{equation}
\begin{equation}
    \WhiteNode(\ForestF) \Meet \BlackNode\Par{\ForestF' \ConcatenateForests \ForestF''}
    := \WhiteNode\Par{\ForestF \Meet \ForestF' \Meet \ForestF''}
\end{equation}
\end{subequations}
and
\begin{subequations}
\begin{equation}
    \ForestF_1 \ConcatenateForests \dots \ConcatenateForests \ForestF_\ell
    \JJoin
    \ForestF_1' \ConcatenateForests \dots \ConcatenateForests \ForestF_\ell'
    := \Par{\ForestF_1 \JJoin \ForestF_1'}
    \ConcatenateForests \dots \ConcatenateForests
    \Par{\ForestF_\ell \JJoin \ForestF_\ell'},
\end{equation}
\begin{equation}
    \AnyNode(\ForestF) \JJoin \AnyNode\Par{\ForestF'}
    := \AnyNode\Par{\ForestF \JJoin \ForestF'},
\end{equation}
\begin{equation}
    \WhiteNode(\ForestF) \JJoin \BlackNode\Par{\ForestF' \ConcatenateForests \ForestF''}
    := \BlackNode\Par{\Par{\ForestF \JJoin \ForestF'}
    \ConcatenateForests \Par{\ForestF \JJoin \ForestF''}}.
\end{equation}
\end{subequations}

\begin{Proposition} \label{prop:lattice_duplicative_forests}
    Given a duplicative forest $\ForestF$, the poset $\SetDuplicativeForests(\ForestF)$ is a
    lattice for the operations~$\Meet$ and~$\JJoin$.
\end{Proposition}
\begin{proof}
    This follows by structural induction on $\ForestF$ by establishing the fact that for any
    $\ForestF', \ForestF'' \in \SetDuplicativeForests(\ForestF)$, the duplicative forest
    $\ForestF' \Meet \ForestF''$ (resp.\ $\ForestF' \JJoin \ForestF''$) is well-defined and
    is the meet (resp.\ the join) of $\ForestF'$ and $\ForestF''$. This uses
    Lemma~\ref{lem:comparison_duplicative_forests} in order to describe $\ForestF'$ and
    $\ForestF''$ from $\ForestF$ knowing that $\ForestF \LeqDuplicative \ForestF'$
    and~$\ForestF \LeqDuplicative \ForestF''$.
\end{proof}

By Proposition~\ref{prop:lattice_duplicative_forests}, for any $\ForestF \in
\SetDuplicativeForests$, each poset $\SetDuplicativeForests(\ForestF)$ is a lattice. We call
it the \Def{duplicative forest lattice} of~$\ForestF$.

Let $\Pruned : \SetDuplicativeForests \to \SetDuplicativeForests$ be the map defined
recursively, for any $\ell \geq 0$, $\ForestF_1, \dots, \ForestF_\ell \in
\SetDuplicativeForests$, and $\ForestF \in \SetDuplicativeForests$, by
\begin{subequations}
\begin{equation}
    \Pruned\Par{\ForestF_1 \ConcatenateForests \dots \ConcatenateForests \ForestF_\ell}
    := \Pruned\Par{\ForestF_1}
    \ConcatenateForests \dots \ConcatenateForests
    \Pruned\Par{\ForestF_\ell},
\end{equation}
\begin{equation}
    \Pruned\Par{\WhiteNode(\ForestF)} := \WhiteNode(\Pruned(\ForestF)),
\end{equation}
\begin{equation}
    \Pruned\Par{\BlackNode(\ForestF)} := \Pruned(\ForestF).
\end{equation}
\end{subequations}
The duplicative forest $\Pruned(\ForestF)$ is the \Def{pruning} of $\ForestF$. This forest
has by construction no occurrence of any $\BlackNode$. For instance,
\begin{equation}
    \begin{tikzpicture}[Centering,xscale=0.25,yscale=0.18]
        \node[WhiteNode](1)at(0.00,-6.75){};
        \node[WhiteNode](3)at(2.00,-6.75){};
        \node[BlackNode](5)at(4.00,-6.75){};
        \node[BlackNode](7)at(5.25,-2.25){};
        \node[WhiteNode](8)at(6.50,-2.25){};
        \node[BlackNode](0)at(0.00,-4.50){};
        \node[WhiteNode](2)at(1.00,-2.25){};
        \node[WhiteNode](4)at(3.00,-4.50){};
        \draw[Edge](0)--(2);
        \draw[Edge](1)--(0);
        \draw[Edge](3)--(4);
        \draw[Edge](4)--(2);
        \draw[Edge](5)--(4);
    \end{tikzpicture}
    \enspace \xmapsto{\Pruned} \enspace
    \begin{tikzpicture}[Centering,xscale=0.25,yscale=0.22]
        \node[WhiteNode](0)at(0.00,-3.00){};
        \node[WhiteNode](3)at(2.00,-4.50){};
        \node[WhiteNode](5)at(4.00,-1.50){};
        \node[WhiteNode](1)at(1.00,-1.50){};
        \node[WhiteNode](2)at(2.00,-3.00){};
        \draw[Edge](0)--(1);
        \draw[Edge](2)--(1);
        \draw[Edge](3)--(2);
    \end{tikzpicture}.
\end{equation}

\begin{Lemma} \label{lem:isomorphic_pruned_forests}
    For any $\ForestF \in \SetDuplicativeForests$, the posets
    $\SetDuplicativeForests(\ForestF)$ and $\SetDuplicativeForests(\Pruned(\ForestF))$ are
    isomorphic.
\end{Lemma}
\begin{proof}
    This follows by structural induction on $\ForestF$ by establishing the fact that
    $\Pruned$ is a one-to-one correspondence between the sets
    $\SetDuplicativeForests(\ForestF)$ and $\SetDuplicativeForests(\Pruned(\ForestF))$ and
    that $\Pruned$ is an order isomorphism. This uses
    Lemma~\ref{lem:comparison_duplicative_forests} in order to have a recursive description
    of the order relation~$\LeqDuplicative$.
\end{proof}

For any $d \geq 0$, the \Def{$d$-ladder} is the duplicative forest $\Ladder_d$ defined
recursively by $\Ladder_0 := \epsilon$ and, for any $d \geq 1$, by $\Ladder_d := \WhiteNode
\Par{\Ladder_{d - 1}}$. Let also
\begin{equation}
    \GreaterLadders := \bigcup_{d \geq 0} \SetDuplicativeForests\Par{\Ladder_d}.
\end{equation}
By definition, $\GreaterLadders$ is the set of the duplicative forests that are equal as or
greater than a $d$-ladder for a $d \geq 0$.

\begin{Lemma} \label{lem:smaller_ladder_pruned}
    For any $\ForestF \in \SetDuplicativeForests$, there exists $\ForestG \in
    \GreaterLadders$ such that $\Pruned(\ForestF) = \Pruned(\ForestG)$.
\end{Lemma}
\begin{proof}
    This follows by structural induction on $\ForestF$ by showing that there is a $d \geq 0$
    such that $\Ladder_d \LeqDuplicative \ForestG$ and $\Pruned(\ForestF) =
    \Pruned(\ForestG)$. This uses Lemma~\ref{lem:comparison_duplicative_forests}.
\end{proof}

An important consequence of Lemmas~\ref{lem:isomorphic_pruned_forests}
and~\ref{lem:smaller_ladder_pruned} is that for any $\ForestF \in \SetDuplicativeForests$,
there exists $d \geq 0$ such that $\SetDuplicativeForests(\ForestF)$ is isomorphic as a
poset to a maximal interval of $\SetDuplicativeForests \Par{\Ladder_d}$.

\subsection{Mockingbird lattices}
We show here that each subposet $\Poset(\TreeT)$, $\TreeT \in \SetTerms(\GeneratingSet)$, of
$\Poset$ is a lattice. This is based on a poset isomorphism between $\Poset(\TreeT)$ and an
interval of a lattice of duplicative forests. We also define for each $d \geq 0$ the
Mockingbird lattice of order $d$ as a particular maximal interval of~$\Poset$.

\subsubsection{From terms to duplicative forests}
Let $\ToForest : \SetTerms(\GeneratingSet) \to \SetDuplicativeForests$ be the map defined
recursively, for any $\VarX_i \in \SetVariables$ and $\TreeT, \TreeT' \in
\SetTerms(\GeneratingSet)$, by
\begin{subequations}
\begin{equation}
    \ToForest\Par{\VarX_i} := \epsilon,
\end{equation}
\begin{equation}
    \ToForest(\M) := \epsilon,
\end{equation}
\begin{equation}
    \ToForest\Par{\TreeT \TreeT'} :=
    \begin{cases}
        \WhiteNode\Par{\ToForest\Par{\TreeT'}}
            & \mbox{ if } \TreeT = \M \mbox{ and } \TreeT' \ne \M, \\
        \BlackNode\Par{\ToForest(\TreeT) \ConcatenateForests \ToForest\Par{\TreeT'}}
            & \mbox{otherwise}.
    \end{cases}
\end{equation}
\end{subequations}
For instance,
\begin{equation}
    \begin{tikzpicture}[Centering,xscale=0.22,yscale=0.13]
        \node[LeafST](0)at(0.00,-11.57){$\M$};
        \node[LeafST](10)at(10.00,-23.14){$\VarX_1$};
        \node[LeafST](12)at(12.00,-11.57){$\VarX_1$};
        \node[LeafST](14)at(14.00,-15.43){$\M$};
        \node[LeafST](16)at(16.00,-15.43){$\M$};
        \node[LeafST](18)at(18.00,-11.57){$\M$};
        \node[LeafST](2)at(2.00,-19.29){$\M$};
        \node[LeafST](20)at(20.00,-15.43){$\M$};
        \node[LeafST](22)at(22.00,-15.43){$\M$};
        \node[LeafST](24)at(24.00,-11.57){$\VarX_3$};
        \node[LeafST](26)at(26.00,-11.57){$\VarX_2$};
        \node[LeafST](4)at(4.00,-19.29){$\VarX_2$};
        \node[LeafST](6)at(6.00,-19.29){$\M$};
        \node[LeafST](8)at(8.00,-23.14){$\M$};
        \node[NodeST](1)at(1.00,-7.71){$\App$};
        \node[NodeST](11)at(11.00,-3.86){$\App$};
        \node[NodeST](13)at(13.00,-7.71){$\App$};
        \node[NodeST](15)at(15.00,-11.57){$\App$};
        \node[NodeST](17)at(17.00,0.00){$\App$};
        \node[NodeST](19)at(19.00,-7.71){$\App$};
        \node[NodeST](21)at(21.00,-11.57){$\App$};
        \node[NodeST](23)at(23.00,-3.86){$\App$};
        \node[NodeST](25)at(25.00,-7.71){$\App$};
        \node[NodeST](3)at(3.00,-15.43){$\App$};
        \node[NodeST](5)at(5.00,-11.57){$\App$};
        \node[NodeST](7)at(7.00,-15.43){$\App$};
        \node[NodeST](9)at(9.00,-19.29){$\App$};
        \draw[Edge](0)--(1);
        \draw[Edge](1)--(11);
        \draw[Edge](10)--(9);
        \draw[Edge](11)--(17);
        \draw[Edge](12)--(13);
        \draw[Edge](13)--(11);
        \draw[Edge](14)--(15);
        \draw[Edge](15)--(13);
        \draw[Edge](16)--(15);
        \draw[Edge](18)--(19);
        \draw[Edge](19)--(23);
        \draw[Edge](2)--(3);
        \draw[Edge](20)--(21);
        \draw[Edge](21)--(19);
        \draw[Edge](22)--(21);
        \draw[Edge](23)--(17);
        \draw[Edge](24)--(25);
        \draw[Edge](25)--(23);
        \draw[Edge](26)--(25);
        \draw[Edge](3)--(5);
        \draw[Edge](4)--(3);
        \draw[Edge](5)--(1);
        \draw[Edge](6)--(7);
        \draw[Edge](7)--(5);
        \draw[Edge](8)--(9);
        \draw[Edge](9)--(7);
        \node(r)at(17.00,2.89){};
        \draw[Edge](r)--(17);
    \end{tikzpicture}
    \quad \xmapsto{\ToForest} \quad
    \begin{tikzpicture}[Centering,xscale=0.11,yscale=0.09]
        \node[WhiteNode](1)at(1.00,-7.71){};
        \node[BlackNode](11)at(11.00,-3.86){};
        \node[BlackNode](13)at(13.00,-7.71){};
        \node[BlackNode](15)at(15.00,-11.57){};
        \node[BlackNode](17)at(17.00,0.00){};
        \node[WhiteNode](19)at(19.00,-7.71){};
        \node[BlackNode](21)at(21.00,-11.57){};
        \node[BlackNode](23)at(23.00,-3.86){};
        \node[BlackNode](25)at(25.00,-7.71){};
        \node[WhiteNode](3)at(3.00,-15.43){};
        \node[BlackNode](5)at(5.00,-11.57){};
        \node[WhiteNode](7)at(7.00,-15.43){};
        \node[WhiteNode](9)at(9.00,-19.29){};
        \draw[Edge](1)--(11);
        \draw[Edge](11)--(17);
        \draw[Edge](13)--(11);
        \draw[Edge](15)--(13);
        \draw[Edge](19)--(23);
        \draw[Edge](21)--(19);
        \draw[Edge](23)--(17);
        \draw[Edge](25)--(23);
        \draw[Edge](3)--(5);
        \draw[Edge](5)--(1);
        \draw[Edge](7)--(5);
        \draw[Edge](9)--(7);
    \end{tikzpicture}.
\end{equation}
Intuitively, $\ToForest(\TreeT)$ is the duplicative forest obtained from $\TreeT$ by the
following process: replace each internal node of $\TreeT$ having $\M$ as left child and
having a right child different from $\M$ by a~$\WhiteNode$, replace each other internal node
of $\TreeT$ by a~$\BlackNode$, and remove all the leaves. Immediately from the definition,
we observe that this map is not injective. Moreover, again directly from the definition, any
duplicative forest $\ForestF$ in the image of $\ToForest$ is such that each white node of
$\ForestF$ has no more than one child.

\begin{Lemma} \label{lem:to_forest_injection}
    For any $\TreeT \in \SetTerms(\GeneratingSet)$, the restriction of the map $\ToForest$
    on the domain $\Poset(\TreeT)$ is injective.
\end{Lemma}
\begin{proof}
    This follows by structural induction on $\TreeT$ by establishing the fact that for any
    $\TreeS, \TreeS' \in \Poset(\TreeT)$, $\ToForest(\TreeS) = \ToForest\Par{\TreeS'}$
    implies $\TreeS = \TreeS'$. This uses in particular Lemma~\ref{lem:comparison_terms} in
    order to describe $\TreeS$ and $\TreeS'$ from~$\TreeT$ knowing that $\TreeT \Leq \TreeS$
    and $\TreeT \Leq \TreeS'$.
\end{proof}

\begin{Lemma} \label{lem:rewrite_relation_and_duplicative_trees}
    Let $\TreeT \in \SetTerms(\GeneratingSet)$ and $\TreeS, \TreeS' \in \SetTerms(\TreeT)$
    such that $\TreeS \ne \TreeS'$. We have $\TreeS \RewContext \TreeS'$ if and only if
    $\ToForest(\TreeS) \RewDuplicative \ToForest\Par{\TreeS'}$.
\end{Lemma}
\begin{proof}
    This follows by structural induction on $\TreeT$. This uses in particular the definition
    of the binary relation $\RewDuplicative$ on $\SetDuplicativeForests$,
    Lemma~\ref{lem:comparison_terms} in order to describe $\TreeS$ and $\TreeS'$
    from~$\TreeT$ knowing that $\TreeT \Leq \TreeS$ and $\TreeT \Leq \TreeS'$, and
    Lemma~\ref{lem:rewrite} in order to describe $\TreeS'$ from $\TreeS$ knowing that
    $\TreeS \ne \TreeS'$ and~$\TreeS \RewContext \TreeS'$.
\end{proof}

\begin{Proposition} \label{prop:poset_isomorphism}
    For any $\TreeT \in \SetTerms(\GeneratingSet)$, the posets $\Poset(\TreeT)$ and
    $\SetDuplicativeForests(\ToForest(\TreeT))$ are isomorphic.
\end{Proposition}
\begin{proof}
    Let $D := \Bra{\ToForest(\TreeT) : \TreeT \in \Poset(\TreeT)}$. By
    Lemma~\ref{lem:to_forest_injection}, the sets $\Poset(\TreeT)$ and $D$ are in one-to-one
    correspondence. Moreover, by Lemma~\ref{lem:rewrite_relation_and_duplicative_trees}, we
    have that $D = \SetDuplicativeForests(\ToForest(\TreeT))$ and that $\ToForest$ is an
    order isomorphism between $\Poset(\TreeT)$ and
    $\SetDuplicativeForests(\ToForest(\TreeT))$.
\end{proof}

\begin{figure}[ht]
    \centering
    \subfloat[][Hasse diagram of $\Poset(\TreeT)$.]{
    \begin{Page}{.55}
        \centering
        \scalebox{.8}{
        \begin{tikzpicture}[Centering,xscale=0.44,yscale=0.44,font=\footnotesize]
            \tikzstyle{GraphEdge2}=[GraphEdge,thick]
            \node(WWW)at(-11,0){$\M \Par{\VarX_1 \Par{\M \VarX_2}} \Par{\M \M}$};
            \node(WBW)at(-11,-8){$\M \Par{\VarX_1 \Par{\VarX_2 \VarX_2}} \Par{\M \M}$};
            \node(BWWW)at(0,-4){$\VarX_1 \Par{\M \VarX_2}
                \Par{\VarX_1 \Par{\M \VarX_2}} \Par{\M \M}$};
            \node(BBWW)at(-4,-8){$\VarX_1 \Par{\VarX_2 \VarX_2}
                \Par{\VarX_1 \Par{\M \VarX_2}} \Par{\M \M}$};
            \node(BWBW)at(4,-8){$\VarX_1 \Par{\M \VarX_2}
                \Par{\VarX_1 \Par{\VarX_2 \VarX_2}} \Par{\M \M}$};
            \node(BBBW)at(0,-12){$\VarX_1 \Par{\VarX_2 \VarX_2}
                \Par{\VarX_1 \Par{\VarX_2 \VarX_2}} \Par{\M \M}$};
            \draw[GraphEdge2](WWW)--(WBW);
            \draw[GraphEdge2](WWW)--(BWWW);
            \draw[GraphEdge2](WBW)--(BBBW);
            \draw[GraphEdge2](BWWW)--(BBWW);
            \draw[GraphEdge2](BWWW)--(BWBW);
            \draw[GraphEdge2](BBWW)--(BBBW);
            \draw[GraphEdge2](BWBW)--(BBBW);
        \end{tikzpicture}}
    \end{Page}
    \label{subfig:poset_terms}}
    %
    %
    \subfloat[][Hasse diagram of $\SetDuplicativeForests(\ToForest(\TreeT))$.]{
    \begin{Page}{.42}
        \centering
        \scalebox{.8}{
        \begin{tikzpicture}[Centering,xscale=0.31,yscale=0.38]
            \tikzstyle{GraphEdge2}=[GraphEdge,thick]
            \node(WWW)at(-12,0){
                \begin{tikzpicture}[Centering,xscale=0.18,yscale=0.24]
                    \node[WhiteNode](2)at(0.00,-3.75){};
                    \node[BlackNode](4)at(2.00,-1.25){};
                    \node[WhiteNode](0)at(0.00,-1.25){};
                    \node[BlackNode](1)at(0.00,-2.50){};
                    \node[BlackNode](3)at(1.00,0.00){};
                    \draw[Edge](0)--(3);
                    \draw[Edge](1)--(0);
                    \draw[Edge](2)--(1);
                    \draw[Edge](4)--(3);
                \end{tikzpicture}};
            \node(WBW)at(-12,-8){
                \begin{tikzpicture}[Centering,xscale=0.18,yscale=0.24]
                    \node[BlackNode](2)at(0.00,-3.75){};
                    \node[BlackNode](4)at(2.00,-1.25){};
                    \node[WhiteNode](0)at(0.00,-1.25){};
                    \node[BlackNode](1)at(0.00,-2.50){};
                    \node[BlackNode](3)at(1.00,0.00){};
                    \draw[Edge](0)--(3);
                    \draw[Edge](1)--(0);
                    \draw[Edge](2)--(1);
                    \draw[Edge](4)--(3);
                \end{tikzpicture}};
            \node(BWWW)at(0,-4){
                \begin{tikzpicture}[Centering,xscale=0.18,yscale=0.17]
                    \node[WhiteNode](1)at(0.00,-5.25){};
                    \node[WhiteNode](4)at(2.00,-5.25){};
                    \node[BlackNode](6)at(4.00,-1.75){};
                    \node[BlackNode](0)at(0.00,-3.50){};
                    \node[BlackNode](2)at(1.00,-1.75){};
                    \node[BlackNode](3)at(2.00,-3.50){};
                    \node[BlackNode](5)at(3.00,0.00){};
                    \draw[Edge](0)--(2);
                    \draw[Edge](1)--(0);
                    \draw[Edge](2)--(5);
                    \draw[Edge](3)--(2);
                    \draw[Edge](4)--(3);
                    \draw[Edge](6)--(5);
                \end{tikzpicture}};
            \node(BBWW)at(-4,-8){
                \begin{tikzpicture}[Centering,xscale=0.18,yscale=0.17]
                    \node[BlackNode](1)at(0.00,-5.25){};
                    \node[WhiteNode](4)at(2.00,-5.25){};
                    \node[BlackNode](6)at(4.00,-1.75){};
                    \node[BlackNode](0)at(0.00,-3.50){};
                    \node[BlackNode](2)at(1.00,-1.75){};
                    \node[BlackNode](3)at(2.00,-3.50){};
                    \node[BlackNode](5)at(3.00,0.00){};
                    \draw[Edge](0)--(2);
                    \draw[Edge](1)--(0);
                    \draw[Edge](2)--(5);
                    \draw[Edge](3)--(2);
                    \draw[Edge](4)--(3);
                    \draw[Edge](6)--(5);
                \end{tikzpicture}};
            \node(BWBW)at(4,-8){
                \begin{tikzpicture}[Centering,xscale=0.18,yscale=0.17]
                    \node[WhiteNode](1)at(0.00,-5.25){};
                    \node[BlackNode](4)at(2.00,-5.25){};
                    \node[BlackNode](6)at(4.00,-1.75){};
                    \node[BlackNode](0)at(0.00,-3.50){};
                    \node[BlackNode](2)at(1.00,-1.75){};
                    \node[BlackNode](3)at(2.00,-3.50){};
                    \node[BlackNode](5)at(3.00,0.00){};
                    \draw[Edge](0)--(2);
                    \draw[Edge](1)--(0);
                    \draw[Edge](2)--(5);
                    \draw[Edge](3)--(2);
                    \draw[Edge](4)--(3);
                    \draw[Edge](6)--(5);
                \end{tikzpicture}};
            \node(BBBW)at(0,-12){
                \begin{tikzpicture}[Centering,xscale=0.18,yscale=0.17]
                    \node[BlackNode](1)at(0.00,-5.25){};
                    \node[BlackNode](4)at(2.00,-5.25){};
                    \node[BlackNode](6)at(4.00,-1.75){};
                    \node[BlackNode](0)at(0.00,-3.50){};
                    \node[BlackNode](2)at(1.00,-1.75){};
                    \node[BlackNode](3)at(2.00,-3.50){};
                    \node[BlackNode](5)at(3.00,0.00){};
                    \draw[Edge](0)--(2);
                    \draw[Edge](1)--(0);
                    \draw[Edge](2)--(5);
                    \draw[Edge](3)--(2);
                    \draw[Edge](4)--(3);
                    \draw[Edge](6)--(5);
                \end{tikzpicture}};
            \draw[GraphEdge2](WWW)--(WBW);
            \draw[GraphEdge2](WWW)--(BWWW);
            \draw[GraphEdge2](WBW)--(BBBW);
            \draw[GraphEdge2](BWWW)--(BBWW);
            \draw[GraphEdge2](BWWW)--(BWBW);
            \draw[GraphEdge2](BBWW)--(BBBW);
            \draw[GraphEdge2](BWBW)--(BBBW);
        \end{tikzpicture}}
    \end{Page}
    \label{subfig:poset_duplicative_forests}}
    \caption{A maximal interval of the Mockingbird poset from the term $\TreeT := \M
    \Par{\VarX_1 \Par{\M \VarX_2}} \Par{\M \M}$ and its realization as a maximal interval in
    the duplicative forest poset.}
    \label{fig:example_poset_isomorphism}
\end{figure}
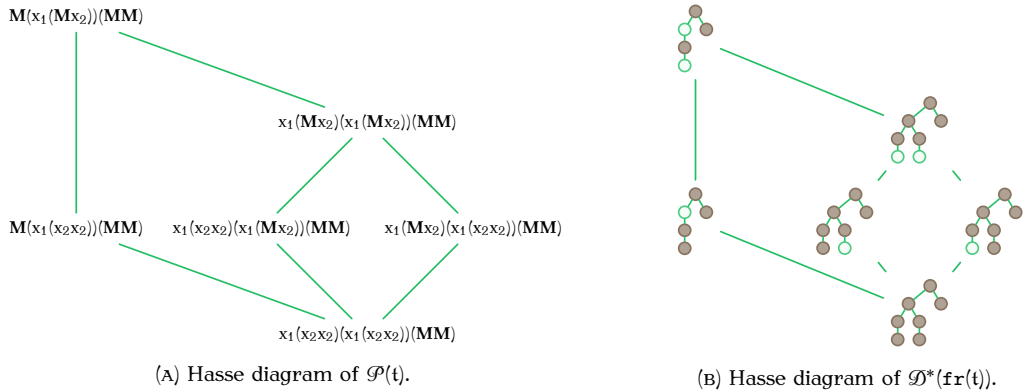

Figure~\ref{fig:example_poset_isomorphism} shows the poset $\Poset(\TreeT)$ for a
$\GeneratingSet$-term $\TreeT$ and the isomorphic poset
$\SetDuplicativeForests(\ToForest(\TreeT))$.

\begin{Theorem} \label{thm:mockingbird_lattices}
    For any $\TreeT \in \SetTerms(\GeneratingSet)$, the poset $\Poset(\TreeT)$ is a finite
    lattice.
\end{Theorem}
\begin{proof}
    By Proposition~\ref{prop:poset_isomorphism}, $\Poset(\TreeT)$ is isomorphic as a poset
    to $\SetDuplicativeForests(\ToForest(\TreeT))$. Hence, and since by
    Proposition~\ref{prop:lattice_duplicative_forests},
    $\SetDuplicativeForests(\ToForest(\TreeT))$ is a lattice, $\Poset(\TreeT)$ also is. The
    finiteness of $\Poset(\TreeT)$ is a consequence of the fact that by
    Proposition~\ref{prop:first_graph_properties}, $\CLS$ is locally finite.
\end{proof}

The \Def{Mockingbird lattice} of order $d \geq 0$ is the lattice $\MockingbirdLattice(d) :=
\Poset\Par{\RightComb_d}$ where $\RightComb_d$ is the term recursively defined by
$\RightComb_0 := \M$ and, for any $d \geq 1$, by $\RightComb_d := \M \RightComb_{d - 1}$.
Figure~\ref{fig:hasse_diagram_mockingbird_lattices} shows the Hasse diagrams of the first
Mockingbird lattices.
\begin{figure}[ht]
    \centering
    \subfloat[][$\MockingbirdLattice(0)$.]{
    \centering
    \begin{Page}{.23}
        \centering
        \begin{tikzpicture}[Centering]
            \node[GraphVertex](1)at(0,0){};
        \end{tikzpicture}
    \end{Page}
    \label{subfig:M_0}}
    \hfill
    \subfloat[][$\MockingbirdLattice(1)$.]{
    \centering
    \begin{Page}{.23}
        \centering
        \begin{tikzpicture}[Centering]
            \node[GraphVertex](1)at(0,0){};
        \end{tikzpicture}
    \end{Page}
    \label{subfig:M_1}}
    \hfill
    \subfloat[][$\MockingbirdLattice(2)$.]{
    \centering
    \begin{Page}{.23}
        \centering
        \begin{tikzpicture}[Centering,yscale=.5]
            \tikzstyle{GraphEdge2}=[GraphEdge,thin]
            \node[GraphVertex](1)at(0,0){};
            \node[GraphVertex](2)at(0,-1){};
            \draw[GraphEdge2](1)--(2);
        \end{tikzpicture}
    \end{Page}
    \label{subfig:M_2}}
    \\
    \subfloat[][$\MockingbirdLattice(3)$.]{
    \centering
    \begin{Page}{.23}
        \centering
        \begin{tikzpicture}[Centering,xscale=.3,yscale=.45]
            \tikzstyle{GraphEdge2}=[GraphEdge,thin]
            \node[GraphVertex](1)at(-1,0){};
            \node[GraphVertex](2)at(-1,-2){};
            \node[GraphVertex](3)at(2,-1){};
            \node[GraphVertex](4)at(1,-2){};
            \node[GraphVertex](5)at(3,-2){};
            \node[GraphVertex](6)at(2,-3){};
            \draw[GraphEdge2](1)--(2);
            \draw[GraphEdge2](1)--(3);
            \draw[GraphEdge2](3)--(4);
            \draw[GraphEdge2](3)--(5);
            \draw[GraphEdge2](4)--(6);
            \draw[GraphEdge2](5)--(6);
            \draw[GraphEdge2](2)--(6);
        \end{tikzpicture}
    \end{Page}
    \label{subfig:M_3}}
    \hfill
    \subfloat[][$\MockingbirdLattice(4)$.]{
    \centering
    \begin{Page}{.73}
        \centering
        \begin{tikzpicture}[Centering,xscale=.075,yscale=.075]
            \tikzstyle{GraphVertex2}=[GraphVertex,minimum size=1mm]
            \tikzstyle{GraphEdge2}=[GraphEdge,thin]
            \node[GraphVertex2](0)at(-28, -26){};
            \node[GraphVertex2](1)at(-22, -28){};
            \node[GraphVertex2](2)at(-24, -30){};
            \node[GraphVertex2](3)at(-20, -30){};
            \node[GraphVertex2](4)at(-22, -32){};
            \node[GraphVertex2](5)at(-28, -30){};
            \node[GraphVertex2](6)at(28, 0){};
            \node[GraphVertex2](7)at(34, -2){};
            \node[GraphVertex2](8)at(32, -4){};
            \node[GraphVertex2](9)at(36, -4){};
            \node[GraphVertex2](10)at(34, -6){};
            \node[GraphVertex2](11)at(28, -4){};
            \node[GraphVertex2](12)at(34, -20){};
            \node[GraphVertex2](13)at(40, -22){};
            \node[GraphVertex2](14)at(38, -24){};
            \node[GraphVertex2](15)at(42, -24){};
            \node[GraphVertex2](16)at(40, -26){};
            \node[GraphVertex2](17)at(34, -24){};
            \node[GraphVertex2](18)at(1, -14){};
            \node[GraphVertex2](19)at(7, -16){};
            \node[GraphVertex2](20)at(5, -18){};
            \node[GraphVertex2](21)at(9, -18){};
            \node[GraphVertex2](22)at(7, -20){};
            \node[GraphVertex2](23)at(1, -18){};
            \node[GraphVertex2](24)at(7, -34){};
            \node[GraphVertex2](25)at(13, -36){};
            \node[GraphVertex2](26)at(11, -38){};
            \node[GraphVertex2](27)at(15, -38){};
            \node[GraphVertex2](28)at(13, -40){};
            \node[GraphVertex2](29)at(7, -38){};
            \node[GraphVertex2](30)at(-7, 8){};
            \node[GraphVertex2](31)at(-1, 6){};
            \node[GraphVertex2](32)at(-3, 4){};
            \node[GraphVertex2](33)at(1, 4){};
            \node[GraphVertex2](34)at(-1, 2){};
            \node[GraphVertex2](35)at(-7, 4){};
            \node[GraphVertex2](36)at(-3, 18){};
            \node[GraphVertex2](37)at(38, 12){};
            \node[GraphVertex2](38)at(9, -4){};
            \node[GraphVertex2](39)at(46, -10){};
            \node[GraphVertex2](40)at(17, -26){};
            \node[GraphVertex2](41)at(-24, -16){};
            \draw[GraphEdge2](0)--(1);
            \draw[GraphEdge2](1)--(2);
            \draw[GraphEdge2](1)--(3);
            \draw[GraphEdge2](3)--(4);
            \draw[GraphEdge2](2)--(4);
            \draw[GraphEdge2](0)--(5);
            \draw[GraphEdge2](5)--(4);
            \draw[GraphEdge2](6)--(7);
            \draw[GraphEdge2](7)--(8);
            \draw[GraphEdge2](7)--(9);
            \draw[GraphEdge2](9)--(10);
            \draw[GraphEdge2](8)--(10);
            \draw[GraphEdge2](6)--(11);
            \draw[GraphEdge2](11)--(10);
            \draw[GraphEdge2](12)--(13);
            \draw[GraphEdge2](13)--(14);
            \draw[GraphEdge2](13)--(15);
            \draw[GraphEdge2](15)--(16);
            \draw[GraphEdge2](14)--(16);
            \draw[GraphEdge2](12)--(17);
            \draw[GraphEdge2](17)--(16);
            \draw[GraphEdge2](18)--(19);
            \draw[GraphEdge2](19)--(20);
            \draw[GraphEdge2](19)--(21);
            \draw[GraphEdge2](21)--(22);
            \draw[GraphEdge2](20)--(22);
            \draw[GraphEdge2](18)--(23);
            \draw[GraphEdge2](23)--(22);
            \draw[GraphEdge2](24)--(25);
            \draw[GraphEdge2](25)--(26);
            \draw[GraphEdge2](25)--(27);
            \draw[GraphEdge2](27)--(28);
            \draw[GraphEdge2](26)--(28);
            \draw[GraphEdge2](24)--(29);
            \draw[GraphEdge2](29)--(28);
            \draw[GraphEdge2](30)--(31);
            \draw[GraphEdge2](31)--(32);
            \draw[GraphEdge2](31)--(33);
            \draw[GraphEdge2](33)--(34);
            \draw[GraphEdge2](32)--(34);
            \draw[GraphEdge2](30)--(35);
            \draw[GraphEdge2](35)--(34);
            \draw[GraphEdge2](36)--(37);
            \draw[GraphEdge2](37)--(38);
            \draw[GraphEdge2](37)--(39);
            \draw[GraphEdge2](39)--(40);
            \draw[GraphEdge2](38)--(40);
            \draw[GraphEdge2](36)--(41);
            \draw[GraphEdge2](41)--(40);
            \draw[GraphEdge2](30)--(0);
            \draw[GraphEdge2](35)--(5);
            \draw[GraphEdge2](31)--(1);
            \draw[GraphEdge2](32)--(2);
            \draw[GraphEdge2](33)--(3);
            \draw[GraphEdge2](34)--(4);
            \draw[GraphEdge2](30)--(6);
            \draw[GraphEdge2](35)--(11);
            \draw[GraphEdge2](32)--(8);
            \draw[GraphEdge2](34)--(10);
            \draw[GraphEdge2](6)--(18);
            \draw[GraphEdge2](11)--(23);
            \draw[GraphEdge2](8)--(20);
            \draw[GraphEdge2](9)--(21);
            \draw[GraphEdge2](7)--(19);
            \draw[GraphEdge2](10)--(22);
            \draw[GraphEdge2](6)--(12);
            \draw[GraphEdge2](7)--(13);
            \draw[GraphEdge2](11)--(17);
            \draw[GraphEdge2](8)--(14);
            \draw[GraphEdge2](9)--(15);
            \draw[GraphEdge2](10)--(16);
            \draw[GraphEdge2](0)--(24);
            \draw[GraphEdge2](18)--(24);
            \draw[GraphEdge2](12)--(24);
            \draw[GraphEdge2](5)--(29);
            \draw[GraphEdge2](23)--(29);
            \draw[GraphEdge2](17)--(29);
            \draw[GraphEdge2](2)--(26);
            \draw[GraphEdge2](20)--(26);
            \draw[GraphEdge2](14)--(26);
            \draw[GraphEdge2](3)--(27);
            \draw[GraphEdge2](21)--(27);
            \draw[GraphEdge2](15)--(27);
            \draw[GraphEdge2](1)--(25);
            \draw[GraphEdge2](19)--(25);
            \draw[GraphEdge2](13)--(25);
            \draw[GraphEdge2](4)--(28);
            \draw[GraphEdge2](22)--(28);
            \draw[GraphEdge2](16)--(28);
            \draw[GraphEdge2](36)--(30);
            \draw[GraphEdge2](41)--(5);
            \draw[GraphEdge2](37)--(7);
            \draw[GraphEdge2](38)--(20);
            \draw[GraphEdge2](39)--(15);
            \draw[GraphEdge2](40)--(28);
            \draw[GraphEdge2](31)--(7);
            \draw[GraphEdge2](33)--(9);
        \end{tikzpicture}
    \end{Page}
    \label{subfig:M_4}}
    \caption{The Hasse diagrams of the Mockingbird lattices $\MockingbirdLattice(d)$ for
    $d \in \HanL{4}$.}
    \label{fig:hasse_diagram_mockingbird_lattices}
\end{figure}

\begin{Theorem} \label{thm:universality_mockingbird_lattices}
    For any $\ForestF \in \SetDuplicativeForests$, the poset
    $\SetDuplicativeForests(\ForestF)$ is isomorphic to a maximal interval of a Mockingbird
    lattice.
\end{Theorem}
\begin{proof}
    By Lemmas~\ref{lem:isomorphic_pruned_forests} and~\ref{lem:smaller_ladder_pruned}, the
    poset $\SetDuplicativeForests(\ForestF)$ is isomorphic to a maximal interval of
    $\SetDuplicativeForests \Par{\Ladder_d}$ for a $d \geq 0$. Since $\Pruned\Par{\ToForest
    \Par{\RightComb_{d + 1}}} = \Ladder_d$, and since by
    Proposition~\ref{prop:poset_isomorphism}, the poset $\Poset \Par{\RightComb_{d + 1}} =
    \MockingbirdLattice(d + 1)$ is isomorphic to $\SetDuplicativeForests
    \Par{\Ladder_d}$, the poset $\SetDuplicativeForests(\ForestF)$ is isomorphic to a
    maximal interval of~$\MockingbirdLattice(d + 1)$.
\end{proof}

Theorem~\ref{thm:universality_mockingbird_lattices} justifies the fact that the study of the
Mockingbird lattices is universal enough because these lattices contain as maximal interval
all duplicative forest lattices.

The Mockingbird lattices are not graded and not self-dual. Since they are not graded, they
are not distributive neither. Moreover, they are not semidistributive. For instance, in
$\MockingbirdLattice(3)$, by setting $\TreeT_1 := \M(\M\M(\M\M))$, $\TreeT_2 :=
\M\M(\M\M)(\M(\M\M))$, and $\TreeT_3 := \M(\M\M)(\M\M(\M\M))$, we have $\TreeT_1 \Meet
\TreeT_2 = \TreeT_1 \Meet \TreeT_3$ but $\TreeT_1 \Meet \Par{\TreeT_2 \JJoin \TreeT_3} \ne
\TreeT_1 \Meet \TreeT_2$, contradicting one of the required relations to have this property.

\section{Enumerative properties} \label{sec:enumerative_properties}
In this last central part of this work, we present some enumerative results concerning the
Mockingbird poset and lattices $\MockingbirdLattice(d)$, $d \geq 0$. We enumerate the
maximal and minimal elements of the Mockingbird poset by degree and by height, the length of
the shortest and longest saturated chains in $\MockingbirdLattice(d)$, the number of
elements in $\MockingbirdLattice(d)$, the number of edges in the Hasse diagram of
$\MockingbirdLattice(d)$, and the number of intervals of $\MockingbirdLattice(d)$. All this
use the lattice isomorphism between the Mockingbird lattices and the duplicative forest
lattices introduced in the previous sections and formal series on terms and on duplicative
forests.

\subsection{Formal power series}
We set here some notions about formal power series and generating series. For the rest of
the text, $\K$ is any field of characteristic zero (as $\Q$ for instance).

\subsubsection{Formal power series over sets}
For any set $X$, let $\K \Angle{X}$ be the linear span of $X$. The dual space of $\K
\Angle{X}$ is denoted by $\K \AAngle{X}$ and is by definition the space of the maps
$\SeriesF : X \to \K$, called \Def{$X$-series}. The coefficient $\SeriesF(x)$ of any $x \in
X$ is denoted by $\Angle{x, \SeriesF}$. The \Def{support} of $\SeriesF$ is the set
$\Support(\SeriesF) := \Bra{x \in X : \Angle{x, \SeriesF} \ne 0}$. The \Def{characteristic
series} of any subset $X'$ of $X$ is the series $\CharacteristicSeries\Par{X'}$ having $X'$
as support and such that the coefficient of each $x \in X'$ is $1$. For any $k \geq 0$,
$\T{k} \K \AAngle{X}$ is the $k$-th tensor power of $\K \AAngle{X}$. Elements of this space
are possibly infinite linear combinations of tensors $x_1 \otimes \dots \otimes x_k$, where
for any $i \in [k]$, $x_i \in X$. The \Def{tensor algebra} of $\K \AAngle{X}$ is the space
\begin{equation}
    \T{*} \K \AAngle{X} := \bigoplus_{k \geq 0} \T{k} \K \AAngle{X}.
\end{equation}
This space is endowed with the tensor product $\otimes$ so that $\Par{\T{*} \K \AAngle{X},
\otimes}$ is a unital associative algebra, admitting $1 \in \K$ as unit.

A linear map $\phi : \T{k_1} \K \AAngle{X} \to \T{k_2} \K \AAngle{X}$, $k_1, k_2 \geq 0$, is
a \Def{$\Par{k_1, k_2}$-operation} on $\K \AAngle{X}$. To lighten the notation, when $\phi$
is a $(2, 1)$-operation on $\K \AAngle{X}$, we shall sometimes use $\phi$ as an infix
operation by writing $\SeriesF_1 \, \phi \, \SeriesF_2$ instead of $\phi\Par{\SeriesF_1
\otimes \SeriesF_2}$ for any $\SeriesF_1, \SeriesF_2 \in \K \AAngle{X}$. The \Def{diagonal
coproduct} is the $(1, 2)$-operation $\Delta$ on $\K \AAngle{X}$ satisfying $\Delta(x) = x
\otimes x$ for any $x \in X$. When $X$ is endowed with an $n$-ary operation $\Product : X^n
\to X$, $n \geq 0$, for any $k \geq 0$, the \Def{$k$-linearization} of $\Product$ is the $(n
k, k)$-operation $\TT{\Product}{k}$ on $\K \AAngle{X}$ satisfying
\begin{equation}
    \TT{\Product}{k} \Par{
        x_{1, 1} \otimes \dots \otimes x_{1, k}
        \otimes \dots \otimes
        x_{n, 1} \otimes \dots \otimes x_{n, k}}
    =
    \Product\Par{x_{1, 1}, \dots, x_{n, 1}}
    \otimes \dots \otimes
    \Product\Par{x_{k, 1}, \dots, x_{n, k}}
\end{equation}
for any $x_{i, j} \in X$, $i \in [n]$, and $j \in [k]$. To lighten the notation, we shall
write $\TT{\Product}{}$ for~$\TT{\Product}{1}$.

\subsubsection{Generating series}
The space of the usual power series on the formal parameter $\VarZ$ is denoted by $\K
\AAngle{\VarZ}$. For any $F, F' \in \K \AAngle{\VarZ}$, $F\Han{\VarZ := F'}$ is the series
of $\K \AAngle{\VarZ}$ obtained by substituting $F'$ for $\VarZ$ in $F$. The \Def{Hadamard
product} is the binary operation $\HadamardProduct$ on $\K \AAngle{\VarZ}$ defined linearly
for any $n_1, n_2 \geq 0$ by
\begin{equation}
    \VarZ^{n_1} \HadamardProduct \VarZ^{n_2} := \Iverson{n_1 = n_2} \, z^{n_1}.
\end{equation}
The \Def{max product} is the binary operation $\MaxProduct$ on $\K \AAngle{\VarZ}$ defined
linearly for any $n_1, n_2 \geq 0$ by $\VarZ^{n_1} \MaxProduct \VarZ^{n_2} := \VarZ^{\max
\Bra{n_1, n_2}}$. Observe that for any $F \in \K \AAngle{\VarZ}$ and $n \geq 0$,
\begin{equation} \label{equ:coefficient_max_product}
    \Angle{\VarZ^n, F \MaxProduct F}
    = \Angle{\VarZ^n, F}^2
    + 2 \Angle{\VarZ^n, F} \sum_{i \in [n]} \Angle{\VarZ^{i - 1}, F}.
\end{equation}

\subsubsection{Formal power series and enumeration}
If $X$ is endowed with a map $\omega : X \to \N$, the \Def{$\omega$-enumeration map} is the
partial map $\EnumerationMap_\omega : \T{*} \K \AAngle{X} \to \K \AAngle{\VarZ}$ defined
linearly for any $k \geq 1$ and $x_1, \dots, x_k \in X$ by
\begin{equation}
    \EnumerationMap_\omega\Par{x_1 \otimes \dots \otimes x_k}
    := \VarZ^{\omega\Par{x_1}} \MaxProduct \dots \MaxProduct \VarZ^{\omega\Par{x_k}}.
\end{equation}
For any $\SeriesF \in \T{*} \K \AAngle{X}$, the generating series
$\EnumerationMap_\omega(\SeriesF)$ is the \Def{$\omega$-enumeration} of $\SeriesF$.

In the sequel, we shall use the following strategy to enumerate a set $X$ w.r.t.\ such a map
$\omega$: we shall provide a description of $\CharacteristicSeries(X)$, then deduce a
description of $\EnumerationMap_\omega(\CharacteristicSeries(X))$, and finally deduce from
this a formula to compute the coefficients $\Angle{\VarZ^n,
\EnumerationMap_\omega(\CharacteristicSeries(X))}$, $n \geq 0$.

\subsection{Maximal, minimal, and isolated terms}
We use here series on terms and some operations on these in order to enumerate the closed
maximal, closed minimal, and closed isolated elements of the poset $\Poset$ w.r.t.\ their
degrees or their heights. These enumerative results can be used to establish the
probabilities for a term $\TreeT$ of a given degree or of a given height generated uniformly
at random to have one of these properties. In particular, the number of maximal (resp.\
minimal, isolated) closed terms gives information to compute the probability for the fact
that $\TreeT$ cannot be rewritten into a different term (resp.\ no term different from
$\TreeT$ rewrites into $\TreeT$, any sequence of rewriting steps involving $\TreeT$ does not
involve different other terms).

\subsubsection{Series of terms and operations}
Recall that $\App$ is the binary operation on $\SetTerms(\GeneratingSet)$ such that for
any $\TreeT_1, \TreeT_2 \in \SetTerms(\GeneratingSet)$, $\TreeT_1 \App \TreeT_2$ is the
$\GeneratingSet$-term $\TreeT_1 \TreeT_2$ having $\TreeT_1$ as left subtree and $\TreeT_2$
as right subtree. Observe that for any $\SeriesF_1, \SeriesF_2 \in \K
\AAngle{\SetTerms(\GeneratingSet)}$,
\begin{equation} \label{equ:deg_enumeration_application}
    \EnumerationMap_\Deg \Par{\SeriesF_1 \TT{\App}{} \SeriesF_2}
    = \VarZ \, \EnumerationMap_\Deg\Par{\SeriesF_1} \, \EnumerationMap_\Deg\Par{\SeriesF_2}
\end{equation}
and
\begin{equation} \label{equ:ht_enumeration_application}
    \EnumerationMap_\Height \Par{\SeriesF_1 \TT{\App}{} \SeriesF_2}
    = \VarZ \, \Par{
        \EnumerationMap_\Height \Par{\SeriesF_1}
        \MaxProduct \EnumerationMap_\Height \Par{\SeriesF_2}}.
\end{equation}
Observe also that for any $\SeriesF \in \K \AAngle{\SetTerms(\GeneratingSet)}$,
\begin{equation} \label{equ:deg_enumeration_delta_application}
    \EnumerationMap_\Deg \Par{\TT{\App}{}(\Delta(\SeriesF))}
    = \VarZ \, \Par{\EnumerationMap_\Deg(\SeriesF)\Han{\VarZ := \VarZ^2}}
\end{equation}
and
\begin{equation} \label{equ:ht_enumeration_delta_application}
    \EnumerationMap_{\Height} \Par{\TT{\App}{}(\Delta(\SeriesF))}
    = \VarZ \, \EnumerationMap_{\Height}(\SeriesF).
\end{equation}

\subsubsection{Three series of terms}
In order to achieve the objectives described above, we begin by providing equations
satisfied by the characteristic series $\SeriesMax$ of the closed maximal terms of $\Poset$,
the characteristic  series $\SeriesMin$ of the closed minimal terms of $\Poset$, and the
characteristic series $\SeriesIsolated$ of the closed isolated terms of~$\Poset$.

\begin{Proposition} \label{prop:maximal_elements_enumeration}
    The characteristic series $\SeriesMax$ satisfies
    \begin{equation} \label{equ:maximal_elements_enumeration}
        \SeriesMax
        = \M + \M \M + \SeriesMax \TT{\App}{} \SeriesMax
        - \M \TT{\App}{} \SeriesMax.
    \end{equation}
\end{Proposition}
\begin{proof}
    By using the description of Proposition~\ref{prop:minimal_maximal} for the set $S$ of
    closed maximal elements of $\Poset$, we have
    \begin{equation}
        \SeriesMax
        =
        \M + \sum_{\TreeS_1, \TreeS_2 \in S} \TreeS_1 \TreeS_2
        - \sum_{\TreeS_3 \in S \setminus \{\M\}} \M \TreeS_3
        = \M + \SeriesMax \TT{\App}{} \SeriesMax
        - \M \TT{\App}{} \Par{\SeriesMax - \M}.
    \end{equation}
    This leads to~\eqref{equ:maximal_elements_enumeration}.
\end{proof}

\begin{Proposition} \label{prop:minimal_elements_enumeration}
    The characteristic series $\SeriesMin$ satisfies
    \begin{equation} \label{equ:minimal_elements_enumeration}
        \SeriesMin
        = \M + \M \M + \SeriesMin \TT{\App}{} \SeriesMin
        - \TT{\App}{}\Par{\Delta(\SeriesMin)}.
    \end{equation}
\end{Proposition}
\begin{proof}
    By using the description of Proposition~\ref{prop:minimal_maximal} for the set $S$ of
    closed minimal elements of $\Poset$, we have
    \begin{equation}
        \SeriesMin
        =
        \M + \sum_{\TreeS_1, \TreeS_2 \in S} \TreeS_1 \TreeS_2
        - \sum_{\TreeS_3 \in S \setminus \{\M\}} \TreeS_3 \TreeS_3
        = \M + \SeriesMin \TT{\App}{} \SeriesMin
        - \TT{\App}{}\Par{\Delta\Par{\SeriesMin - \M}}.
    \end{equation}
    This leads to~\eqref{equ:minimal_elements_enumeration}.
\end{proof}

\begin{Proposition} \label{prop:isolated_elements_enumeration}
    The characteristic series $\SeriesIsolated$ satisfies
    \begin{equation} \label{equ:isolated_elements_enumeration}
        \SeriesIsolated
        = \M + 2 \M \M + \SeriesIsolated \TT{\App}{} \SeriesIsolated
        - \M \TT{\App}{} \SeriesIsolated
        - \TT{\App}{}\Par{\Delta\Par{\SeriesIsolated}}.
    \end{equation}
\end{Proposition}
\begin{proof}
    By using the description of Proposition~\ref{prop:minimal_maximal} for the set $S_1$
    (resp.\ $S_2$) of the closed maximal (resp.\ minimal) elements of $\Poset$, the set of
    the closed isolated elements of $\Poset$ is the set $S = S_1 \cap S_2$ and we have
    \begin{equation} \begin{split}
        \SeriesIsolated
        & =
        \M + \sum_{\TreeS_1, \TreeS_2 \in S} \TreeS_1 \TreeS_2
        - \sum_{\TreeS_3 \in S \setminus \{\M\}} \M \TreeS_3
        - \sum_{\TreeS_4 \in S \setminus \{\M\}} \TreeS_4 \TreeS_4 \\
        & = \M + \SeriesIsolated \TT{\App}{} \SeriesIsolated
        - \M \TT{\App}{} \Par{\SeriesIsolated - \M}
        - \TT{\App}{}\Par{\Delta\Par{\SeriesIsolated - \M}}.
    \end{split} \end{equation}
    This leads to~\eqref{equ:isolated_elements_enumeration}.
\end{proof}

\subsubsection{$\Deg$-enumerations}
A consequence of Proposition~\ref{prop:maximal_elements_enumeration} and
of~\eqref{equ:deg_enumeration_application} is that the $\Deg$-enumeration
$\GeneratingSeriesMaxDeg$ of $\SeriesMax$, enumerating the closed maximal elements of
$\Poset$ w.r.t.\ their degrees, satisfies
\begin{equation}
    \GeneratingSeriesMaxDeg
    = 1 + \VarZ + \VarZ \GeneratingSeriesMaxDeg^2 - \VarZ \GeneratingSeriesMaxDeg.
\end{equation}
We deduce from this that the number of these terms of degree $d \geq 0$ is
$\SequenceMaxDeg(d)$ where $\SequenceMaxDeg$ is the integer sequence satisfying
$\SequenceMaxDeg(0) = \SequenceMaxDeg(1) = 1$ and, for any $d \geq 2$,
\begin{equation}
    \SequenceMaxDeg(d)
    = \sum_{i \in \HanL{d - 2}} \SequenceMaxDeg(i) \ \SequenceMaxDeg(d - 1 - i).
\end{equation}
The first numbers are
\begin{equation}
    1, 1, 1, 2, 4, 9, 21, 51,
\end{equation}
form Sequence~\OEIS{A001006} of~\cite{Slo}, and are Motzkin numbers. A \Def{Motzkin tree} is
a planar rooted tree made of nodes of arity $0$, $1$, or $2$. It is well-known that the sets
of these trees having exactly $d \geq 1$ nodes are enumerated by Motzkin numbers
(see~\cite{DP02} for instance). Given a closed maximal term $\TreeT$, it is easy to see that
$\ToForest(\TreeT)$ contains only black nodes and that it is a Motzkin tree. It is also easy
to check that $\ToForest$ is a bijection between the set of the closed maximal terms of
degree $d \geq 1$ and the set of the Motzkin trees with $d$ nodes.

Besides, a consequence of Proposition~\ref{prop:minimal_elements_enumeration},
\eqref{equ:deg_enumeration_application}, and~\eqref{equ:deg_enumeration_delta_application}
is that the $\Deg$-enumeration $\GeneratingSeriesMinDeg$ of $\SeriesMin$, enumerating the
closed minimal elements of $\Poset$ w.r.t.\ their degrees, satisfies
\begin{equation}
    \GeneratingSeriesMinDeg
    = 1 + \VarZ + \VarZ \GeneratingSeriesMinDeg^2
    - \VarZ \Par{\GeneratingSeriesMinDeg\Han{\VarZ := \VarZ^2}}.
\end{equation}
We deduce from this that the number of these terms of degree $d \geq 0$ is
$\SequenceMinDeg(d)$ where $\SequenceMinDeg$ is the integer sequence satisfying
$\SequenceMinDeg(0) = \SequenceMinDeg(1) = 1$ and, for any $d \geq 2$,
\begin{equation}
    \SequenceMinDeg(d) =
    \sum_{i \in \HanL{d - 1}} \SequenceMinDeg(i) \ \SequenceMinDeg(d - 1 - i)
    - \Iverson{d \mbox{ is odd}} \, \SequenceMinDeg\Par{\frac{d - 1}{2}}.
\end{equation}
The first numbers are
\begin{equation}
    1, 1, 2, 4, 12, 34, 108, 344
\end{equation}
and form Sequence~\OEIS{A343663} of~\cite{Slo}. A \Def{semi-identity tree} is a binary tree
having no subtrees that are brother and equal, unless they are leaves. It is known that the
sets of these trees having exactly $d \geq 0$ internal nodes are enumerated by the
previous sequence. It is immediate to see that the identity map is a bijection between the
set of the closed minimal terms of degree $d \geq 0$ and the set of the semi-identity trees
with exactly $d$ internal nodes.

Besides, a consequence of Proposition~\ref{prop:isolated_elements_enumeration},
\eqref{equ:deg_enumeration_application}, and~\eqref{equ:deg_enumeration_delta_application}
is that the $\Deg$-enumeration $\GeneratingSeriesIsolatedDeg$ of $\SeriesIsolated$,
enumerating the closed isolated elements of $\Poset$ w.r.t.\ their degrees, satisfies
\begin{equation}
    \GeneratingSeriesIsolatedDeg
    = 1 + 2\VarZ + \VarZ \GeneratingSeriesIsolatedDeg^2 - \VarZ \GeneratingSeriesIsolatedDeg
    - \VarZ \Par{\GeneratingSeriesIsolatedDeg\Han{\VarZ := \VarZ^2}}.
\end{equation}
We deduce from this that the number of these terms of degree $d \geq 0$ is
$\SequenceIsolatedDeg(d)$ where $\SequenceIsolatedDeg$ is the integer sequence satisfying
$\SequenceIsolatedDeg(0) = \SequenceIsolatedDeg(1) = 1$ and, for any $d \geq 2$,
\begin{equation}
    \SequenceIsolatedDeg(d) =
    \sum_{i \in \HanL{d - 2}} \SequenceIsolatedDeg(i) \ \SequenceIsolatedDeg(d - 1 - i)
    - \Iverson{d \mbox{ is odd}} \, \SequenceIsolatedDeg\Par{\frac{d - 1}{2}}.
\end{equation}
The first numbers are
\begin{equation}
    1, 1, 1, 1, 3, 5, 13, 29, 71, 171, 427, 1067, 2709.
\end{equation}
This sequence does not appear in~\cite{Slo} for the time being.

\subsubsection{$\Height$-enumerations}
Due to the fact that $\CLS$ is hierarchical and, by
Proposition~\ref{prop:first_graph_properties}, $\Poset$ is rooted,
Proposition~\ref{prop:finite_equivalence_classes} implies that $\Poset$ satisfies the
properties exposed at the very end of Section~\ref{subsubsec:properties_CLS}. In particular,
in each $\Equiv$-equivalence class of closed terms, there is exactly one closed maximal
term, exactly one closed minimal term, and these terms have the same height. For this
reason, the $\Height$-enumerations of $\SeriesMax$ and $\SeriesMin$ are equal and enumerate
also the $\Equiv$-equivalence classes of closed terms w.r.t.\ the heights of their elements.
By denoting by $\GeneratingSeriesMinMaxHt$ this generating series, by
Propositions~\ref{prop:maximal_elements_enumeration}
and~\ref{prop:minimal_elements_enumeration}, $\GeneratingSeriesMinMaxHt$ satisfies
\begin{equation} \label{equ:number_equivalence_classes_height}
    \GeneratingSeriesMinMaxHt
    = 1 + \VarZ
    + \VarZ \Par{\GeneratingSeriesMinMaxHt \MaxProduct \GeneratingSeriesMinMaxHt}
    - \VarZ \GeneratingSeriesMinMaxHt.
\end{equation}
We deduce from this and~\eqref{equ:coefficient_max_product} that the number of these
$\Equiv$-equivalence classes of terms of height $h \geq 0$ is $\SequenceMinMaxHt(h)$ where
$\SequenceMinMaxHt$ is the integer sequence satisfying $\SequenceMinMaxHt(0) =
\SequenceMinMaxHt(1) = 1$ and, for any $h \geq 2$,
\begin{equation}
    \SequenceMinMaxHt(h) =
    \SequenceMinMaxHt(h - 1)^2 - \SequenceMinMaxHt(h - 1)
    + 2 \SequenceMinMaxHt(h - 1) \sum_{i \in [h - 1]} \SequenceMinMaxHt(i - 1).
\end{equation}
The first numbers are
\begin{equation}
    1, 1, 2, 10, 170, 33490, 1133870930, 1285739648704587610
\end{equation}
and form Sequence~\OEIS{A063573} of~\cite{Slo}.

Besides, a consequence of Proposition~\ref{prop:isolated_elements_enumeration},
\eqref{equ:ht_enumeration_application}, and~\eqref{equ:ht_enumeration_delta_application} is
that the $\Height$-enumeration $\GeneratingSeriesIsolatedHt$ of $\SeriesIsolated$,
enumerating the closed isolated elements of $\Poset$ w.r.t.\ their heights, satisfies
\begin{equation}
    \GeneratingSeriesIsolatedHt =
    1 + 2 z + z \Par{\GeneratingSeriesIsolatedHt \MaxProduct \GeneratingSeriesIsolatedHt}
    - 2z \GeneratingSeriesIsolatedHt.
\end{equation}
We deduce from this and~\eqref{equ:coefficient_max_product} that the number of these terms
of height $h \geq 0$ is $\SequenceIsolatedHt(h)$ where $\SequenceIsolatedHt$ is the integer
sequence satisfying $\SequenceIsolatedHt(0) = \SequenceIsolatedHt(1) = 1$ and, for any $h
\geq 2$,
\begin{equation}
    \SequenceIsolatedHt(h) =
    \SequenceIsolatedHt(h - 1)^2 - 2 \SequenceIsolatedHt(h - 1)
    + 2 \SequenceIsolatedHt(h - 1) \sum_{i \in [h - 1]} \SequenceIsolatedHt(i - 1).
\end{equation}
The first numbers are
\begin{equation}
    1, 1, 1, 3, 21, 651, 457653, 210065930571
\end{equation}
and form Sequence~\OEIS{A001699} of~\cite{Slo}, which enumerates binary trees w.r.t.\ the
height (but with a shift). Therefore, there is a one-to-one correspondence between the set
of the closed isolated elements of $\Poset$ of height $h \geq 1$ and the set
of the binary trees of height~$h - 1$.

\subsection{Shortest and longest saturated chains}
Let $\MaxLength : \SetDuplicativeForests \to \N$ be the statistics defined for any $\ell
\geq 0$, $\ForestF_1, \dots, \ForestF_\ell \in \SetDuplicativeForests$, and $\ForestF \in
\SetDuplicativeForests$ by
\begin{subequations}
\begin{equation} \label{equ:max_length_1}
    \MaxLength\Par{\ForestF_1 \ConcatenateForests \dots \ConcatenateForests \ForestF_\ell}
    := 1 - \ell + \sum_{i \in [\ell]} \MaxLength \Par{\ForestF_i},
\end{equation}
\begin{equation} \label{equ:max_length_2}
    \MaxLength\Par{\BlackNode\Par{\ForestF}} := \MaxLength(\ForestF),
\end{equation}
\begin{equation} \label{equ:max_length_3}
    \MaxLength\Par{\WhiteNode\Par{\ForestF}} := 2 \MaxLength(\ForestF).
\end{equation}
\end{subequations}
For instance, we have
\begin{equation}
    \begin{tikzpicture}[Centering,xscale=0.29,yscale=0.19]
        \node[WhiteNode](0)at(-1.00,-4.40){};
        \node[WhiteNode](10)at(8.00,-4.40){};
        \node[WhiteNode](2)at(1.00,-4.40){};
        \node[BlackNode](4)at(2.00,-8.80){};
        \node[WhiteNode](6)at(4.00,-6.60){};
        \node[WhiteNode](8)at(6.00,-2.20){};
        \node[BlackNode](1)at(1.00,-2.20){};
        \node[WhiteNode](3)at(2.00,-6.60){};
        \node[WhiteNode](5)at(3.00,-4.40){};
        \node[BlackNode](9)at(8.00,-2.20){};
        \draw[Edge](0)--(1);
        \draw[Edge](10)--(9);
        \draw[Edge](2)--(1);
        \draw[Edge](3)--(5);
        \draw[Edge](4)--(3);
        \draw[Edge](5)--(1);
        \draw[Edge](6)--(5);
        \node[left of=0,node distance=3mm,font=\footnotesize]{$2$};
        \node[left of=1,node distance=3mm,font=\footnotesize]{$8$};
        \node[left of=2,node distance=3mm,font=\footnotesize]{$2$};
        \node[left of=3,node distance=3mm,font=\footnotesize]{$2$};
        \node[left of=4,node distance=3mm,font=\footnotesize]{$1$};
        \node[left of=5,node distance=3mm,font=\footnotesize]{$6$};
        \node[left of=6,node distance=3mm,font=\footnotesize]{$2$};
        \node[left of=8,node distance=3mm,font=\footnotesize]{$2$};
        \node[left of=9,node distance=3mm,font=\footnotesize]{$2$};
        \node[left of=10,node distance=3mm,font=\footnotesize]{$2$};
    \end{tikzpicture}
    \enspace \xmapsto{\MaxLength} \enspace
    10,
\end{equation}
where the integer decoration of each node refers to the image by $\MaxLength$ of the
duplicative tree rooted at this node.

\begin{Proposition} \label{prop:paths_lengths}
    For any $\TreeT \in \SetTerms_0(\GeneratingSet)$,
    \begin{enumerate}[label=(\roman*)]
        \item \label{item:paths_lengths_1}
        a shortest saturated chain from $\TreeT$ to the maximal element of $\Poset(\TreeT)$
        has as length the number of $\WhiteNode$ in $\ToForest(\TreeT)$;
        \item \label{item:paths_lengths_2}
        a longest saturated chain from $\TreeT$ to the maximal element of $\Poset(\TreeT)$
        has length $\MaxLength(\ToForest(\TreeT))$.
    \end{enumerate}
\end{Proposition}
\begin{proof}
    Let $\ForestF := \ToForest(\TreeT)$. By Proposition~\ref{prop:poset_isomorphism}, the
    posets $\Poset(\TreeT)$ and $\SetDuplicativeForests(\ForestF)$ are isomorphic.
    Therefore, to prove the statement, we consider the lengths of the shortest and longest
    saturated chains in $\SetDuplicativeForests(\ForestF)$ from $\ForestF$ to $\ForestG$
    where $\ForestG$ is the maximal element of $\SetDuplicativeForests(\ForestF)$.

    First, immediately from the definition of $\LeqDuplicative$, a shortest saturated chain
    from $\ForestF$ to $\ForestG$ consists in selecting at each step a white node that has
    no white nodes as descendants and turn it into black. This
    establishes~\ref{item:paths_lengths_1}.

    Let us prove~\ref{item:paths_lengths_2} by structural induction on $\ForestF$. We use
    Lemma~\ref{lem:comparison_duplicative_forests} in order to exhibit a longest saturated
    chain from $\ForestF$ to $\ForestG$. If $\ForestF = \ForestF_1 \ConcatenateForests \dots
    \ConcatenateForests \ForestF_\ell$, $\ell \geq 0$, a longest saturated chain from
    $\ForestF$ passes through the duplicative forest $\ForestF_1' \ConcatenateForests
    \ForestF_2 \ConcatenateForests \dots \ConcatenateForests \ForestF_\ell$ and then passes
    through $\ForestF_1' \ConcatenateForests \ForestF_2' \ConcatenateForests \dots
    \ConcatenateForests \ForestF_\ell$, \dots, and ends at $\ForestF_1' \ConcatenateForests
    \ForestF_2' \ConcatenateForests \dots \ConcatenateForests \ForestF_\ell' = \ForestG$,
    where for any $i \in [\ell]$, $\ForestF_i'$ is the greatest element of
    $\SetDuplicativeForests\Par{\ForestF_i}$. By induction hypothesis,
    $\MaxLength\Par{\ForestF_i}$ is the length of a longest saturated chain from
    $\ForestF_i$ to $\ForestF_i'$. We deduce that the length of the former chain from
    $\ForestF$ to $\ForestG$ is $1 + \sum_{i \in [\ell]} \Par{\MaxLength\Par{\ForestF_i} -
    1}$. Therefore, \eqref{equ:max_length_1} is consistent. If $\ForestF =
    \BlackNode\Par{\ForestF'}$, $\ForestF' \in \SetDuplicativeForests$, a longest saturated
    saturated chain from $\ForestF$ to $\ForestG$ has the same length as a longest saturated
    chain from $\ForestF'$ to $\ForestG'$ where $\ForestG'$ is the greatest element of
    $\SetDuplicativeForests\Par{\ForestF'}$. Hence, by induction hypothesis,
    $\MaxLength(\ForestF)$ is consistent with~\eqref{equ:max_length_2}. If $\ForestF =
    \WhiteNode\Par{\ForestF'}$, $\ForestF' \in \SetDuplicativeForests$, a longest saturated
    chain from $\ForestF$ to $\ForestG$ passes through the covering
    $\BlackNode\Par{\ForestF' \ConcatenateForests \ForestF'}$ of $\ForestF$. Therefore, by
    induction hypothesis and by~\eqref{equ:max_length_2} and~\eqref{equ:max_length_1}, a
    longest saturated chain from $\ForestF$ to $\ForestG$ has length $1 +
    \MaxLength\Par{\BlackNode\Par{\ForestF' \ConcatenateForests \ForestF'}} = 1 + 1 - 2 + 2
    \MaxLength\Par{\ForestF'} = 2 \MaxLength\Par{\ForestF'}$, which is consistent
    with~\eqref{equ:max_length_3}.
\end{proof}

By Proposition~\ref{prop:paths_lengths}, for any $d \geq 1$, in $\MockingbirdLattice(d)$,
shortest saturated chains are of length $d$ and longest saturated chains are of length
$2^{d - 1}$.

\subsection{Elements, covering pairs, and intervals}
We use here series on duplicative forests and a palette of operations on these to enumerate
the elements, covering pairs, and intervals of the Mockingbird lattices. We obtain for these
three recursive formulas leading to expressions for the generating series of these numbers.
These enumerative results can be used to establish the probability for a pair
$\Par{\TreeT_1, \TreeT_2}$ of terms both of a given height and generated uniformly at random
to have one of these properties. In particular, the number of covering pairs (resp.\ of
intervals) gives information to compute the probability for the fact that $\TreeT_1$ can be
rewritten in one step (resp.\ in any number of steps) into~$\TreeT_2$.

\subsubsection{Series of duplicative forests and operations}
We shall consider in the sequel the $k$-linearization $\TT{\ConcatenateForests}{k}$, $k \geq
0$, of the concatenation product $\ConcatenateForests$ of duplicative forests. Observe that
for any $\SeriesF_1, \SeriesF_2 \in \K
\AAngle{\SetDuplicativeForests}$, we have in particular
\begin{equation} \label{equ:ht_enumeration_concatenation}
    \EnumerationMap_{\Height}\Par{\SeriesF_1 \TT{\ConcatenateForests}{} \SeriesF_2}
    = \EnumerationMap_{\Height} \Par{\SeriesF_1}
        \MaxProduct \EnumerationMap_{\Height} \Par{\SeriesF_2}.
\end{equation}
We shall also consider the $k$-linearization $\TT{\AnyNode}{k}$, $k \geq 0$ of the grafting
product~$\AnyNode$. Observe that for any $\SeriesF \in\K \AAngle{\SetDuplicativeForests}$,
we have in particular
\begin{equation} \label{equ:ht_enumeration_graft_products}
    \EnumerationMap_{\Height}\Par{\TT{\AnyNode}{}(\SeriesF)}
    = \VarZ \, \EnumerationMap_{\Height}(\SeriesF).
\end{equation}

For any $k \geq 1$ and $u \in \Bra{\WhiteNode, \BlackNode}^k$, the \Def{merging product} is
the $\Par{k + |u|_{\BlackNode}, k}$-operation $\Merge_u$ on $\K
\AAngle{\SetDuplicativeForests}$ satisfying, for any $\ForestF_1, \dots, \ForestF_{k +
|u|_{\BlackNode}} \in \SetDuplicativeForests$,
\begin{subequations}
\begin{equation}
    \Merge_{\WhiteNode}\Par{\ForestF_1} = \WhiteNode \Par{\ForestF_1},
\end{equation}
\begin{equation}
    \Merge_{\WhiteNode u'}\Par{\ForestF_1 \otimes \dots \otimes \ForestF_k}
        = \Merge_{\WhiteNode}\Par{\ForestF_1}
            \otimes \Merge_{u'}\Par{\ForestF_2 \otimes \dots \otimes \ForestF_k},
\end{equation}
\begin{equation}
    \Merge_{\BlackNode}\Par{\ForestF_1 \otimes \ForestF_2}
    = \BlackNode \Par{\ForestF_1 \ConcatenateForests \ForestF_2},
\end{equation}
\begin{equation}
    \Merge_{\BlackNode u'}\Par{\ForestF_1 \otimes \dots \otimes \ForestF_k}
        = \Merge_{\BlackNode}\Par{\ForestF_1 \otimes \ForestF_2}
        \otimes \Merge_{u'}\Par{\ForestF_3 \otimes \dots \otimes \ForestF_k},
\end{equation}
\end{subequations}
where $u' \in \Bra{\WhiteNode, \BlackNode}^*$. For instance, with $k = 3$,
\begin{equation}
    \Merge_{\WhiteNode \, \BlackNode \, \BlackNode}
    \Par{
    \scalebox{.75}{
    \begin{tikzpicture}[Centering,xscale=0.17,yscale=0.3]
        \node[BlackNode](1)at(0.00,-1.00){};
        \node[BlackNode](0)at(0.00,0.00){};
        \draw[Edge](1)--(0);
    \end{tikzpicture}}
    \otimes
    \scalebox{.75}{
    \begin{tikzpicture}[Centering,xscale=0.17,yscale=0.2]
        \node[WhiteNode](0)at(0.00,-1.50){};
        \node[WhiteNode](2)at(2.00,-1.50){};
        \node[WhiteNode](1)at(1.00,0.00){};
        \draw[Edge](0)--(1);
        \draw[Edge](2)--(1);
    \end{tikzpicture}}
    \otimes
    \scalebox{.75}{
    \begin{tikzpicture}[Centering,xscale=0.17,yscale=0.24]
        \node[WhiteNode](0)at(0.00,-1.33){};
        \node[BlackNode](3)at(2.00,-2.67){};
        \node[WhiteNode](2)at(2.00,-1.33){};
        \draw[Edge](3)--(2);
    \end{tikzpicture}}
    \otimes
    \scalebox{.75}{
    \begin{tikzpicture}[Centering,xscale=0.14,yscale=0.2]
        \node[WhiteNode](0)at(0.00,-1.50){};
        \node[BlackNode](2)at(2.00,-1.50){};
    \end{tikzpicture}}
    \otimes
    \scalebox{.75}{
    \begin{tikzpicture}[Centering,xscale=0.14,yscale=0.2]
        \node[WhiteNode](0)at(0.00,-1.50){};
        \node[WhiteNode](2)at(2.00,-1.50){};
    \end{tikzpicture}}}
    =
    \scalebox{.75}{
    \begin{tikzpicture}[Centering,xscale=0.17,yscale=0.33]
        \node[BlackNode](2)at(0.00,-2.00){};
        \node[WhiteNode](0)at(0.00,0.00){};
        \node[BlackNode](1)at(0.00,-1.00){};
        \draw[Edge](1)--(0);
        \draw[Edge](2)--(1);
    \end{tikzpicture}}
    \otimes
    \scalebox{.75}{
    \begin{tikzpicture}[Centering,xscale=0.17,yscale=0.14]
        \node[WhiteNode](0)at(0.00,-4.67){};
        \node[WhiteNode](2)at(2.00,-4.67){};
        \node[WhiteNode](4)at(3.00,-2.33){};
        \node[BlackNode](6)at(5.00,-4.67){};
        \node[WhiteNode](1)at(1.00,-2.33){};
        \node[BlackNode](3)at(3.00,0.00){};
        \node[WhiteNode](5)at(5.00,-2.33){};
        \draw[Edge](0)--(1);
        \draw[Edge](1)--(3);
        \draw[Edge](2)--(1);
        \draw[Edge](4)--(3);
        \draw[Edge](5)--(3);
        \draw[Edge](6)--(5);
    \end{tikzpicture}}
    \otimes
    \scalebox{.75}{
    \begin{tikzpicture}[Centering,xscale=0.13,yscale=0.17]
        \node[WhiteNode](0)at(-1.00,-2.50){};
        \node[BlackNode](1)at(1.00,-2.50){};
        \node[WhiteNode](3)at(3.00,-2.50){};
        \node[WhiteNode](4)at(5.00,-2.50){};
        \node[BlackNode](2)at(2.00,0.00){};
        \draw[Edge](0)--(2);
        \draw[Edge](1)--(2);
        \draw[Edge](3)--(2);
        \draw[Edge](4)--(2);
    \end{tikzpicture}}.
\end{equation}
Intuitively, this product consists in grafting a $\WhiteNode$ onto one forest or a
$\BlackNode$ onto a forest and its right neighbor following the letters of $u$. Observe that
for any $u \in \Bra{\WhiteNode, \BlackNode}^*$ and $\SeriesF \in \T{k + |u|_{\BlackNode}} \K
\AAngle{\SetDuplicativeForests}$, we have
\begin{equation} \label{equ:ht_enumeration_merging_product}
    \EnumerationMap_{\Height}\Par{\Merge_u(\SeriesF)}
    = \VarZ \, \EnumerationMap_{\Height}(\SeriesF).
\end{equation}

Let us finally define the \Def{series of ladders} as the $\SetDuplicativeForests$-series
\begin{equation}
    \SeriesLadders := \sum_{h \geq 0} \Ladder_h
    =
    \epsilon
    +
    \scalebox{.75}{
    \begin{tikzpicture}[Centering,xscale=0.17,yscale=0.3]
        \node[WhiteNode](0)at(0.00,0.00){};
    \end{tikzpicture}}
    +
    \scalebox{.75}{
    \begin{tikzpicture}[Centering,xscale=0.17,yscale=0.3]
        \node[WhiteNode](1)at(0.00,-1.00){};
        \node[WhiteNode](0)at(0.00,0.00){};
    \draw[Edge](1)--(0);
    \end{tikzpicture}}
    +
    \scalebox{.75}{
    \begin{tikzpicture}[Centering,xscale=0.17,yscale=0.3]
        \node[WhiteNode](2)at(0.00,-2.00){};
        \node[WhiteNode](0)at(0.00,0.00){};
        \node[WhiteNode](1)at(0.00,-1.00){};
        \draw[Edge](1)--(0);
        \draw[Edge](2)--(1);
    \end{tikzpicture}}
    +
    \scalebox{.75}{
    \begin{tikzpicture}[Centering,xscale=0.17,yscale=0.3]
        \node[WhiteNode](3)at(0.00,-3.00){};
        \node[WhiteNode](0)at(0.00,0.00){};
        \node[WhiteNode](1)at(0.00,-1.00){};
        \node[WhiteNode](2)at(0.00,-2.00){};
        \draw[Edge](1)--(0);
        \draw[Edge](2)--(1);
        \draw[Edge](3)--(2);
    \end{tikzpicture}}
    + \cdots.
\end{equation}
Observe that
\begin{equation} \label{equ:series_ladders_decomposition}
    \SeriesLadders = \epsilon + \TT{\WhiteNode}{}(\SeriesLadders).
\end{equation}

\subsubsection{Number of elements} \label{subsubsec:number_elements}
Let $\SeriesGreater$ be the $(1, 1)$-operation on $\K \AAngle{\SetDuplicativeForests}$
satisfying, for any $\ForestF \in \SetDuplicativeForests$,
\begin{equation}
    \SeriesGreater(\ForestF)
    = \sum_{\ForestF' \in \SetDuplicativeForests(\ForestF)} \ForestF'.
\end{equation}
By definition, $\SeriesGreater(\ForestF)$ is the characteristic series of
$\SetDuplicativeForests(\ForestF)$ and is therefore also the formal sum of all duplicative
forests $\ForestF'$ such that $\ForestF \LeqDuplicative \ForestF'$. For instance,
\begin{equation}
    \SeriesGreater \Par{
    \scalebox{.75}{
    \begin{tikzpicture}[Centering,xscale=0.19,yscale=0.18]
        \node[WhiteNode](0)at(0.00,-4.00){};
        \node[BlackNode](3)at(2.00,-6.00){};
        \node[BlackNode](7)at(4.00,-6.00){};
        \node[BlackNode](1)at(1.00,-2.00){};
        \node[WhiteNode](2)at(2.00,-4.00){};
        \node[WhiteNode](5)at(4.00,-2.00){};
        \node[BlackNode](6)at(4.00,-4.00){};
        \draw[Edge](0)--(1);
        \draw[Edge](2)--(1);
        \draw[Edge](3)--(2);
        \draw[Edge](6)--(5);
        \draw[Edge](7)--(6);
    \end{tikzpicture}}}
    =
    \scalebox{.75}{
    \begin{tikzpicture}[Centering,xscale=0.2,yscale=0.18]
        \node[WhiteNode](0)at(0.00,-4.00){};
        \node[BlackNode](3)at(2.00,-6.00){};
        \node[BlackNode](7)at(4.00,-6.00){};
        \node[BlackNode](1)at(1.00,-2.00){};
        \node[WhiteNode](2)at(2.00,-4.00){};
        \node[WhiteNode](5)at(4.00,-2.00){};
        \node[BlackNode](6)at(4.00,-4.00){};
        \draw[Edge](0)--(1);
        \draw[Edge](2)--(1);
        \draw[Edge](3)--(2);
        \draw[Edge](6)--(5);
        \draw[Edge](7)--(6);
    \end{tikzpicture}}
    +
    \scalebox{.75}{
    \begin{tikzpicture}[Centering,xscale=0.2,yscale=0.18]
        \node[BlackNode](0)at(0.00,-4.00){};
        \node[BlackNode](3)at(2.00,-6.00){};
        \node[BlackNode](7)at(4.00,-6.00){};
        \node[BlackNode](1)at(1.00,-2.00){};
        \node[WhiteNode](2)at(2.00,-4.00){};
        \node[WhiteNode](5)at(4.00,-2.00){};
        \node[BlackNode](6)at(4.00,-4.00){};
        \draw[Edge](0)--(1);
        \draw[Edge](2)--(1);
        \draw[Edge](3)--(2);
        \draw[Edge](6)--(5);
        \draw[Edge](7)--(6);
    \end{tikzpicture}}
    +
    \scalebox{.75}{
    \begin{tikzpicture}[Centering,xscale=0.2,yscale=0.18]
        \node[WhiteNode](0)at(0.00,-4.00){};
        \node[BlackNode](3)at(2.75,-6.00){};
        \node[BlackNode](3')at(1.25,-6.00){};
        \node[BlackNode](7)at(5.00,-6.00){};
        \node[BlackNode](1)at(1.00,-2.00){};
        \node[BlackNode](2)at(2.00,-4.00){};
        \node[WhiteNode](5)at(5.00,-2.00){};
        \node[BlackNode](6)at(5.00,-4.00){};
        \draw[Edge](0)--(1);
        \draw[Edge](2)--(1);
        \draw[Edge](3)--(2);
        \draw[Edge](3')--(2);
        \draw[Edge](6)--(5);
        \draw[Edge](7)--(6);
    \end{tikzpicture}}
    +
    \scalebox{.75}{
    \begin{tikzpicture}[Centering,xscale=0.2,yscale=0.18]
        \node[WhiteNode](0)at(0.00,-4.00){};
        \node[BlackNode](3)at(2.00,-6.00){};
        \node[BlackNode](7)at(5.75,-6.00){};
        \node[BlackNode](7')at(4.25,-6.00){};
        \node[BlackNode](1)at(1.00,-2.00){};
        \node[WhiteNode](2)at(2.00,-4.00){};
        \node[BlackNode](5)at(5.00,-2.00){};
        \node[BlackNode](6)at(5.75,-4.00){};
        \node[BlackNode](6')at(4.25,-4.00){};
        \draw[Edge](0)--(1);
        \draw[Edge](2)--(1);
        \draw[Edge](3)--(2);
        \draw[Edge](6)--(5);
        \draw[Edge](6')--(5);
        \draw[Edge](7)--(6);
        \draw[Edge](7')--(6');
    \end{tikzpicture}}
    +
    \scalebox{.75}{
    \begin{tikzpicture}[Centering,xscale=0.2,yscale=0.18]
        \node[BlackNode](0)at(0.00,-4.00){};
        \node[BlackNode](3)at(2.75,-6.00){};
        \node[BlackNode](3')at(1.25,-6.00){};
        \node[BlackNode](7)at(5.00,-6.00){};
        \node[BlackNode](1)at(1.00,-2.00){};
        \node[BlackNode](2)at(2.00,-4.00){};
        \node[WhiteNode](5)at(5.00,-2.00){};
        \node[BlackNode](6)at(5.00,-4.00){};
        \draw[Edge](0)--(1);
        \draw[Edge](2)--(1);
        \draw[Edge](3)--(2);
        \draw[Edge](3')--(2);
        \draw[Edge](6)--(5);
        \draw[Edge](7)--(6);
    \end{tikzpicture}}
    +
    \scalebox{.75}{
    \begin{tikzpicture}[Centering,xscale=0.2,yscale=0.18]
        \node[BlackNode](0)at(0.00,-4.00){};
        \node[BlackNode](3)at(2.00,-6.00){};
        \node[BlackNode](7)at(5.75,-6.00){};
        \node[BlackNode](7')at(4.25,-6.00){};
        \node[BlackNode](1)at(1.00,-2.00){};
        \node[WhiteNode](2)at(2.00,-4.00){};
        \node[BlackNode](5)at(5.00,-2.00){};
        \node[BlackNode](6)at(5.75,-4.00){};
        \node[BlackNode](6')at(4.25,-4.00){};
        \draw[Edge](0)--(1);
        \draw[Edge](2)--(1);
        \draw[Edge](3)--(2);
        \draw[Edge](6)--(5);
        \draw[Edge](6')--(5);
        \draw[Edge](7)--(6);
        \draw[Edge](7')--(6');
    \end{tikzpicture}}
    +
    \scalebox{.75}{
    \begin{tikzpicture}[Centering,xscale=0.2,yscale=0.18]
        \node[WhiteNode](0)at(0.00,-4.00){};
        \node[BlackNode](3)at(2.75,-6.00){};
        \node[BlackNode](3')at(1.25,-6.00){};
        \node[BlackNode](7)at(5.75,-6.00){};
        \node[BlackNode](7')at(4.25,-6.00){};
        \node[BlackNode](1)at(1.00,-2.00){};
        \node[BlackNode](2)at(2.00,-4.00){};
        \node[BlackNode](5)at(5.00,-2.00){};
        \node[BlackNode](6)at(5.75,-4.00){};
        \node[BlackNode](6')at(4.25,-4.00){};
        \draw[Edge](0)--(1);
        \draw[Edge](2)--(1);
        \draw[Edge](3)--(2);
        \draw[Edge](3')--(2);
        \draw[Edge](6)--(5);
        \draw[Edge](6')--(5);
        \draw[Edge](7)--(6);
        \draw[Edge](7')--(6');
    \end{tikzpicture}}
    +
    \scalebox{.75}{
    \begin{tikzpicture}[Centering,xscale=0.2,yscale=0.18]
        \node[BlackNode](0)at(0.00,-4.00){};
        \node[BlackNode](3)at(2.75,-6.00){};
        \node[BlackNode](3')at(1.25,-6.00){};
        \node[BlackNode](7)at(5.75,-6.00){};
        \node[BlackNode](7')at(4.25,-6.00){};
        \node[BlackNode](1)at(1.00,-2.00){};
        \node[BlackNode](2)at(2.00,-4.00){};
        \node[BlackNode](5)at(5.00,-2.00){};
        \node[BlackNode](6)at(5.75,-4.00){};
        \node[BlackNode](6')at(4.25,-4.00){};
        \draw[Edge](0)--(1);
        \draw[Edge](2)--(1);
        \draw[Edge](3)--(2);
        \draw[Edge](3')--(2);
        \draw[Edge](6)--(5);
        \draw[Edge](6')--(5);
        \draw[Edge](7)--(6);
        \draw[Edge](7')--(6');
    \end{tikzpicture}}.
\end{equation}

\begin{Proposition} \label{prop:greatest_series}
    For any $\ell \geq 0$, $\ForestF_1, \dots, \ForestF_\ell \in \SetDuplicativeForests$,
    and $\ForestF \in \SetDuplicativeForests$,
    \begin{subequations}
    \begin{equation} \label{equ:greatest_series_1}
        \SeriesGreater \Par{
        \ForestF_1 \ConcatenateForests \dots \ConcatenateForests \ForestF_\ell} =
        \SeriesGreater \Par{\ForestF_1}
        \TT{\ConcatenateForests}{} \dots \TT{\ConcatenateForests}{}
        \SeriesGreater\Par{\ForestF_\ell},
    \end{equation}
    \begin{equation} \label{equ:greatest_series_2}
        \SeriesGreater \Par{\BlackNode(\ForestF)}
        = \TT{\BlackNode}{}(\SeriesGreater(\ForestF)),
    \end{equation}
    \begin{equation} \label{equ:greatest_series_3}
        \SeriesGreater \Par{\WhiteNode(\ForestF)}
        = \TT{\WhiteNode}{}(\SeriesGreater(\ForestF))
        + \TT{\BlackNode}{}(\SeriesGreater(\ForestF \ConcatenateForests \ForestF)).
    \end{equation}
    \end{subequations}
\end{Proposition}
\begin{proof}
    This follows by structural induction on duplicative forests and is a consequence of
    Lemma~\ref{lem:comparison_duplicative_forests} describing a necessary and sufficient
    condition for the fact that a duplicative forest $\ForestF'$ belong to
    $\SetDuplicativeForests(\ForestF)$ and the definitions of the operations
    $\TT{\ConcatenateForests}{}$, $\TT{\WhiteNode}{}$, and $\TT{\BlackNode}{}$ on
    $\SetDuplicativeForests$-series.
\end{proof}

Observe that $\SeriesGreater(\SeriesLadders)$ is the characteristic series of the set
$\GreaterLadders$ (defined in Section~\ref{subsubsec:lattices_duplicative_forests}).
Moreover, since for any $h \geq 0$, all elements of $\SetDuplicativeForests\Par{\Ladder_h}$
have $h$ as height, we have
\begin{equation}
    \EnumerationMap_\Height(\SeriesGreater(\SeriesLadders))
    = \sum_{h \geq 0} \# \SetDuplicativeForests\Par{\Ladder_h} \VarZ^h,
\end{equation}
so that $\EnumerationMap_\Height(\SeriesGreater(\SeriesLadders))$ is the generating series
of the cardinalities of the lattices $\SetDuplicativeForests\Par{\Ladder_h}$, enumerated
w.r.t.\ $h \geq 0$.

\begin{Theorem} \label{thm:series_elements}
    The series $\SeriesGreater(\SeriesLadders)$ satisfies
    \begin{equation}
        \SeriesGreater(\SeriesLadders)
        = \epsilon + \TT{\WhiteNode}{}(\SeriesGreater(\SeriesLadders))
        + \TT{\BlackNode}{}\Par{\SeriesGreater\Par{
            \TT{\ConcatenateForests}{}\Par{\Delta\Par{\SeriesLadders}}}}.
    \end{equation}
\end{Theorem}
\begin{proof}
    By~\eqref{equ:series_ladders_decomposition} and by
    Relation~\eqref{equ:greatest_series_3} of Proposition~\ref{prop:greatest_series}, we
    have
    \begin{equation}
        \SeriesGreater(\SeriesLadders)
        = \SeriesGreater(\epsilon + \TT{\WhiteNode}{}(\SeriesLadders))
        =
        \epsilon
        + \sum_{\substack{h \geq 0}}
        \Par{\TT{\WhiteNode}{}\Par{\SeriesGreater\Par{\Ladder_h}}
        + \TT{\BlackNode}{}\Par{
            \SeriesGreater\Par{\Ladder_h \ConcatenateForests \Ladder_h}}}
    \end{equation}
    and the relation of the statement follows.
\end{proof}

\begin{Proposition} \label{prop:number_elements}
    The $\Height$-enumeration $\GeneratingSeriesGreaterHt$ of
    $\SeriesGreater(\SeriesLadders)$ satisfies
    \begin{equation}
        \GeneratingSeriesGreaterHt
        = 1 + \VarZ \GeneratingSeriesGreaterHt
        + \VarZ \Par{
            \GeneratingSeriesGreaterHt \HadamardProduct \GeneratingSeriesGreaterHt}.
    \end{equation}
\end{Proposition}
\begin{proof}
    By~\eqref{equ:ht_enumeration_concatenation} and by
    Relation~\eqref{equ:greatest_series_1} of Proposition~\ref{prop:greatest_series}, we
    have
    \begin{equation} \begin{split} \label{equ:number_elements_1}
        \EnumerationMap_{\Height}(\SeriesGreater(\TT{\ConcatenateForests}{}
            (\Delta(\SeriesLadders))))
        & =
        \sum_{h \geq 0} \EnumerationMap_{\Height}\Par{
                \SeriesGreater\Par{\Ladder_h \ConcatenateForests \Ladder_h}}
        \\
        & =
        \sum_{h \geq 0} \EnumerationMap_{\Height}\Par{
            \SeriesGreater\Par{\Ladder_h}
                \TT{\ConcatenateForests}{} \SeriesGreater\Par{\Ladder_h}}
        \\
        & =
        \sum_{h \geq 0} \Par{
            \EnumerationMap_{\Height}\Par{\SeriesGreater\Par{\Ladder_h}}
                \MaxProduct \EnumerationMap_{\Height}\Par{\SeriesGreater\Par{\Ladder_h}}}.
    \end{split} \end{equation}
    Since for any $h \geq 0$, all the duplicative forests appearing in
    $\SeriesGreater\Par{\Ladder_h}$ have $h$ as height, the last member
    of~\eqref{equ:number_elements_1} is equal to $\GeneratingSeriesGreaterHt
    \HadamardProduct \GeneratingSeriesGreaterHt$. Now, by using this identity together
    with~\eqref{equ:ht_enumeration_graft_products} and Theorem~\ref{thm:series_elements}, we
    obtain the stated expression for~$\GeneratingSeriesGreaterHt$.
\end{proof}

By Proposition~\ref{prop:poset_isomorphism}, for any $d \geq 1$, the cardinality of
$\MockingbirdLattice(d)$ is
\begin{math}
    \SequenceGreater(h) := \Angle{\VarZ^h, \GeneratingSeriesGreaterHt)}
\end{math}
where $h = d - 1$. By Proposition~\ref{prop:number_elements}, $\SequenceGreater$ is the
integer sequence satisfying $\SequenceGreater(0) = 1$ and, for any $h \geq 1$,
\begin{equation}
    \SequenceGreater(h) = \SequenceGreater(h - 1) + \SequenceGreater(h - 1)^2.
\end{equation}
The sequence of the cardinalities of $\MockingbirdLattice(d)$, $d \geq 0$, starts by
\begin{equation}
    1, 1, 2, 6, 42, 1806, 3263442, 10650056950806
\end{equation}
and forms Sequence~\OEIS{A007018} of~\cite{Slo}. It is known that the sets of the planar
rooted trees $\TreeS$ made of nodes of arity $0$, $1$, or $2$ and such that any path
connecting the root of $\TreeS$ to any leaf of $\TreeS$ consists in $d \geq 1$ nodes are
enumerated by the previous sequence. Let $\phi$ be the map sending any such planar rooted
tree $\TreeS$ to the duplicative tree obtained by replacing each node of $\TreeS$ having a
single child (resp.\ two children) by a white (resp.\ black) node and by removing all nodes
which have no children. It is easy to see that $\phi$ is a bijection between the set of the
planar rooted trees described above and $\SetDuplicativeForests\Par{\Ladder_h}$ where~$h = d
- 1$.

\subsubsection{Number of covering pairs} \label{subsubsec:number_coverings}
Let $\SeriesCoverings$ be the $(1, 1)$-operation on $\K \AAngle{\SetDuplicativeForests}$
satisfying, for any $\ForestF \in \SetDuplicativeForests$,
\begin{equation}
    \SeriesCoverings(\ForestF)
    = \sum_{\substack{
        \ForestF' \in \SetDuplicativeForests \\
        \ForestF \RewDuplicative \ForestF'
    }}
    \ForestF'.
\end{equation}
Immediately from the definition of the covering relation $\RewDuplicative$ of
$\LeqDuplicative$, it follows that for any $\ell \geq 0$, $\ForestF_1, \dots, \ForestF_\ell
\in \SetDuplicativeForests$, and $\ForestF \in \SetDuplicativeForests$,
\begin{subequations}
\begin{equation} \label{equ:coverings_1}
    \SeriesCoverings\Par{\ForestF_1 \ConcatenateForests \dots \ConcatenateForests
        \ForestF_\ell}
    = \sum_{i \in [\ell]}
        \ForestF_1 \ConcatenateForests \dots \ConcatenateForests \ForestF_{i - 1}
        \TT{\ConcatenateForests}{}
        \SeriesCoverings\Par{\ForestF_i}
        \TT{\ConcatenateForests}{}
        \ForestF_{i + 1} \ConcatenateForests \dots \ConcatenateForests \ForestF_\ell,
\end{equation}
\begin{equation} \label{equ:coverings_2}
    \SeriesCoverings\Par{\BlackNode(\ForestF)}
    = \TT{\BlackNode}{}\Par{\SeriesCoverings(\ForestF)},
\end{equation}
\begin{equation} \label{equ:coverings_3}
    \SeriesCoverings\Par{\WhiteNode(\ForestF)}
    = \TT{\WhiteNode}{}(\SeriesCoverings(\ForestF))
    + \TT{\BlackNode}{}(\ForestF \ConcatenateForests \ForestF).
\end{equation}
\end{subequations}
Let $\SeriesNbInputs$ be the $(1, 1)$-operation on $\K \AAngle{\SetDuplicativeForests}$
satisfying, for any $\ForestF \in \SetDuplicativeForests$,
\begin{equation}
    \SeriesNbInputs(\ForestF) = \SeriesCoverings(\SeriesGreater(\ForestF)).
\end{equation}
By a straightforward computation, we obtain
\begin{equation}
    \SeriesNbInputs(\ForestF) =
    \sum_{\ForestF' \in \SetDuplicativeForests(\ForestF)}
    \# \Bra{\ForestF'' \in \SetDuplicativeForests(\ForestF) : \ForestF'' \RewDuplicative
    \ForestF'} \, \ForestF',
\end{equation}
so that the coefficient of each $\ForestF' \in \SetDuplicativeForests(\ForestF)$ in
$\SeriesNbInputs(\ForestF)$ is the number of duplicative forests admitting $\ForestF'$ as
covering in $\SetDuplicativeForests(\ForestF)$. For instance (see at the same time
Figure~\ref{fig:example_poset_duplicative_forests}),
\begin{equation}
    \SeriesNbInputs\Par{
    \scalebox{.75}{
    \begin{tikzpicture}[Centering,xscale=0.18,yscale=0.26]
        \node[WhiteNode](1)at(0.00,-2.67){};
        \node[WhiteNode](3)at(2.00,-1.33){};
        \node[WhiteNode](0)at(0.00,-1.33){};
        \draw[Edge](1)--(0);
    \end{tikzpicture}}}
    =
    \scalebox{.75}{
    \begin{tikzpicture}[Centering,xscale=0.18,yscale=0.26]
        \node[WhiteNode](1)at(0.00,-2.67){};
        \node[BlackNode](3)at(2.00,-1.33){};
        \node[WhiteNode](0)at(0.00,-1.33){};
        \draw[Edge](1)--(0);
    \end{tikzpicture}}
    +
    \scalebox{.75}{
    \begin{tikzpicture}[Centering,xscale=0.18,yscale=0.26]
        \node[BlackNode](1)at(0.00,-2.67){};
        \node[WhiteNode](3)at(2.00,-1.33){};
        \node[WhiteNode](0)at(0.00,-1.33){};
        \draw[Edge](1)--(0);
    \end{tikzpicture}}
    +
    2
    \scalebox{.75}{
    \begin{tikzpicture}[Centering,xscale=0.18,yscale=0.26]
        \node[BlackNode](1)at(0.00,-2.67){};
        \node[BlackNode](3)at(2.00,-1.33){};
        \node[WhiteNode](0)at(0.00,-1.33){};
        \draw[Edge](1)--(0);
    \end{tikzpicture}}
    +
    \scalebox{.75}{
    \begin{tikzpicture}[Centering,xscale=0.18,yscale=0.21]
        \node[WhiteNode](0)at(0.00,-3.33){};
        \node[WhiteNode](2)at(2.00,-3.33){};
        \node[WhiteNode](4)at(4.00,-1.67){};
        \node[BlackNode](1)at(1.00,-1.67){};
        \draw[Edge](0)--(1);
        \draw[Edge](2)--(1);
    \end{tikzpicture}}
    +
    2
    \scalebox{.75}{
    \begin{tikzpicture}[Centering,xscale=0.18,yscale=0.21]
        \node[WhiteNode](0)at(0.00,-3.33){};
        \node[WhiteNode](2)at(2.00,-3.33){};
        \node[BlackNode](4)at(4.00,-1.67){};
        \node[BlackNode](1)at(1.00,-1.67){};
        \draw[Edge](0)--(1);
        \draw[Edge](2)--(1);
    \end{tikzpicture}}
    +
    \scalebox{.75}{
    \begin{tikzpicture}[Centering,xscale=0.18,yscale=0.21]
        \node[WhiteNode](0)at(0.00,-3.33){};
        \node[BlackNode](2)at(2.00,-3.33){};
        \node[WhiteNode](4)at(4.00,-1.67){};
        \node[BlackNode](1)at(1.00,-1.67){};
        \draw[Edge](0)--(1);
        \draw[Edge](2)--(1);
    \end{tikzpicture}}
    +
    \scalebox{.75}{
    \begin{tikzpicture}[Centering,xscale=0.18,yscale=0.21]
        \node[BlackNode](0)at(0.00,-3.33){};
        \node[WhiteNode](2)at(2.00,-3.33){};
        \node[WhiteNode](4)at(4.00,-1.67){};
        \node[BlackNode](1)at(1.00,-1.67){};
        \draw[Edge](0)--(1);
        \draw[Edge](2)--(1);
    \end{tikzpicture}}
    +
    2
    \scalebox{.75}{
    \begin{tikzpicture}[Centering,xscale=0.18,yscale=0.21]
        \node[WhiteNode](0)at(0.00,-3.33){};
        \node[BlackNode](2)at(2.00,-3.33){};
        \node[BlackNode](4)at(4.00,-1.67){};
        \node[BlackNode](1)at(1.00,-1.67){};
        \draw[Edge](0)--(1);
        \draw[Edge](2)--(1);
    \end{tikzpicture}}
    +
    2
    \scalebox{.75}{
    \begin{tikzpicture}[Centering,xscale=0.18,yscale=0.21]
        \node[BlackNode](0)at(0.00,-3.33){};
        \node[WhiteNode](2)at(2.00,-3.33){};
        \node[BlackNode](4)at(4.00,-1.67){};
        \node[BlackNode](1)at(1.00,-1.67){};
        \draw[Edge](0)--(1);
        \draw[Edge](2)--(1);
    \end{tikzpicture}}
    +
    3
    \scalebox{.75}{
    \begin{tikzpicture}[Centering,xscale=0.18,yscale=0.21]
        \node[BlackNode](0)at(0.00,-3.33){};
        \node[BlackNode](2)at(2.00,-3.33){};
        \node[WhiteNode](4)at(4.00,-1.67){};
        \node[BlackNode](1)at(1.00,-1.67){};
        \draw[Edge](0)--(1);
        \draw[Edge](2)--(1);
    \end{tikzpicture}}
    +
    4
    \scalebox{.75}{
    \begin{tikzpicture}[Centering,xscale=0.18,yscale=0.21]
        \node[BlackNode](0)at(0.00,-3.33){};
        \node[BlackNode](2)at(2.00,-3.33){};
        \node[BlackNode](4)at(4.00,-1.67){};
        \node[BlackNode](1)at(1.00,-1.67){};
        \draw[Edge](0)--(1);
        \draw[Edge](2)--(1);
    \end{tikzpicture}}.
\end{equation}

\begin{Proposition} \label{prop:number_inputs}
    For any $\ell \geq 0$, $\ForestF_1, \dots, \ForestF_\ell \in \SetDuplicativeForests$,
    and $\ForestF \in \SetDuplicativeForests$,
    \begin{subequations}
    \begin{equation} \label{equ:number_inputs_1}
        \SeriesNbInputs\Par{
            \ForestF_1 \ConcatenateForests \dots \ConcatenateForests \ForestF_\ell}
        = \sum_{i \in [\ell]}
        \SeriesGreater\Par{
            \ForestF_1 \ConcatenateForests \dots \ConcatenateForests \ForestF_{i - 1}}
        \TT{\ConcatenateForests}{}
        \SeriesNbInputs\Par{\ForestF_i}
        \TT{\ConcatenateForests}{}
        \SeriesGreater\Par{
            \ForestF_{i + 1} \ConcatenateForests \dots \ConcatenateForests \ForestF_\ell},
    \end{equation}
    \begin{equation} \label{equ:number_inputs_2}
        \SeriesNbInputs\Par{\BlackNode(\ForestF)}
        = \TT{\BlackNode}{}(\SeriesNbInputs(\ForestF)),
    \end{equation}
    \begin{equation} \label{equ:number_inputs_3}
        \SeriesNbInputs\Par{\WhiteNode(\ForestF)}
        = \TT{\WhiteNode}{}(\SeriesNbInputs(\ForestF))
        + \TT{\BlackNode}{}(\SeriesNbInputs(\ForestF \ConcatenateForests \ForestF))
        + \TT{\BlackNode}{}\Par{\TT{\ConcatenateForests}{}\Par{
            \Delta\Par{\SeriesGreater\Par{\ForestF}}}}.
    \end{equation}
    \end{subequations}
\end{Proposition}
\begin{proof}
    Relations~\eqref{equ:number_inputs_1} and~\eqref{equ:number_inputs_2} are direct
    consequences of Relations~\eqref{equ:coverings_1} and \eqref{equ:coverings_2}, and of
    Proposition~\ref{prop:greatest_series}. By Proposition~\ref{prop:greatest_series} and
    Relations~\eqref{equ:coverings_2} and~\eqref{equ:coverings_3}, we have
    \begin{equation} \begin{split}
        \SeriesNbInputs\Par{\WhiteNode(\ForestF)}
        & =
        \SeriesCoverings\Par{\TT{\WhiteNode}{}(\SeriesGreater(\ForestF))
        + \TT{\BlackNode}{}(\SeriesGreater(\ForestF \ConcatenateForests \ForestF))}
        \\
        & =
        \sum_{\ForestF' \in \SetDuplicativeForests(\ForestF)}
        \SeriesCoverings\Par{\WhiteNode\Par{\ForestF'}}
        +
        \TT{\BlackNode}{}(\SeriesNbInputs(\ForestF \ConcatenateForests \ForestF))
        \\
        & =
        \sum_{\ForestF' \in \SetDuplicativeForests(\ForestF)}
        \Par{\TT{\WhiteNode}{}\Par{\SeriesCoverings\Par{\ForestF'}}
            + \TT{\BlackNode}{}\Par{\ForestF' \ConcatenateForests \ForestF'}}
        +
        \TT{\BlackNode}{}(\SeriesNbInputs(\ForestF \ConcatenateForests \ForestF)).
    \end{split} \end{equation}
    This establishes~\eqref{equ:number_inputs_3}.
\end{proof}

Observe that
\begin{equation}
    \Support(\SeriesNbInputs(\SeriesLadders))
    = \GreaterLadders \setminus \Bra{\Ladder_h : h \geq 0}
\end{equation}
and that the coefficient of each duplicative forest $\ForestF$ of this set is the number of
duplicative forests of this same set covered by $\ForestF$. Moreover, since for any $h \geq
0$, all elements of $\SetDuplicativeForests\Par{\Ladder_h}$ have $h$ as height,
\begin{equation}
    \EnumerationMap_{\Height}(\SeriesNbInputs(\SeriesLadders))
    = \sum_{h \geq 0} \sum_{\ForestF \in \SetDuplicativeForests\Par{\Ladder_h}}
    \# \Bra{\ForestF' \in \SetDuplicativeForests\Par{\Ladder_h} :
        \ForestF' \RewDuplicative \ForestF} \, \VarZ^h,
\end{equation}
so that $\EnumerationMap_{\Height}(\SeriesNbInputs(\SeriesLadders))$ is the generating
series of the number of edges of the Hasse diagrams of the lattices
$\SetDuplicativeForests\Par{\Ladder_h}$, enumerated w.r.t.\ $h \geq 0$.

\begin{Theorem} \label{thm:series_number_coverings}
    The series $\SeriesNbInputs(\SeriesLadders)$ satisfies
    \begin{equation}
        \SeriesNbInputs(\SeriesLadders)
        = \TT{\WhiteNode}{}(\SeriesNbInputs(\SeriesLadders))
        + \TT{\BlackNode}{}\Par{\SeriesNbInputs\Par{\TT{\ConcatenateForests}{}\Par{
            \Delta\Par{\SeriesLadders}}}}
        + \TT{\BlackNode}{}\Par{\TT{\ConcatenateForests}{}\Par{\Delta\Par{
            \SeriesGreater\Par{\SeriesLadders}}}}.
    \end{equation}
\end{Theorem}
\begin{proof}
    By~\eqref{equ:series_ladders_decomposition} and by Relation~\eqref{equ:number_inputs_3}
    of Proposition~\ref{prop:number_inputs}, we have
    \begin{equation} \begin{split}
        \SeriesNbInputs(\SeriesLadders)
        & =
        \SeriesNbInputs\Par{\epsilon + \TT{\WhiteNode}{}(\SeriesLadders)}
        \\
        & =
        0 + \sum_{h \geq 0} \SeriesNbInputs\Par{\WhiteNode\Par{\Ladder_h}}
        \\
        & =
        \sum_{h \geq 0} \Par{
        \TT{\WhiteNode}{}\Par{\SeriesNbInputs\Par{\Ladder_h}}
        + \TT{\BlackNode}{}\Par{\SeriesNbInputs\Par{
            \Ladder_h \ConcatenateForests \Ladder_h}}
        + \TT{\BlackNode}{}\Par{\TT{\ConcatenateForests}{}\Par{\Delta\Par{
            \SeriesGreater\Par{\Ladder_h}}}}}
    \end{split} \end{equation}
    and the relation of the statement follows.
\end{proof}

\begin{Proposition} \label{prop:number_coverings}
    The $\Height$-enumeration $\GeneratingSeriesNbInputsHt$ of
    $\SeriesNbInputs(\SeriesLadders)$ satisfies
    \begin{equation}
        \GeneratingSeriesNbInputsHt
        = \VarZ \GeneratingSeriesNbInputsHt
        + \VarZ \GeneratingSeriesGreaterHt
        + 2 \VarZ
        \Par{\GeneratingSeriesNbInputsHt \HadamardProduct \GeneratingSeriesGreaterHt}.
    \end{equation}
\end{Proposition}
\begin{proof}
    Observe first that we have
    \begin{equation} \begin{split} \label{equ:number_coverings_1}
        \EnumerationMap_{\Height}\Par{\TT{\ConcatenateForests}{}
            (\Delta(\SeriesGreater(\SeriesLadders)))}
        & =
        \sum_{h \geq 0}
        \sum_{\ForestF \in \SetDuplicativeForests\Par{\Ladder_h}}
        \EnumerationMap_{\Height}(\ForestF \ConcatenateForests \ForestF)
        \\
        & =
        \sum_{h \geq 0}
        \sum_{\ForestF \in \SetDuplicativeForests\Par{\Ladder_h}}
        \EnumerationMap_{\Height}(\ForestF)
        \\
        & =
        \GeneratingSeriesGreaterHt.
    \end{split} \end{equation}
    Moreover, by~\eqref{equ:ht_enumeration_concatenation} and by
    Relation~\eqref{equ:number_inputs_1} of Proposition~\ref{prop:number_inputs}, we have
    \begin{equation} \begin{split} \label{equ:number_coverings_2}
        \EnumerationMap_{\Height}\Par{\SeriesNbInputs\Par{\TT{\ConcatenateForests}{}\Par{
            \Delta\Par{\SeriesLadders}}}}
        & =
        \sum_{h \geq 0}
        \EnumerationMap_{\Height}\Par{\SeriesNbInputs\Par{
            \Ladder_h \ConcatenateForests \Ladder_h}}
        \\
        & =
        \sum_{h \geq 0}
        \EnumerationMap_{\Height}\Par{
            \SeriesGreater\Par{\Ladder_h}
                \TT{\ConcatenateForests}{} \SeriesNbInputs\Par{\Ladder_h}
            +
            \SeriesNbInputs\Par{\Ladder_h}
                \TT{\ConcatenateForests}{} \SeriesGreater\Par{\Ladder_h}}
        \\
        & =
        \sum_{h \geq 0} \Par{
        \EnumerationMap_{\Height}\Par{\SeriesGreater\Par{\Ladder_h}}
            \MaxProduct \EnumerationMap_{\Height}\Par{\SeriesNbInputs\Par{\Ladder_h}}
        +
        \EnumerationMap_{\Height}\Par{\SeriesNbInputs\Par{\Ladder_h}}
            \MaxProduct \EnumerationMap_{\Height}\Par{\SeriesGreater\Par{\Ladder_h}}}.
    \end{split} \end{equation}
    Since for any $h \geq 0$, all duplicative forests appearing in
    $\SeriesGreater\Par{\Ladder_h}$ and in $\SeriesNbInputs\Par{\Ladder_h}$ have $h$ as
    height, the last member of~\eqref{equ:number_coverings_2} is equal to
    \begin{math}
        \GeneratingSeriesGreaterHt \HadamardProduct \GeneratingSeriesNbInputsHt
        + \GeneratingSeriesNbInputsHt \HadamardProduct \GeneratingSeriesGreaterHt.
    \end{math}
    Now by using this identity together with~\eqref{equ:number_coverings_1},
    \eqref{equ:ht_enumeration_graft_products}, and
    Theorem~\ref{thm:series_number_coverings}, we obtain the stated expression
    for~$\GeneratingSeriesNbInputsHt$.
\end{proof}

By Proposition~\ref{prop:poset_isomorphism}, for any $d \geq 1$, the number of edges in the
Hasse diagram of $\MockingbirdLattice(d)$ is
\begin{math}
    \SequenceNbInputs(h) := \Angle{\VarZ^h, \GeneratingSeriesNbInputsHt}
\end{math}
where $h = d - 1$. By Proposition~\ref{prop:number_coverings}, $\SequenceNbInputs$ is the
integer sequence satisfying $\SequenceNbInputs(0) = 0$ and, for any $h \geq 1$,
\begin{equation}
    \SequenceNbInputs(h)
    = \SequenceNbInputs(h - 1) + \SequenceGreater(h - 1)
    + 2 \SequenceNbInputs(h - 1) \, \SequenceGreater(h - 1),
\end{equation}
where $\SequenceGreater$ is the integer sequence defined in
Section~\ref{subsubsec:number_elements}. The sequence of the number of edges of the Hasse
diagram of $\MockingbirdLattice(d)$, $d \geq 0$, starts by
\begin{equation}
    0, 0, 1, 7, 97, 8287, 29942737, 195432804247687.
\end{equation}
This sequence does not appear in~\cite{Slo} for the time being.

\subsubsection{Number of intervals}
Let $\SeriesNbSmaller$ be the $(1, 1)$-operation on $\K \AAngle{\SetDuplicativeForests}$
satisfying, for any $\ForestF \in \SetDuplicativeForests$,
\begin{equation}
    \SeriesNbSmaller(\ForestF) = \SeriesGreater(\SeriesGreater(\ForestF)).
\end{equation}
By a straightforward computation, we obtain
\begin{equation}
    \SeriesNbSmaller(\ForestF) =
    \sum_{\ForestF' \in \SetDuplicativeForests(\ForestF)}
    \# \Han{\ForestF, \ForestF'} \ForestF',
\end{equation}
so that the coefficient of each $\ForestF' \in \SetDuplicativeForests(\ForestF)$ in
$\SeriesNbSmaller(\ForestF)$ is the number of duplicative forests smaller than or equal as
$\ForestF'$ in $\SetDuplicativeForests(\ForestF)$. For instance (see at the same time
Figure~\ref{fig:example_poset_duplicative_forests}),
\begin{multline}
    \SeriesNbSmaller\Par{
    \scalebox{.75}{
    \begin{tikzpicture}[Centering,xscale=0.18,yscale=0.26]
        \node[WhiteNode](1)at(0.00,-2.67){};
        \node[WhiteNode](3)at(2.00,-1.33){};
        \node[WhiteNode](0)at(0.00,-1.33){};
        \draw[Edge](1)--(0);
    \end{tikzpicture}}}
    =
    \scalebox{.75}{
    \begin{tikzpicture}[Centering,xscale=0.18,yscale=0.26]
        \node[WhiteNode](1)at(0.00,-2.67){};
        \node[WhiteNode](3)at(2.00,-1.33){};
        \node[WhiteNode](0)at(0.00,-1.33){};
        \draw[Edge](1)--(0);
    \end{tikzpicture}}
    +
    2
    \scalebox{.75}{
    \begin{tikzpicture}[Centering,xscale=0.18,yscale=0.26]
        \node[WhiteNode](1)at(0.00,-2.67){};
        \node[BlackNode](3)at(2.00,-1.33){};
        \node[WhiteNode](0)at(0.00,-1.33){};
        \draw[Edge](1)--(0);
    \end{tikzpicture}}
    +
    2
    \scalebox{.75}{
    \begin{tikzpicture}[Centering,xscale=0.18,yscale=0.26]
        \node[BlackNode](1)at(0.00,-2.67){};
        \node[WhiteNode](3)at(2.00,-1.33){};
        \node[WhiteNode](0)at(0.00,-1.33){};
        \draw[Edge](1)--(0);
    \end{tikzpicture}}
    +
    4
    \scalebox{.75}{
    \begin{tikzpicture}[Centering,xscale=0.18,yscale=0.26]
        \node[BlackNode](1)at(0.00,-2.67){};
        \node[BlackNode](3)at(2.00,-1.33){};
        \node[WhiteNode](0)at(0.00,-1.33){};
        \draw[Edge](1)--(0);
    \end{tikzpicture}}
    +
    2
    \scalebox{.75}{
    \begin{tikzpicture}[Centering,xscale=0.18,yscale=0.21]
        \node[WhiteNode](0)at(0.00,-3.33){};
        \node[WhiteNode](2)at(2.00,-3.33){};
        \node[WhiteNode](4)at(4.00,-1.67){};
        \node[BlackNode](1)at(1.00,-1.67){};
        \draw[Edge](0)--(1);
        \draw[Edge](2)--(1);
    \end{tikzpicture}}
    +
    4
    \scalebox{.75}{
    \begin{tikzpicture}[Centering,xscale=0.18,yscale=0.21]
        \node[WhiteNode](0)at(0.00,-3.33){};
        \node[WhiteNode](2)at(2.00,-3.33){};
        \node[BlackNode](4)at(4.00,-1.67){};
        \node[BlackNode](1)at(1.00,-1.67){};
        \draw[Edge](0)--(1);
        \draw[Edge](2)--(1);
    \end{tikzpicture}}
    +
    3
    \scalebox{.75}{
    \begin{tikzpicture}[Centering,xscale=0.18,yscale=0.21]
        \node[WhiteNode](0)at(0.00,-3.33){};
        \node[BlackNode](2)at(2.00,-3.33){};
        \node[WhiteNode](4)at(4.00,-1.67){};
        \node[BlackNode](1)at(1.00,-1.67){};
        \draw[Edge](0)--(1);
        \draw[Edge](2)--(1);
    \end{tikzpicture}}
    +
    3
    \scalebox{.75}{
    \begin{tikzpicture}[Centering,xscale=0.18,yscale=0.21]
        \node[BlackNode](0)at(0.00,-3.33){};
        \node[WhiteNode](2)at(2.00,-3.33){};
        \node[WhiteNode](4)at(4.00,-1.67){};
        \node[BlackNode](1)at(1.00,-1.67){};
        \draw[Edge](0)--(1);
        \draw[Edge](2)--(1);
    \end{tikzpicture}}
    \\
    +
    6
    \scalebox{.75}{
    \begin{tikzpicture}[Centering,xscale=0.18,yscale=0.21]
        \node[WhiteNode](0)at(0.00,-3.33){};
        \node[BlackNode](2)at(2.00,-3.33){};
        \node[BlackNode](4)at(4.00,-1.67){};
        \node[BlackNode](1)at(1.00,-1.67){};
        \draw[Edge](0)--(1);
        \draw[Edge](2)--(1);
    \end{tikzpicture}}
    +
    6
    \scalebox{.75}{
    \begin{tikzpicture}[Centering,xscale=0.18,yscale=0.21]
        \node[BlackNode](0)at(0.00,-3.33){};
        \node[WhiteNode](2)at(2.00,-3.33){};
        \node[BlackNode](4)at(4.00,-1.67){};
        \node[BlackNode](1)at(1.00,-1.67){};
        \draw[Edge](0)--(1);
        \draw[Edge](2)--(1);
    \end{tikzpicture}}
    +
    6
    \scalebox{.75}{
    \begin{tikzpicture}[Centering,xscale=0.18,yscale=0.21]
        \node[BlackNode](0)at(0.00,-3.33){};
        \node[BlackNode](2)at(2.00,-3.33){};
        \node[WhiteNode](4)at(4.00,-1.67){};
        \node[BlackNode](1)at(1.00,-1.67){};
        \draw[Edge](0)--(1);
        \draw[Edge](2)--(1);
    \end{tikzpicture}}
    +
    12
    \scalebox{.75}{
    \begin{tikzpicture}[Centering,xscale=0.18,yscale=0.21]
        \node[BlackNode](0)at(0.00,-3.33){};
        \node[BlackNode](2)at(2.00,-3.33){};
        \node[BlackNode](4)at(4.00,-1.67){};
        \node[BlackNode](1)at(1.00,-1.67){};
        \draw[Edge](0)--(1);
        \draw[Edge](2)--(1);
    \end{tikzpicture}}.
\end{multline}

Contrary to what we have undertaken in Sections~\ref{subsubsec:number_elements}
and~\ref{subsubsec:number_coverings} where we have provided direct recursive expressions to
compute $\SeriesGreater(\ForestF)$ and $\SeriesNbInputs(\ForestF)$, we fail to provide
similar expressions for $\SeriesNbSmaller(\ForestF)$. The trick here consists in considering
first a slightly different series depending on a parameter $k \geq 1$ which can be seen as a
catalytic parameter. With this in mind, let for any $k \geq 1$, $\SeriesMeetDecomposition_k$
be the $(1, k)$-operation on $\K \AAngle{ \SetDuplicativeForests}$ satisfying, for any
$\ForestF \in \SetDuplicativeForests$,
\begin{equation}
    \SeriesMeetDecomposition_k(\ForestF) =
    \sum_{\substack{
        \ForestG_1, \dots, \ForestG_k \in \SetDuplicativeForests(\ForestF) \\
        \ForestG_1 \Meet \dots \Meet \ForestG_k = \ForestF
    }}
    \ForestG_1 \otimes \dots \otimes \ForestG_k,
\end{equation}
where $\Meet$ is the meet operation of the duplicative forest lattices introduced in
Section~\ref{subsubsec:lattices_duplicative_forests}. We call
$\SeriesMeetDecomposition_k(\ForestF)$ the \Def{$k$-meet decomposition} of $\ForestF$.
Observe that $\SeriesMeetDecomposition_1$ is the identity map.

\begin{Lemma} \label{lem:meet_decompositions_nb_smaller}
    For any $k \geq 1$ and $\ForestF \in \SetDuplicativeForests$,
    \begin{equation}
        \SeriesMeetDecomposition_k \Par{\SeriesNbSmaller(\ForestF)}
        =
        \sum_{\substack{
            \ForestF' \in \SetDuplicativeForests(\ForestF) \\
            \ForestG_1, \dots, \ForestG_k \in \SetDuplicativeForests\Par{\ForestF'}
        }}
        \ForestG_1 \otimes \dots \otimes \ForestG_k.
    \end{equation}
\end{Lemma}
\begin{proof}
    We have
    \begin{equation} \begin{split}
        \SeriesMeetDecomposition_k \Par{\SeriesNbSmaller(\ForestF)}
        & =
        \sum_{\substack{
            \ForestF' \in \SetDuplicativeForests(\ForestF) \\
            \ForestF'' \in \SetDuplicativeForests \Par{\ForestF'}
        }}
        \sum_{\substack{
            \ForestG_1, \dots, \ForestG_k \in \SetDuplicativeForests\Par{\ForestF} \\
            \ForestG_1 \Meet \dots \Meet \ForestG_k = \ForestF''
        }}
        \ForestG_1 \otimes \dots \otimes \ForestG_k
        \\
        & =
        \sum_{\ForestF' \in \SetDuplicativeForests(\ForestF)}
        \sum_{\substack{
            \ForestG_1, \dots, \ForestG_k \in \SetDuplicativeForests\Par{\ForestF} \\
            \ForestG_1 \Meet \dots \Meet \ForestG_k = \ForestF'
        }}
        \# \Han{\ForestF, \ForestF'}
        \ForestG_1 \otimes \dots \otimes \ForestG_k
        \\
        & =
        \sum_{\ForestG_1, \dots, \ForestG_k \in \SetDuplicativeForests(\ForestF)}
        \# \Han{\ForestF, \ForestG_1 \Meet \dots \Meet \ForestG_k}
        \ForestG_1 \otimes \dots \otimes \ForestG_k
        \\
        & =
        \sum_{\ForestG_1, \dots, \ForestG_k \in \SetDuplicativeForests(\ForestF)}
        \sum_{\substack{
            \ForestF' \in \SetDuplicativeForests(\ForestF) \\
            \ForestF' \LeqDuplicative \ForestG_1, \dots,
            \ForestF' \LeqDuplicative \ForestG_k
        }}
        \ForestG_1 \otimes \dots \otimes \ForestG_k,
    \end{split} \end{equation}
    showing the stated identity.
\end{proof}

\begin{Proposition} \label{prop:meet_decompositions_nb_smaller_recursive}
    For any $k \geq 1$, $\ell \geq 0$, $\ForestF_1, \dots, \ForestF_\ell \in
    \SetDuplicativeForests$, and $\ForestF \in \SetDuplicativeForests$,
    \begin{subequations}
    \begin{equation} \label{equ:meet_decompositions_nb_smaller_recursive_1}
        \SeriesMeetDecomposition_k\Par{\SeriesNbSmaller\Par{
            \ForestF_1 \ConcatenateForests \dots \ConcatenateForests \ForestF_\ell}}
        =
        \SeriesMeetDecomposition_k\Par{\SeriesNbSmaller\Par{\ForestF_1}}
        \TT{\ConcatenateForests}{k} \dots \TT{\ConcatenateForests}{k}
        \SeriesMeetDecomposition_k\Par{\SeriesNbSmaller\Par{\ForestF_\ell}},
    \end{equation}
    \begin{equation} \label{equ:meet_decompositions_nb_smaller_recursive_2}
        \SeriesMeetDecomposition_k\Par{\SeriesNbSmaller\Par{\BlackNode(\ForestF)}}
        = \TT{\BlackNode}{k} \Par{
            \SeriesMeetDecomposition_k\Par{\SeriesNbSmaller\Par{\ForestF}}},
    \end{equation}
    \begin{equation} \label{equ:meet_decompositions_nb_smaller_recursive_3}
        \SeriesMeetDecomposition_k\Par{\SeriesNbSmaller\Par{\WhiteNode(\ForestF)}}
        = \sum_{u \in \Bra{\WhiteNode, \BlackNode}^k}
        \Merge_u
        \Par{\SeriesMeetDecomposition_{k + |u|_{\BlackNode}}
            \Par{\SeriesNbSmaller(\ForestF)}}
        +
        \TT{\BlackNode}{k}\Par{
            \SeriesMeetDecomposition_k\Par{
                \SeriesNbSmaller(\ForestF \ConcatenateForests \ForestF)}}.
    \end{equation}
    \end{subequations}
\end{Proposition}
\begin{proof}
    Relations~\eqref{equ:meet_decompositions_nb_smaller_recursive_1}
    and~\eqref{equ:meet_decompositions_nb_smaller_recursive_2} are direct consequences
    of Lemmas~\ref{lem:comparison_duplicative_forests}
    and~\ref{lem:meet_decompositions_nb_smaller}.
    By Lemmas~\ref{lem:comparison_duplicative_forests}
    and~\ref{lem:meet_decompositions_nb_smaller},
    \begin{equation} \begin{split}
        \SeriesMeetDecomposition_k \Par{\SeriesNbSmaller \Par{\WhiteNode(\ForestF)}}
        & =
        \sum_{\substack{
            \ForestF' \in \SetDuplicativeForests\Par{\WhiteNode(\ForestF)} \\
            \ForestG_1, \dots, \ForestG_k \in \SetDuplicativeForests\Par{\ForestF'}
        }}
        \ForestG_1 \otimes \dots \otimes \ForestG_k
        \\
        & =
        \sum_{\substack{
            \WhiteNode\Par{\ForestF'}
                \in \SetDuplicativeForests\Par{\WhiteNode(\ForestF)} \\
            \ForestG_1, \dots, \ForestG_k
                \in \SetDuplicativeForests\Par{\WhiteNode\Par{\ForestF'}}
        }}
        \ForestG_1 \otimes \dots \otimes \ForestG_k
        +
        \sum_{\substack{
            \BlackNode\Par{\ForestF'}
                \in \SetDuplicativeForests\Par{\WhiteNode(\ForestF)} \\
            \ForestG_1, \dots, \ForestG_k
                \in \SetDuplicativeForests\Par{\BlackNode\Par{\ForestF'}}
        }}
        \ForestG_1 \otimes \dots \otimes \ForestG_k
        \\
        & =
        \sum_{\substack{
            \ForestF' \in \SetDuplicativeForests\Par{\ForestF} \\
            \ForestG_1, \dots, \ForestG_k
                \in \SetDuplicativeForests\Par{\WhiteNode\Par{\ForestF'}}
        }}
        \ForestG_1 \otimes \dots \otimes \ForestG_k
        +
        \sum_{\substack{
            \ForestF'_1, \ForestF'_2
                \in \SetDuplicativeForests\Par{\ForestF} \\
            \ForestG_1, \dots, \ForestG_k
                \in \SetDuplicativeForests\Par{\BlackNode\Par{\ForestF'_1 \ForestF'_2}}
        }}
        \ForestG_1 \otimes \dots \otimes \ForestG_k
        \\
        & =
        \sum_{\substack{
            u \in \Bra{\WhiteNode, \BlackNode}^k \\
            \ForestF' \in \SetDuplicativeForests\Par{\ForestF} \\
            \ForestG_1, \dots, \ForestG_k
                \in \SetDuplicativeForests\Par{\WhiteNode\Par{\ForestF'}} \\
            \Roots\Par{\ForestG_1 \dots \ForestG_k} = u
        }}
        \ForestG_1 \otimes \dots \otimes \ForestG_k
        +
        \sum_{\substack{
            \ForestF'_1, \ForestF'_2
                \in \SetDuplicativeForests\Par{\ForestF} \\
            \ForestG_1', \dots, \ForestG_k' \in \SetDuplicativeForests\Par{\ForestF'_1} \\
            \ForestG_1'', \dots, \ForestG_k'' \in \SetDuplicativeForests\Par{\ForestF'_2} \\
        }}
        \BlackNode\Par{\ForestG'_1 \ConcatenateForests \ForestG''_1} \otimes \dots \otimes
            \BlackNode\Par{\ForestG'_k \ConcatenateForests \ForestG''_k}
        \\
        & =
        \sum_{u \in \Bra{\WhiteNode, \BlackNode}^k}
        \Merge_u \Par{\SeriesMeetDecomposition_{k + |u|_{\BlackNode}}
            \Par{\SeriesNbSmaller(\ForestF)}}
        +
        \TT{\BlackNode}{k}
            \Par{\SeriesMeetDecomposition_k\Par{
                \SeriesNbSmaller(\ForestF \ConcatenateForests \ForestF)}},
    \end{split} \end{equation}
    where for any duplicative forest $\ForestF$, $\Roots(\ForestF)$ is the word on
    $\Bra{\WhiteNode, \BlackNode}$ containing from left to right the roots of the trees
    forming $\ForestF$. This shows~\eqref{equ:meet_decompositions_nb_smaller_recursive_3}.
\end{proof}

Observe that $\Support\Par{\SeriesNbSmaller(\SeriesLadders)} = \GreaterLadders$ and that the
coefficient of each duplicative forest $\ForestF$ of this set is the number of duplicative
forests of this same set smaller than or equal as $\ForestF$. Moreover, since for any $h
\geq 0$, all elements of $\SetDuplicativeForests\Par{\Ladder_h}$ have $h$ as height,
\begin{equation}
    \EnumerationMap_{\Height}(\SeriesNbSmaller(\SeriesLadders))
    = \sum_{h \geq 0}
    \sum_{\ForestF \in \SetDuplicativeForests\Par{\Ladder_h}}
    \# \Han{\Ladder_h, \ForestF} \, \VarZ^h,
\end{equation}
so that $\EnumerationMap_{\Height}(\SeriesNbSmaller(\SeriesLadders))$ is the generating
series of the number of intervals of the lattices $\SetDuplicativeForests\Par{\Ladder_h}$,
enumerated w.r.t.\ $h \geq 0$.

\begin{Theorem} \label{thm:series_nb_smaller}
    The series $\SeriesNbSmaller(\SeriesLadders)$ satisfies
    $\SeriesNbSmaller(\SeriesLadders) =
    \SeriesMeetDecomposition_1(\SeriesNbSmaller(\SeriesLadders))$ where, for any $k \geq 1$,
    the series $\SeriesMeetDecomposition_k(\SeriesNbSmaller(\SeriesLadders))$ satisfies
    \begin{equation}
        \SeriesMeetDecomposition_k(\SeriesNbSmaller(\SeriesLadders))
        = \underbrace{\epsilon \otimes \dots \otimes \epsilon}_
            {k \mbox{ \footnotesize terms}}
        + \sum_{u \in \Bra{\WhiteNode, \BlackNode}^k}
            \Merge_u\Par{\SeriesMeetDecomposition_{k + |u|_{\BlackNode}}
                (\SeriesNbSmaller(\SeriesLadders))}
        + \TT{\BlackNode}{k} \Par{
            \SeriesMeetDecomposition_k\Par{\SeriesNbSmaller\Par{\TT{\ConcatenateForests}{}
            \Par{\Delta\Par{\SeriesLadders}}}}}.
    \end{equation}
\end{Theorem}
\begin{proof}
    First, the relation $\SeriesNbSmaller(\SeriesLadders) =
    \SeriesMeetDecomposition_1(\SeriesNbSmaller(\SeriesLadders))$ holds since
    $\SeriesMeetDecomposition_1$ is the identity map.
    By~\eqref{equ:series_ladders_decomposition} and by
    Proposition~\ref{prop:meet_decompositions_nb_smaller_recursive}, we have
    \begin{equation} \begin{split}
        \SeriesMeetDecomposition_k(\SeriesNbSmaller(\SeriesLadders))
        & =
        \SeriesMeetDecomposition_k\Par{
            \SeriesNbSmaller\Par{\epsilon + \TT{\WhiteNode}{}(\SeriesLadders)}}
        \\
        & =
        \underbrace{\epsilon \otimes \dots \otimes \epsilon}_{k \mbox{ \footnotesize terms}}
        + \sum_{h \geq 0}
        \Par{
        \sum_{u \in \Bra{\WhiteNode, \BlackNode}^k}
        \Merge_u\Par{\SeriesMeetDecomposition_{k + |u|_{\BlackNode}}\Par{
            \SeriesNbSmaller\Par{\Ladder_h}}}
        +
        \TT{\BlackNode}{k}\Par{\SeriesMeetDecomposition_k\Par{\SeriesNbSmaller\Par{
            \Ladder_h \ConcatenateForests \Ladder_h}}}}
    \end{split} \end{equation}
    and the relation of the statement follows.
\end{proof}

\begin{Proposition} \label{prop:number_intervals}
    The $\Height$-enumeration $\GeneratingSeriesNbSmallerHt$ of
    $\SeriesNbSmaller(\SeriesLadders)$ satisfies $\GeneratingSeriesNbSmallerHt =
    \GeneratingSeriesNbSmallerHt^{(1)}$ where, for any $k \geq 1$,
    $\GeneratingSeriesNbSmallerHt^{(k)}$ is the $\Height$-enumeration of
    $\SeriesMeetDecomposition_k(\SeriesNbSmaller(\SeriesLadders))$ which satisfies
    \begin{equation}
        \GeneratingSeriesNbSmallerHt^{(k)}
        = 1 + \VarZ
        \Par{\GeneratingSeriesNbSmallerHt^{(k)}
        \HadamardProduct \GeneratingSeriesNbSmallerHt^{(k)}}
        + \VarZ \sum_{i \in \HanL{k}} \binom{k}{i} \GeneratingSeriesNbSmallerHt^{(k + i)}.
    \end{equation}
\end{Proposition}
\begin{proof}
    By Relation~\eqref{equ:meet_decompositions_nb_smaller_recursive_1} of
    Proposition~\ref{prop:meet_decompositions_nb_smaller_recursive}, we have
    \begin{equation} \begin{split} \label{equ:number_intervals_1}
        \EnumerationMap_{\Height}
        \Par{\SeriesMeetDecomposition_k\Par{\SeriesNbSmaller\Par{\TT{\ConcatenateForests}{}
            \Par{\Delta\Par{\SeriesLadders}}}}}
        & =
        \sum_{h \geq 0} \EnumerationMap_{\Height}\Par{
            \SeriesMeetDecomposition_k\Par{\SeriesNbSmaller\Par{
                \Ladder_h \ConcatenateForests \Ladder_h}}}
        \\
        & =
        \sum_{h \geq 0} \EnumerationMap_{\Height}\Par{
            \SeriesMeetDecomposition_k\Par{\SeriesNbSmaller\Par{\Ladder_h}}
            \TT{\ConcatenateForests}{k}
            \SeriesMeetDecomposition_k\Par{\SeriesNbSmaller\Par{\Ladder_h}}
        }
        \\
        & =
        \sum_{h \geq 0} \Par{\EnumerationMap_{\Height}\Par{
            \SeriesMeetDecomposition_k\Par{\SeriesNbSmaller\Par{\Ladder_h}}}
            \MaxProduct
            \EnumerationMap_{\Height}\Par{
            \SeriesMeetDecomposition_k\Par{\SeriesNbSmaller\Par{\Ladder_h}}}}.
    \end{split} \end{equation}
    Since for any $h \geq 0$, all the duplicative forests of the tensors $\ForestF_1 \otimes
    \dots \otimes \ForestF_k$ appearing in
    $\SeriesMeetDecomposition_k\Par{\SeriesNbSmaller\Par{\Ladder_h}}$ have $h$ as height,
    the last member of~\eqref{equ:number_intervals_1} is equal to
    $\GeneratingSeriesNbSmallerHt^{(k)} \HadamardProduct
    \GeneratingSeriesNbSmallerHt^{(k)}$. Moreover,
    by~\eqref{equ:ht_enumeration_merging_product}, for any $u \in \Bra{\WhiteNode,
    \BlackNode}^k$,
    \begin{equation}
        \EnumerationMap_{\Height}\Par{
            \Merge_u\Par{\SeriesMeetDecomposition_{k + |u|_{\BlackNode}}
                (\SeriesNbSmaller(\SeriesLadders))}
        }
        =
        \VarZ \, \EnumerationMap_{\Height}\Par{
        \SeriesMeetDecomposition_{k + |u|_{\BlackNode}}(\SeriesNbSmaller(\SeriesLadders)}
        =
        \VarZ \GeneratingSeriesNbSmallerHt^{(k + |u|_{\BlackNode})}.
    \end{equation}
    Now, by using these two identities together with Theorem~\ref{thm:series_nb_smaller},
    and by the fact that the number of words of length $k$ on $\Bra{\WhiteNode, \BlackNode}$
    having exactly $i \in \HanL{k}$ occurrences of $\BlackNode$ is $\binom{k}{i}$, we obtain
    the stated expression for $\GeneratingSeriesNbSmallerHt^{(k)}$. Finally, we have
    $\GeneratingSeriesNbSmallerHt = \GeneratingSeriesNbSmallerHt^{(1)}$ since
    $\SeriesNbSmaller(\SeriesLadders) =
    \SeriesMeetDecomposition_1(\SeriesNbSmaller(\SeriesLadders))$.
\end{proof}

By Proposition~\ref{prop:poset_isomorphism}, for any $d \geq 1$, the number of
intervals in $\MockingbirdLattice(d)$ is
\begin{math}
    \SequenceNbSmaller^{(1)}(h) := \Angle{\VarZ^h, \GeneratingSeriesNbSmallerHt}
\end{math}
where $h = d - 1$. By Proposition~\ref{prop:number_intervals}, for any $k \geq 1$,
$\SequenceNbSmaller^{(k)}$ is the integer sequence satisfying $\SequenceNbSmaller^{(k)}(0) =
1$ and, for any $h \geq 1$,
\begin{equation}
    \SequenceNbSmaller^{(k)}(h) =
    \SequenceNbSmaller^{(k)}(h - 1)^2
    + \sum_{i \in \HanL{k}} \binom{k}{i} \SequenceNbSmaller^{(k + i)}(h - 1).
\end{equation}
The sequence of the cardinalities of $\MockingbirdLattice(d)$, $d \geq 0$, starts by
\begin{equation}
    1, 1, 3, 17, 371, 144513, 20932611523, 438176621806663544657.
\end{equation}
This sequence does not appear in~\cite{Slo} for the time being.

\section*{Open questions and future work}
We have studied in this work a CLS having a lot of rich combinatorial properties despite its
simplicity. This can be considered as the prototypical example for this kind of
investigation. We could expect to have also similar combinatorial properties (related to
poset and lattice properties, and quantitative data) for more complex CLS. In addition to
the obvious possible future investigations consisting in studying in a similar fashion some
other CLS, we describe here three open questions raised by this work.

The description of minimal and maximal elements of $\Poset$ uses a notion of pattern
avoidance in terms (see Proposition~\ref{prop:minimal_maximal}). This is a general fact:
when a CLS $(\GeneratingSet, \Rew)$ has the poset property, its minimal (resp.\ maximal)
elements are the terms avoiding some terms deduced from the ones appearing as right-hand
(resp.\ left-hand) members of $\Rew$. Such an enumerative problem has been considered
in~\cite{KP15,Gir20} for the particular case of linear terms. The question here concerns the
general enumeration of terms avoiding a set of other terms wherein multiple occurrences of
a same variable are allowed.

We have shown that the Mockingbird CLS has the poset property and is rooted
(Proposition~\ref{prop:first_graph_properties}), and has the lattice property
(Theorem~\ref{thm:mockingbird_lattices}) by employing some specific reasoning from the
definition of the Mockingbird basic combinator $\M$. A question here concerns the existence
of a general criterion to decide if a CLS has the poset (resp.\ lattice) property and if it
is rooted. All this seems independent from the property of a CLS to be hierarchical because,
among others, the CLS containing a basic combinator $\Combinator{X}$ such that
$\Combinator{X}$ has order $3$ and $\TreeT_{\Combinator{X}} := \VarX_3 \Par{\VarX_1 \VarX_1
\VarX_2}$ has not the poset property because we have $\TreeT \RewContext \TreeT'$ and
$\TreeT' \RewContext \TreeT$ with
\begin{math}
    \TreeT :=
    \Combinator{X}\Combinator{X}\Combinator{X}
    \Par{\Combinator{X}\Combinator{X}\Par{\Combinator{X}\Combinator{X}}}
\end{math}
and
\begin{math}
    \TreeT' :=
    \Combinator{X} \Combinator{X}
    \Par{\Combinator{X}\Combinator{X}}
    \Par{\Combinator{X}\Combinator{X}\Combinator{X}}
\end{math}
while $\Combinator{X}$ is hierarchical.

Finally, in Proposition~\ref{prop:finite_equivalence_classes}, we have shown that being
hierarchical is a sufficient condition for a CLS $\CLS$ to be locally finite. The question
in this context consists in obtaining a necessary and sufficient condition for this last
property. Observe that being nonerasing is necessary because, by assuming that $\CLS$
contains a basic combinator $\Combinator{X}$ such that $\Combinator{X}$ has order $n \geq 2$
and that there is $i \in [n]$ such that $\VarX_i$ does not appear in
$\TreeT_{\Combinator{X}}$, we have
\begin{math}
    \Combinator{X} \VarX_1 \dots \VarX_{i - 1} \TreeS \VarX_{i + 1} \dots \VarX_n
    \RewContext \TreeT_{\Combinator{X}}
\end{math}
for any term $\TreeS$, showing that the connected component of $\TreeT_{\Combinator{X}}$ in
$\RewGraph_\CLS$ is infinite.

\bibliographystyle{alpha}
\bibliography{Bibliography}

\end{document}